\documentclass[final,3p,times]{elsarticle}

\usepackage{amssymb}
\usepackage{amsthm}

\usepackage{mathtools, float}
\usepackage{booktabs} 
\usepackage{multirow}
\usepackage{makecell}
\usepackage{multicol}
\usepackage{tabularx}
\usepackage{xcolor}
\usepackage{caption}
\usepackage{subcaption}

\usepackage[colorlinks=false, citecolor=red, linkcolor=blue, pdfborder={0 0 0}]{hyperref}
\usepackage[T1]{fontenc}
\usepackage{selinput}

\journal{Journal of Sound and Vibration}

\begin{document}

\begin{frontmatter}



\title{Including connecting elements into the Lagrange Multiplier State-Space Substructuring Formulation}


\author[inst1]{R.S.O. Dias
\texorpdfstring{\corref{cor1}}{Lg}}
\cortext[cor1]{Corresponding author.}

\ead{r.dasilva@staff.univpm.it}
\affiliation[inst1]{organization={Department of Industrial Engineering and Mathematical Sciences, Universit\`a Politecnica delle Marche},
            addressline={Via Brecce Bianche}, 
            city={Ancona},
            postcode={60131}, 
            state={Marche},
            country={Italy}}

\author[inst1]{M. Martarelli}
\ead{m.martarelli@staff.univpm.it}

\author[inst2]{P. Chiariotti}
\ead{paolo.chiariotti@polimi.it}

\affiliation[inst2]{organization={Department of Mechanical Engineering, Politecnico di Milano},
            addressline={Via Privata Giuseppe La Masa}, 
            city={Milan},
            postcode={20156}, 
            state={Lombardia},
            country={Italy}}
            

\begin{abstract}
\textcolor{black}{This paper extends the inverse substructuring (IS) approach into the state-space domain and presents a novel state-space substructuring (SSS) technique that embeds the dynamics of connecting elements (CEs) in the Lagrange Multiplier State-Space Substructuring (LM-SSS) formulation via compatibility relaxation. This coupling approach makes it possible to incorporate into LM-SSS connecting elements that are suitable for being characterized by inverse substructuring (e.g. rubber mounts) by simply using information from one of its \textcolor{black}{off} diagonal apparent mass terms. Therefore, the information obtained from an in-situ experimental characterization of the CEs is enough to include them into the coupling formulation. Moreover, LM-SSS with compatibility relaxation makes it possible to couple an unlimited number of components and CEs, requiring only one matrix inversion to compute the coupled state-space model (SSM). Two post-processing procedures to enable the computation of minimal-order coupled models by using this approach are also presented. Numerical and experimental substructuring applications are exploited to prove the validity of the proposed state-space realization of IS and the inclusion of CEs into the LM-SSS formulation via compatibility relaxation. It is found that the IS approach can be accurately applied on state-space models representative of components linked by CEs to identify models representative of the diagonal apparent mass terms of the CEs, provided that the CEs can be accurately characterized by the underlying assumptions of IS. In this way, state-space models representative of experimentally characterized CEs can be found without performing decoupling operations. Hence, these models are not contaminated with spurious states. Furthermore, it was found that the developed coupling approach is reliable, when the dynamics of the CEs can be accurately characterized by IS, thus making it possible to compute reliable coupled models that are not composed by spurious states.}  
\end{abstract}


\begin{highlights}
\color{black}
\item The state-space realization of IS is presented.
\item A novel SSS technique is developed via compatibility relaxation of LM-SSS.
\item An unlimited number of components and connecting elements (CEs) can be coupled jointly.
\item Two post-processing procedures to obtain the coupled SSM minimal form are discussed.
\item Reliable minimal order coupled SSMs can be computed.
\item Numerical and experimental validation of the discussed approaches is provided.
\end{highlights}
\color{black}

\end{frontmatter}



\section{Nomenclature}

For easier understanding of all the parameters used in this article, the nomenclature is given in table \ref{table:Nomenclature}.

\begin{table}[H]
\centering
\caption{Nomenclature}
\begin{tabular}{@{}llll@{}}
\toprule
$A$ & state matrix & $L$ & Boolean localization matrix\\
$B$ & input matrix & $M$ & mass matrix\\
$B_{C}$ & mapping matrix & $u/U^{\ast}$ & input vector\\
$B_{T}$ & state mapping matrix & $V$ & damping matrix\\
$C$ & output matrix & $x/X^{\ast}$ & state vector\\
$D$ & feed-through matrix & $y/Y^{\ast}$ & displacement output vector\\
$g/G^{\ast}$ & connecting forces vector & $z$ & state vector transformed into coupling form\\
$H$ & Accelerance matrix &  $Z^{A}$ & apparent mass matrix\\
$K$ & stiffness matrix & $Z^{\alpha}$ & dynamic stiffness matrix\\
\vspace{1mm} & & &\\
$\gamma$ & Gaussian distributed independent stochastic variable & $\lambda/\Lambda^{\ast}$ & Lagrange multipliers vector\\ 
$\theta$ & Gaussian distributed independent stochastic variable & \\ 
\vspace{1mm} & & &\\
$\bullet_{D}$ & block diagonal matrix &  $\bullet_{k}$              & k-th discrete frequency\\
$\bullet_{F_{A}}$ & fixture $A$ & $\bullet_{M}$ & variable associated to connecting elements\\
$\bullet_{F_{B}}$ & fixture $B$ & $\bullet_{neg}$ & matrix transformed into negative form\\
$\bullet_{F_{A}MF_{B}}$ & structure composed by $F_{A}$ and $F_{B}$ linked by $M$ & $\bullet_{p}$ & p-th substructure \\
$\bullet_{i}$ & i-th output & $\bullet_{S}$ & matrix associated to structures  \\
$\bullet_{j}$ & j-th input & &\\
$\bullet_{\alpha}$ & substructure $\alpha$ & $\bullet_{M\alpha}$ & side of $M$ which is linked to substructure $\alpha$\\
$\bullet_{\alpha M \beta}$ & structure composed by $\alpha$ and $\beta$ linked by $M$ & $\bullet_{M\beta}$  & side of $M$ which is linked to substructure $\beta$ \\
$\bullet_{\beta}$ & substructure $\beta$ & & \\
\vspace{1mm} & & &\\

$\bullet^{-1}$            & inverse of a matrix & $\bullet^{\dag}$  & pseudoinverse of a matrix\\ 
$\bullet^{accel}$             & acceleration &  $\bullet^{J}$ & interface DOF \\ 
$\bullet^{disp}$               & displacement state-space model & $\bullet^{T}$              & transpose of a vector/matrix \\ 
$\bullet^{I}$ & internal DOF & $\bullet^{vel}$ & velocity state-space model\\ 
$\bullet^{inv}$ & matrix of an inverted state-space model & & \\

\vspace{1mm} & & &\\
$\dot{\bullet}$  & first order time derivative & $\ddot{\bullet}$  & second order time derivative\\
$\bar{\bullet}$  & coupling vector/matrix & &\\ \midrule
\noalign{\makecell{$^{\ast}$ Lower case and capital letters correspond to time and frequency dependent variables, respectively.}}\\
\label{table:Nomenclature}
\end{tabular}
\end{table}

\section{Introduction}\label{Introduction}

Dynamic substructuring revealed to be a valuable tool to analyze complex structures. This concept was firstly introduced in the last century and relies on the dynamic characterization of complex structures by separate analyses of its components. This concept has been continuously investigated and several formulations have been proposed by the scientific community, all ranging from modal, frequency and time domain approaches \citep{AEM_2020}.

The focus of this paper will be on time domain methods, specifically on the class of state-space substructuring techniques. Literature provides several state-space substructuring methods \citep{SJ_94},\citep{SJO_20072697},\citep{mg_2013},\citep{MG_2014},\citep{RD_2021} (just to cite some), however none of them discusses the possibility to include CEs into the formulations. Yet, nowadays, the use of CEs, e.g. rubber mounts, is frequent in many mechanical systems and engineering areas (e.g. automotive industry \citep{AEM_2020}). \textcolor{black}{When exploited in the experimental substructuring field, SSS formulations hold the disadvantage of requiring the use of system identification algorithms to compute state-space models representative of the components under study.
Nonetheless, for some applications, the representation of the dynamics of both the connecting elements, and the coupled structure wherein they are included, through state-space models is advantageous.
When dealing, for instance, with connecting elements presenting time-varying performances (e.g. rubber mount submitted to different temperatures over time), the state-space representation has great potentials and there are several publications in the literature treating the state-space identification of components presenting varying performances (see for example \cite{JC_2009},\cite{JC_2014},\cite{FF_2017},\cite{FF_2007}). These approaches enable the computation of linear-parameter varying (LPV) models that describe the variation of the dynamics of substructures presenting varying dynamic behaviour. Furthermore, by using SSS formulations we are capable of computing different coupled models for each variation on the dynamics of the CEs, thus making it possible to perform time-domain analyses with assembled structures presenting time-varying performances. Last but not least, as by using SSS formulations we are working in time-domain, these approaches are more suited to deal with highly damped CEs than the FBS formulations.}


In the experimental substructuring field, the inclusion of CEs into FBS methods is generally performed by using two different approaches, either by considering the CE as another structure to be coupled or by performing coupling by relaxation of the compatibility conditions \citep{AEM_2020}.

The first approach mentioned relies on testing the CEs incorporated in an assembly (for instance, by testing the connecting element with two fixtures attached at its ends). Hence, after the dynamic characterization of the assembly, the fixtures used to test it must be removed by using decoupling, thus making the identification of the CE possible \citep{SWK_2020},\citep{MH_2020}. This operation might be performed by using either primal or dual formulations \citep{DK_20081169}. To perform the decoupling of the fixtures, one may use dynamic information at the interfaces. However, measuring the dynamic response at the interface of the structure might be difficult when testing the component with fixtures attached. One possibility to overcome this difficulty is testing the CE attached to fixtures manufactured on purpose. For instance, M. H\"auessler et al. \citep{MH_2020} proposed to use aluminium crosses specifically designed to behave as rigid bodies within the frequency range of interest. The use of those fixtures made it possible to excite both translational and rotational DOFs of the CE leading to a twelve DOFs characterization. Furthermore, the rigid body behaviour made it possible the use of virtual point transformation (VPT) \citep{MVS_2016},\citep{MV_2013},\citep{MH_2017} to calculate the FRFs at the interfaces from measuring and exciting at internal locations. The whole approach still requires the testing of the assembly (CE with attached fixtures) and the fixtures isolated. 

The second approach embeds the dynamics of the CE into the coupling formulation by relaxing the interface compatibility conditions \citep{AEM_2020},\citep{EB_2014},\citep{SWK_2020}. By following this coupling approach, the dynamics of the CE is included by using just one of its diagonal dynamic stiffness terms. For this reason, the use of this coupling method is equivalent to characterize the CE by using the inverse substructuring (IS) method \citep{SWK_2020}, which was discussed in several publications, \citep{JZ_2004},\citep{JW_2014},\citep{JM_2015},\citep{JM_T_2017} (in some publications also denominated as in-situ characterization). IS assumes that the CE is massless and that there are no cross couplings between DOFs. This means that each DOF on one side of the connecting element is only coupled to one DOF on the other side (being both DOFs associated to the exact same direction). As proven in \citep{JM_2015}, the \textcolor{black}{off} diagonal terms of the dynamic stiffness of a system composed by two substructures connected by a CE are a property of the CE. \textcolor{black}{Hence, to identify the off diagonal terms of the CE, we are not required to dismount it from the structure were it is included and we are not required to know the dynamics of the structures to which the CE is coupled. Moreover, if inverse substructuring assumptions are valid the identified off diagonal dynamic stiffness terms can be directly used to calculate the diagonal ones \citep{JM_2017},\citep{JM_T_2017},\citep{MH_2020}.}  

The first approach presented to couple CEs can be performed by using the SSS methods already well known from literature. \textcolor{black}{However, in an experimental scenario, to include a CE by treating it as a regular component, we are in general demanded to perform decoupling operations to identify a state-space model representative of its dynamics. When decoupling the state-space models representative of the fixtures from the one representative of the CE with fixtures, the obtained state-space model representative of the CE will present $n_{fmf}+2\times n_{f}-2n_{J}$ states (assuming that a minimal-order coupled state-space model is obtained \citep{SJ_94}), where $n_{fmf}$ is the number of states of the state-space model representative of the CE with fixtures attached to its edges, $n_{f}$ is the number of states of the model representative of a single fixture and $n_{J}$ is the number of interface DOFs.}

For state-space models estimated from experimentally acquired data is usually verified that $n>>n_{o}$ (being $n$ and $n_{o}$ the number of states and outputs of a generic state-space model, respectively), hence the decoupling operation might generate an important increase of the number of states. This increment of the number of states will definitely result in a greater computational effort when performing calculations. In fact, it has been well discussed in  \citep{MS_2015} that the increment of the number of states is caused by the double inclusion of the dynamics of the decoupled structure, since the same dynamics was already present on the initial state-space model and was reintroduced to perform decoupling. Consequently, pairs of spurious modes will be present on the obtained state-space model. These modes are easy to identify and eliminate when there is previous knowledge about the substructure intended to be identified. However, in a real case scenario this method is used to characterize CEs from which no a-priori information is known, turning the identification and elimination of the spurious modes, which might be present on the obtained state-space model, hard to be accomplished.

\color{black}

In this way, the identified state-space models representative of the CEs will in general be composed by spurious states, hence leading to the computation of coupled models spoiled by spurious states. The calculation of these coupled state-space models holds two negative consequences. Firstly, a less elegant representation of the coupled structure is obtained. Secondly, it leads to the computation of state-space models composed by a large number of states, with a consequential increase in the computational cost characterizing these models in time-domain analyses. The increment of the number of states will be more significant as the number of state-space models to be decoupled in order to identify the intended CE model increases and as these state-space models are composed by more states. Consequently, the suitability of the computed coupled models to be exploited in some applications, for example in a real-time context, might be compromised.

To overcome, the mentioned negative effects of including experimentally characterized CEs by treating them as regular structures, we aim at developing a novel SSS technique by relaxing the compatibility condition of the LM-SSS formulation (see \cite{BK_2020},\cite{RD_2021},\cite{RD_MSSP_2022}), enabling the inclusion of the dynamics of CEs by using one of its diagonal terms. This formulation would be tailored to couple CEs suitable to be characterized by IS, because its diagonal terms can be identified from the state-space model representative of the structure were they are included. Therefore, there is no need of performing decoupling operations. Thereby, we also aim at presenting the state-space realization of the IS approach, which will provide us with the possibility to identify models composed by a significantly lower number of states without the presence of spurious states. Hence, when including the dynamics of the CEs by using state-space models representative of their diagonal terms, we are not spoiling the coupled model with spurious states. Thus, the computed coupled state-space models will be a more elegant representation of the coupled structure and will in general be composed by a significantly lower number of states, opening the perspectives of exploiting them on new applications. Furthermore, the development of this approach will allow the user to benefit from the practical advantages of exploiting IS to characterize CEs, i.e. we are not required to dismount the CE from the structure where it is included and we are not demanded to know the dynamics of the structures to which the CE is coupled. Consequently, dismounting and mounting operations are avoided and we are not required to perform further experimental tests to characterize the dynamic behaviour of the structures to which the CE is coupled.

\color{black}

\textcolor{black}{The LM-SSS with compatibility relaxation method is theoretically developed in section \ref{LM_SSS_Non_Rigid_Formulation}, while section \ref{Minimal-order coupled state-space models} discusses two post-processing procedures to erase the extra states that are present in the coupled modes obtained by the developed coupling technique. In section \ref{Determining state-space models representative of diagonal apparent mass terms of connecting elements} procedures to determine analytically, numerically and experimentally the state-space models representative of the CEs to be included in the LM-SSS with compatibility relaxation formulation are outlined. All the discussed approaches are validated on both numerical and experimental substructuring scenarios in sections \ref{Numerical Example} and \ref{Experimental Validation}, respectively. Finally, the conclusions of this paper are presented in section \ref{Conclusion}}.

\section{\textcolor{black}{Including CEs in the LM-SSS formulation via compatibility relaxation}}\label{LM_SSS_Non_Rigid_Formulation}

\begin{figure}
\centering
\begin{subfigure}{0.5\textwidth}
  \centering
  \includegraphics[width=0.7\linewidth]{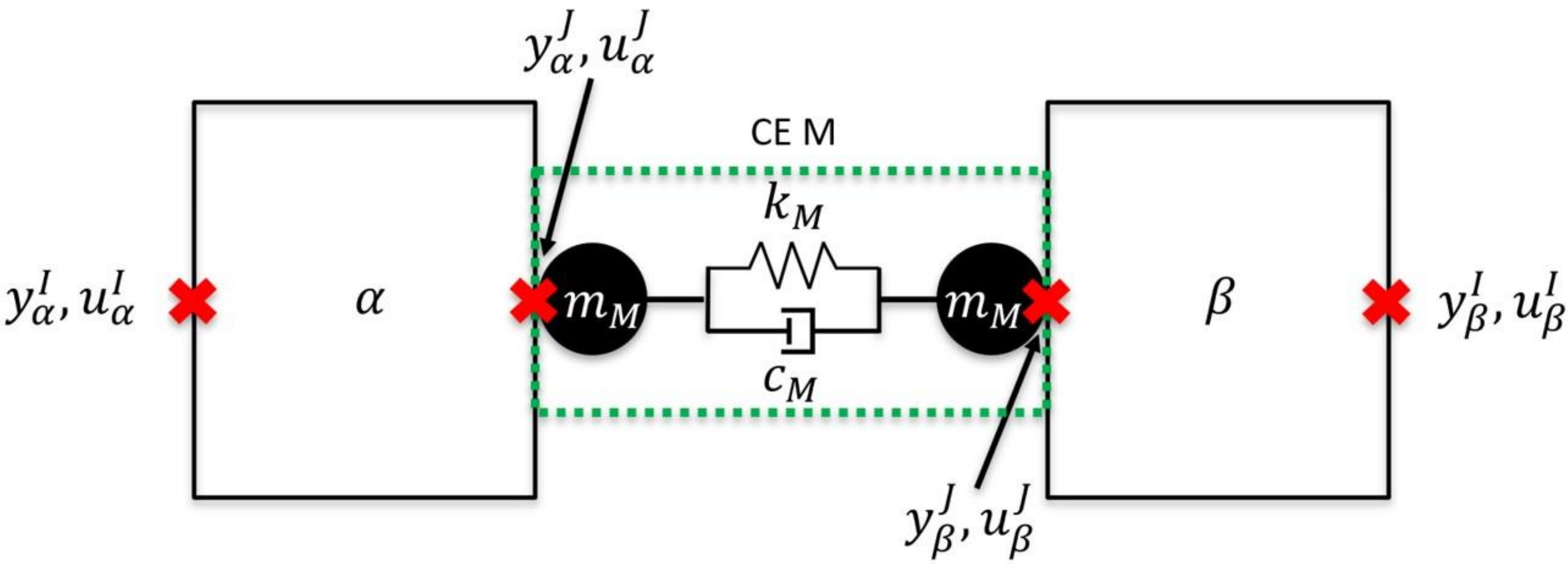}
  \caption{\textcolor{black}{Connecting element M, whose mass cannot be ignored.}}
  \label{fig:Numerical_assembly_1}
\end{subfigure}%
\begin{subfigure}{0.5\textwidth}
  \centering
  \includegraphics[width=0.7\linewidth]{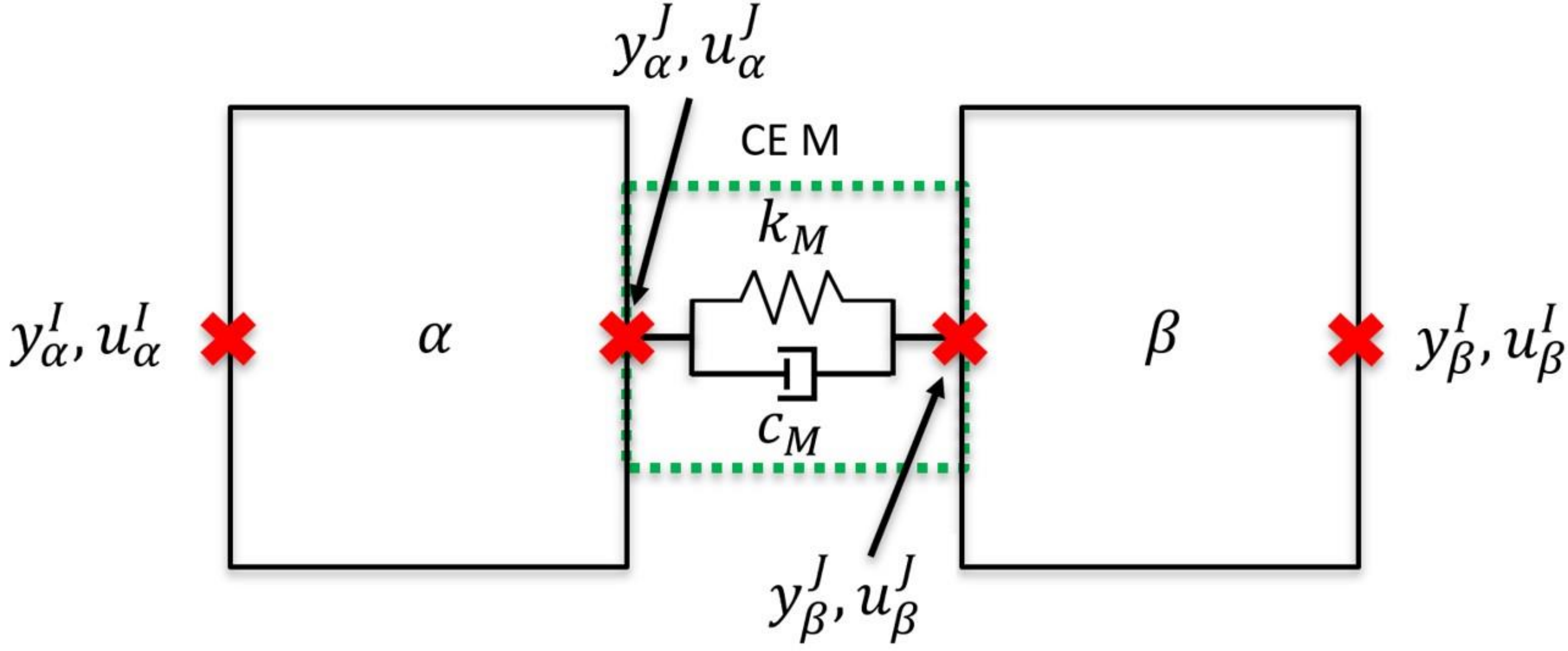}
  \caption{\textcolor{black}{Massless connecting element M.}}
  \label{fig:Numerical_assembly_2}
\end{subfigure}
\caption{\textcolor{black}{Assembled structure composed by two substructures $\alpha$ and $\beta$ connected by a connecting element M.}}
\label{fig:Numerical_assembly}
\end{figure}

In this section, the LM-SSS compatibility conditions will be relaxed (this can be seen as a weakening in the interface \citep{SWK_2020}) to develop a novel SSS technique. \textcolor{black}{To start, let us consider an assembled structure composed by two substructures ($\alpha$ and $\beta$) connected by a CE. Depending on the importance of the contribution of the mass of the CE to its dynamic behaviour in the frequency range of interest, we may represent the mentioned assembled structure in two different ways. Assuming that the mass of the CE cannot be ignored, the structure must be represented as depicted in figure \ref{fig:Numerical_assembly_1}. Whereas, if the CE is massless, we may represent it as given in figure \ref{fig:Numerical_assembly_2}.}

\color{black}

\begin{figure}[ht]
\centering
    \includegraphics[width=0.8\textwidth]{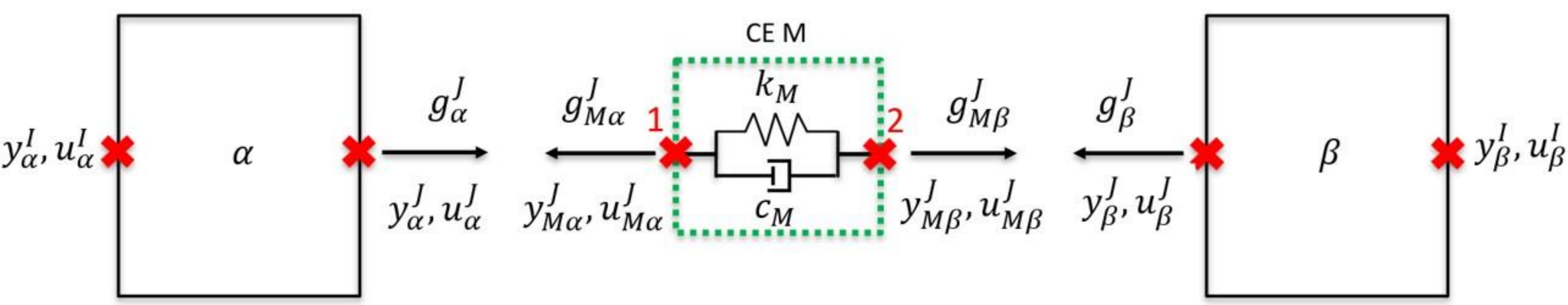}
    \caption{\textcolor{black}{Structure of figure \ref{fig:Numerical_assembly_2} with separated components.}}
     \label{fig:Free_Body_Diagram_Non_Rigid_Coupling}
\end{figure}

Assuming that the assembly under analysis is composed by a massless CE, if we separate the substructures and the CE, connecting forces ($g^{J}_{\alpha}$, $g^{J}_{\beta}$, $g^{J}_{M\alpha}$ and $g^{J}_{M\beta}$) will appear at their interfaces as depicted by figure \ref{fig:Free_Body_Diagram_Non_Rigid_Coupling}.

Further assuming that the CE does not present cross couplings between its DOFs (hence, that it is suitable to be characterized by IS), its off diagonal dynamic stiffness terms and the diagonal ones are related as follows \citep{JM_2017},\citep{JM_T_2017},\citep{MH_2020}:

\begin{equation}\label{eq:Z_IS_relation_1}
[Z^{\alpha}_{M,11}(j\omega)]= -[Z^{\alpha}_{M,12}(j\omega)]= -[Z^{\alpha}_{M,21}(j\omega)]= [Z^{\alpha}_{M,22}(j\omega)] 
\end{equation}

where, $[Z^{\alpha}]$ represents a dynamic stiffness matrix and the subscript $M$ denotes variables related to the CE. As for the other numbers in the subscript, the first refers to the outputs and the second to the inputs of each variable. The number $1$ means that the outputs or/and inputs are associated to the DOFs on the side $M\alpha$ of the CE, the number $2$ tags the outputs or/and inputs associated to the DOFs on the $M\beta$ side (see figure \ref{fig:Free_Body_Diagram_Non_Rigid_Coupling}).

When a structure is composed by two substructures connected by a CE, whose diagonal and off diagonal stiffness terms are related as given in expression \ref{eq:Z_IS_relation_1}, it can be demonstrated that the difference between the displacements on the interface of both substructures is given as follows (see, \cite{EB_2014})

\begin{equation}\label{eq:Compatibility_Condition_Non_Rigid_Formulation}
[B_{C}]\{\begin{matrix}
Y(j\omega)
\end{matrix}\}=[Z^{\alpha}_{M,22}(j\omega)]^{-1}\{\Lambda(j\omega)\}
\end{equation}

where, $[B_{C}]$ is a signed Boolean mapping matrix, which is constructed by giving unitary negative coefficients to the interface DOFs of one substructure and positive ones to the interface DOFs of the other substructure \citep{DK06}. The term $\{Y(j\omega)\}$ is given as follows

\begin{equation}\label{eq:Y(jomega)}
\{Y(j\omega)\}=\left\{\begin{matrix}
Y^{I}_{\alpha}(j\omega) & Y^{J}_{\alpha}(j\omega) & Y^{J}_{\beta}(j\omega) & Y^{I}_{\beta}(j\omega)
\end{matrix}\right\}^{T}
\end{equation}

while, $\{\Lambda(j\omega)\}$ is the vector of Lagrange Multipliers that represent connecting forces \citep{DK06},\citep{DK_20081169}. The relation between $\{\Lambda(j\omega)\}$ and the connecting forces acting on the interface of substructures $\alpha$ and $\beta$ is given as follows (see figure \ref{fig:Free_Body_Diagram_Non_Rigid_Coupling})

\begin{equation}\label{eq:Relation_connecting_forces_lambda}
\left\{\begin{matrix}
\{G_{\alpha}(j\omega)\}\\
\{G_{\beta}(j\omega)\}
\end{matrix}\right\}=-[B_{C}]^{T}\{\Lambda(j\omega)\}
\end{equation}

where, vectors $\{G_{\alpha}(j\omega)\}$ and $\{G_{\beta}(j\omega)\}$ are given below.

    \begin{subequations}\label{eq:connecting_forces_vector}
        \noindent
        \begin{tabularx}{\linewidth}{XX}
        \begin{equation}
        \{G_{\alpha}(j\omega)\}= \left\{\begin{matrix}
        \{0\}\\ 
        \{G^{J}_{\alpha}(j\omega)\} 
        \end{matrix} \right\} \label{eq:connecting_forces_vector_alpha}
        \end{equation}
        &
        \begin{equation}
        \{G_{\beta}(j\omega)\}= \left\{\begin{matrix}
        \{0\}\\
        \{G^{J}_{\beta}(j\omega)\}
        \end{matrix}\right\} \label{eq:connecting_forces_vector_beta}
        \end{equation}
    \end{tabularx}
    \end{subequations}



Note that, the relations given by equation \eqref{eq:Relation_connecting_forces_lambda} continue to be valid in time domain.

Expression \eqref{eq:Compatibility_Condition_Non_Rigid_Formulation} establishes the relative motion between the interface of substructures $\alpha$ and $\beta$, therefore it represents the compatibility condition of our problem. However, this condition is expressed by using the dynamic stiffness terms of the CE. In order to arrive to a compatibility condition written in time domain and given by state-space model terms, we must understand the physical meaning of the term $[Z^{\alpha}_{M,22}(j\omega)]^{-1}\{\Lambda(j\omega)\}$.

\color{black}

The term $[Z^{\alpha}_{M,22}(j\omega)]^{-1}\{\Lambda(j\omega)\}$ represents the displacement at the interface $M\beta$ of the CE (see figure \ref{fig:Free_Body_Diagram_Non_Rigid_Coupling}) when the connecting forces $\{\Lambda(j\omega)\}$ are applied at that location.
Hence, to represent $[Z^{\alpha}_{M,22}(j\omega)]^{-1}\{\Lambda(j\omega)\}$ in time domain we may establish a state-space model for the CE, whose outputs will be the displacements at the interface $M\beta$ and the inputs will be the connecting forces applied at that location (see figure \ref{fig:Free_Body_Diagram_Non_Rigid_Coupling}). This model must be computed as follows:

\begin{equation}\label{eq:disp_state_space_model_Mount_beta_interface_force}
\begin{gathered}
\{\dot{x}_{M}(t)\}=[A_{M}]\{x_{M}(t)\}+\left[\begin{matrix}
B_{M,2}^{J}
\end{matrix}
\right]\{\begin{matrix}
\lambda(t)
\end{matrix}\}\\
\left\{\begin{matrix}
y_{M\beta}^{J}(t)
\end{matrix}\right\}=\left[\begin{matrix}
C_{M,2}^{disp,J}
\end{matrix}\right]\{x_{M}(t)\}
\end{gathered}
\end{equation}

where, superscript $disp$ denotes a state-space matrix of a displacement state-space model (a state-space
model whose output vector elements are displacements).

\color{black}

Expression \eqref{eq:disp_state_space_model_Mount_beta_interface_force} is the state-space representation of the term $[Z^{\alpha}_{M,22}(j\omega)]^{-1}\{\Lambda(j\omega)\}$, hence it is the inverted state-space model representative of the diagonal dynamic stiffness term of the CE (the procedure to invert state-space models is presented in \ref{Inverting_state_space_models}). Thus, by using expressions \eqref{eq:Compatibility_Condition_Non_Rigid_Formulation} and \eqref{eq:disp_state_space_model_Mount_beta_interface_force} we may define the identity given as follows:

\begin{equation}\label{eq:relation_yM_and_yM_coupled}
\{y_{M\beta}^{J}(t)\}=[B_{C}]\{y(t)\}
\end{equation}

where $\{y(t)\}$ is given by equation \eqref{eq:y(t)}.

\begin{equation}\label{eq:y(t)}
\{y(t)\}=\left\{\begin{matrix}
y^{I}_{\alpha}(t) & y^{J}_{\alpha}(t) & y^{J}_{\beta}(t) & y^{I}_{\beta}(t)
\end{matrix}\right\}^{T}
\end{equation}

By using equation \eqref{eq:relation_yM_and_yM_coupled}, we may establish the state-space model given below. 

\color{black}

\begin{equation}\label{eq:disp_state_space_model_Mount_beta_interface_force_disp}
\begin{gathered}
\{\dot{x}_{M}(t)\}=[A_{M}]\{x_{M}(t)\}+\left[\begin{matrix}
B_{M,2}^{J}
\end{matrix}
\right]\left\{\begin{matrix}
\lambda(t)
\end{matrix}
\right\}\\
[B_{C}]\{\begin{matrix}
y(t)
\end{matrix}\}=\left[\begin{matrix}
C_{M,2}^{disp,J}
\end{matrix}
\right]\{x_{M}(t)\}
\end{gathered}
\end{equation}

The output equation of the state-space model given by expression \eqref{eq:disp_state_space_model_Mount_beta_interface_force_disp} represents the compatibility condition of our problem written in state-space model terms. However, to perform coupling by following the LM-SSS formulation the compatibility condition should be defined in terms of acceleration and therefore the state-space model given in expression \eqref{eq:disp_state_space_model_Mount_beta_interface_force_disp} must be double differentiated.

By calculating the second order time derivative of the state-space model given by expression \eqref{eq:disp_state_space_model_Mount_beta_interface_force_disp} (see \citep{FL_1988}), we obtain:

\begin{equation}\label{eq:disp_state_space_model_Mount_beta_interface_force_2}
\begin{gathered}
\{\dot{x}_{M}(t)\}=[A_{M}]\{x_{M}(t)\}+\left[\begin{matrix}
B_{M,2}^{J}
\end{matrix}
\right]\{\begin{matrix}
\lambda(t)
\end{matrix}\}\\
[B_{C}]\{\begin{matrix}
\ddot{y}(t)
\end{matrix}\}=\left[\begin{matrix}
C_{M,2}^{accel,J}
\end{matrix}
\right]\{x_{M}(t)\}+\left[\begin{matrix}
D_{M,22}^{accel,JJ}
\end{matrix}
\right]\{\begin{matrix}
\lambda(t)
\end{matrix}\}
\end{gathered}
\end{equation}

where,

    \begin{subequations}\label{eq:acceleration_ss_matrices}
        \noindent
        \begin{tabularx}{\linewidth}{XX}
        \begin{equation}
        \left[C_{M,2}^{accel,J}\right]= \left[\begin{matrix} C_{M,2}^{disp,J}A_{M}A_{M} 
        \end{matrix} \right]  \label{eq:acceleration_output_matrix}
        \end{equation}
        &
        \begin{equation}
        \left[D_{M,22}^{accel,JJ}\right]=
        \left[\begin{matrix}
        C_{M,2}^{disp,J}A_{M}B^{J}_{M,2}
        \end{matrix} \right] \label{eq:acceleration_feed_through_matrix}
        \end{equation}
    \end{tabularx}
    \end{subequations}

and superscript $accel$ denotes a state-space matrix of an acceleration state-space model (a state-space model whose output vector elements are acclerations). Note that, the state-space model given in expression \eqref{eq:disp_state_space_model_Mount_beta_interface_force_2} is the inverted state-space model representative of the diagonal dynamic mass term of the CE. 

The output equation of the state-space model given by equation \eqref{eq:disp_state_space_model_Mount_beta_interface_force_2} represents the compatibility conditions of our problem in terms of accelerations, when one CE is intended to be included into coupling. However, we can generalize this expression to provide the compatibility conditions, when an unlimited number of CEs is intended to be coupled. To perform such achievement, we must construct a diagonal connecting elements state-space model as follows:

\begin{equation}\label{eq:ss_Mounts_explicit_form}
\begin{gathered}
\left\{\begin{matrix}
\{\dot{x}_{M_{1}}(t)\}\\
\{\dot{x}_{M_{2}}(t)\}\\
\vdots
\end{matrix}
\right\}=\left[
A_{M,D}
\right]\left\{\begin{matrix}
\{x_{M_{1}}(t)\}\\
\{x_{M_{2}}(t)\}\\
\vdots
\end{matrix}
\right\}+\left[
B_{M,D} 
\right]\left\{\begin{matrix}
\{\lambda_{\alpha}(t)\}\\
\{\lambda_{\beta}(t)\}\\
\vdots
\end{matrix}
\right\}\\
\left[B_{C}\right]\left\{\begin{matrix}
\{\ddot{y}_{\alpha}(t)\}\\
\{\ddot{y}_{\beta}(t)\}\\
\vdots
\end{matrix}
\right\}=\left[
C^{accel}_{M,D}
\right]\left\{\begin{matrix}
\{x_{M_{1}}(t)\}\\
\{x_{M_{2}}(t)\}\\
\vdots
\end{matrix}
\right\}+\left[
D^{accel}_{M,D}
\right]\left\{\begin{matrix}
\{\lambda_{\alpha}(t)\}\\
\{\lambda_{\beta}(t)\}\\
\vdots
\end{matrix}
\right\}
\end{gathered}
\end{equation}

where,

\begin{alignat}{2}\label{eq:MatrixABCD}
\begin{split}
[A_{M,D}]=\left[\begin{matrix}
A_{M_{1}} & &\\
& A_{M_{2}} & \\
& & \ddots
\end{matrix}
\right],\ \ \ & [B_{M,D}]=\left[\begin{matrix}
B^{J}_{M_{1},2} & &\\
& B^{J}_{M_{2},2} & \\
& & \ddots
\end{matrix}
\right]\\ 
[C^{accel}_{M,D}]=\left[\begin{matrix}
C^{accel,J}_{M_{1},2} & &\\
& C^{accel,J}_{M_{2},2} & \\
& & \ddots
\end{matrix}
\right],\ \ \ & [D^{accel}_{M,D}]=\left[\begin{matrix}
D^{accel,JJ}_{M_{1},22} & &\\
& D^{accel,JJ}_{M_{2},22} & \\
& & \ddots
\end{matrix}
\right]
\end{split}
\end{alignat}

$\left\{ x(t) \right\} \in \mathbb{R}^{n \times 1}$ represents the state vector, which presents $n$ states and $\left\{ u(t) \right\} \in \mathbb{R}^{n_{i} \times 1}$ represents the input vector, whose elements are forces. Subscripts $M_{1}$ and $M_{2}$ denote variables associated to different CEs, while subscript $D$ denotes a block diagonal matrix.

By observing figure \ref{fig:Free_Body_Diagram_Non_Rigid_Coupling}, it is well evident that a diagonal coupled state-space model similar to the one presented in \citep{RD_2021} can be developed as follows:

\begin{equation}\label{eq:simplified_coupled_ss_model_from_uncoupled_ss_0}
\begin{gathered}
\left\{\begin{matrix}
\{\dot{x}_{\alpha}(t)\}\\
\{\dot{x}_{\beta}(t)\}\\
\vdots
\end{matrix}
\right\}=[A_{S,D}]\left\{\begin{matrix}
\{x_{\alpha}(t)\}\\
\{x_{\beta}(t)\}\\
\vdots
\end{matrix}
\right\}+[B_{S,D}]\left(\left\{\begin{matrix}
\{u_{\alpha}(t)\}\\
\{u_{\beta}(t)\}\\
\vdots
\end{matrix}
\right\}-[B_{C}]^{T}\left\{\begin{matrix}
\{\lambda_{\alpha}(t)\}\\
\{\lambda_{\beta}(t)\}\\
\vdots
\end{matrix}
\right\}\right)\\
\left\{\begin{matrix}
\{\ddot{y}_{\alpha}(t)\}\\
\{\ddot{y}_{\beta}(t)\}\\
\vdots
\end{matrix}
\right\}=[C^{accel}_{S,D}]\left\{\begin{matrix}
\{x_{\alpha}(t)\}\\
\{x_{\beta}(t)\}\\
\vdots
\end{matrix}
\right\}+[D^{accel}_{S,D}]\left(\left\{\begin{matrix}
\{u_{\alpha}(t)\}\\
\{u_{\beta}(t)\}\\
\vdots
\end{matrix}
\right\}-[B_{C}]^{T}\left\{\begin{matrix}
\{\lambda_{\alpha}(t)\}\\
\{\lambda_{\beta}(t)\}\\
\vdots
\end{matrix}
\right\}\right)
\end{gathered}
\end{equation}

where,

\begin{alignat}{2}\label{eq:MatrixABCD2}
\begin{split}
[A_{S,D}]=\left[
\begin{matrix}
A_{\alpha} & &\\
 & A_{\beta} &\\
 & & \ddots
\end{matrix}
\right],\ \ \ & [B_{S,D}]=\left[
\begin{matrix}
B_{\alpha} & &\\
 & B_{\beta} &\\
 & & \ddots
\end{matrix}
\right]\\ 
[C^{accel}_{S,D}]=\left[\begin{matrix}
C^{accel}_{\alpha} & &\\
 & C^{accel}_{\beta} &\\
 & & \ddots
\end{matrix}
\right],\ \ \ & [D^{accel}_{S,D}]=\left[
\begin{matrix}
D^{accel}_{\alpha} & &\\
 & D^{accel}_{\beta} &\\
 & & \ddots
\end{matrix}
\right]
\end{split}
\end{alignat}

and subscript $S$ denote variables related to the substructures.

The diagonal coupled state-space model given by expression \eqref{eq:simplified_coupled_ss_model_from_uncoupled_ss_0} can be rewritten in a more compact form as in equation \eqref{eq:simplified_coupled_ss_model_from_uncoupled_ss}.

\begin{equation}\label{eq:simplified_coupled_ss_model_from_uncoupled_ss}
\begin{gathered}
\{\dot{x}_{S}(t)\}=[A_{S,D}]\{x_{S}(t)\}+[B_{S,D}](\{u(t)\}-[B_{C}]^{T}
\{\lambda(t)\})\\
\{\ddot{y}(t)\}=[C^{accel}_{S,D}]\{x_{S}(t)\}+[D^{accel}_{S,D}](\{u(t)\}-[B_{C}]^{T}\{\lambda(t)\})
\end{gathered}
\end{equation}

When coupled, the substructures must fulfill two different conditions, compatibility and equilibrium. The compatibility condition was already computed, being given by the output equations of the diagonal connecting element state-space model (see expression \eqref{eq:ss_Mounts_explicit_form}). The equilibrium conditions can be fulfilled by forcing the substructures to verify the local equilibrium conditions given by the output equations of the diagonal coupled state-space model (see equation \eqref{eq:simplified_coupled_ss_model_from_uncoupled_ss}). Rearranging the second equation, the following expression is obtained. 

\begin{equation}\label{eq:equilibrium_condition_rearranged}
\{u(t)\}=[D^{accel}_{S,D}]^{-1}(\{\ddot{y}(t)\}-[C^{accel}_{S,D}]\{x_{S}(t)\})+[B_{C}]^{T}\{\lambda(t)\}
\end{equation}

Using together the equilibrium and compatibility conditions (equation \eqref{eq:equilibrium_condition_rearranged} and the output equations \eqref{eq:ss_Mounts_explicit_form}, respectively), and dropping $\{\bullet\}$, $[\bullet]$ and $(t)$ for ease of readability, we can establish the following system of equations.

\begin{equation}\label{eq:coupling_conditions_non_rigid_formulation}\begin{cases}
(D^{accel}_{S,D})^{-1}(\ddot{y}-C^{accel}_{S,D}x_{S})+B_{C}^{T}\lambda=u\\
B_{C}\ddot{y}=C^{accel}_{M,D}x_{M}+D^{accel}_{M,D}\lambda
\end{cases}
\end{equation}

After some mathematical manipulations, the system of equations \eqref{eq:coupling_conditions_non_rigid_formulation}, can be rewritten as in equation \eqref{eq:coupling_conditions_non_rigid_formulation_1}.

\begin{equation}\label{eq:coupling_conditions_non_rigid_formulation_1}
\begin{cases}
\lambda=(B_{C}D^{accel}_{S,D}B_{C}^{T}+D^{accel}_{M,D})^{-1}(-C^{accel}_{M,D}x_{M}+B_{C}C^{accel}_{S,D}x_{S}+B_{C}D^{accel}_{S,D}u)\\
\ddot{y}=(C^{accel}_{S,D}-D^{accel}_{S,D}B_{C}^{T}(B_{C}D^{accel}_{S,D}B_{C}^{T}+D^{accel}_{M,D})^{-1}B_{C}C^{accel}_{S,D})x_{S}+D^{accel}_{S,D}B_{C}^{T}(B_{C}D^{accel}_{S,D}B_{C}^{T}+D^{accel}_{M,D})^{-1}C^{accel}_{M,D}x_{M}\\
\ \ \ +(D^{accel}_{S,D}-D^{accel}_{S,D}B_{C}^{T}(B_{C}D^{accel}_{S,D}B_{C}^{T}+D^{accel}_{M,D})^{-1}B_{C}D^{accel}_{S,D})u
\end{cases}
\end{equation}

By using equation \eqref{eq:equilibrium_condition_rearranged} and the bottom equation of the system of equations \eqref{eq:coupling_conditions_non_rigid_formulation_1}, the state equation of the diagonal coupled state-space model (equation \eqref{eq:simplified_coupled_ss_model_from_uncoupled_ss}) can be rewritten as follows.

\begin{equation}\label{eq:coupled_state_equation_structures}
\begin{split}
\dot{x}_{S} & =(A_{S,D}-B_{S,D}B_{C}^{T}(B_{C}D^{accel}_{S,D}B_{C}^{T}+D^{accel}_{M,D})^{-1}B_{C}C^{accel}_{S,D})x_{S}+B_{S,D}B^{T}(B_{C}D^{accel}_{S,D}B_{C}^{T}+D^{accel}_{M,D})^{-1}C^{accel}_{M,D}x_{M}\\
& +(B_{S,D}-B_{S,D}B_{C}^{T}(B_{C}D^{accel}_{S,D}B_{C}^{T}+D^{accel}_{M,D})^{-1}B_{C}D^{accel}_{S,D})u
\end{split}
\end{equation}

Using the value of $\{\lambda(t)\}$ given by the upper equation of the system of equations \eqref{eq:coupling_conditions_non_rigid_formulation_1}, the state equation of the state-space model given by expression \eqref{eq:ss_Mounts_explicit_form} can be rewritten as in equation \eqref{eq:coupled_state_equation_mount}.

\begin{equation}\label{eq:coupled_state_equation_mount}
\begin{split}
\dot{x}_{M} & =B_{M,D}(B_{C}D^{accel}_{S,D}B_{C}^{T}+D^{accel}_{M,D})^{-1}B_{C}C^{accel}_{S,D}x_{S}+(A_{M,D}-B_{M,D}(B_{C}D^{accel}_{S,D}B_{C}^{T}+D^{accel}_{M,D})^{-1}C^{accel}_{M,D})x_{M}\\
& +B_{M,D}(B_{C}D^{accel}_{S,D}B_{C}^{T}+D^{accel}_{M,D})^{-1}B_{C}D^{accel}_{S,D}u
\end{split}
\end{equation}

From equations \eqref{eq:coupled_state_equation_structures}, \eqref{eq:coupled_state_equation_mount} and the bottom equation of the system of equations \eqref{eq:coupling_conditions_non_rigid_formulation_1}, we obtain the following coupled state-space model:

\begin{equation}\label{eq:coupled_ss_non_rigid_formulation}
\begin{gathered}
\left\{\begin{matrix}
\dot{\bar{x}}_{S}(t)\\
\dot{\bar{x}}_{M}(t)
\end{matrix}\right\}=[\bar{A}]\left\{\begin{matrix}
\bar{x}_{S}(t)\\
\bar{x}_{M}(t)
\end{matrix}\right\}+[\bar{B}]\{\bar{u}(t)\}\\
\{\ddot{\bar{y}}(t)\}=[\bar{C}^{accel}]\left\{\begin{matrix}
\bar{x}_{S}(t)\\
\bar{x}_{M}(t)
\end{matrix}\right\}+[\bar{D}^{accel}]\{\bar{u}(t)\}
\end{gathered}
\end{equation}

where, 

\begin{alignat}{2}\label{eq:MatrixABCD_coup}
\begin{split}
[\bar{A}]=\left[\begin{matrix}
\bar{A}_{SS} & \bar{A}_{SM} \\
\bar{A}_{MS} & \bar{A}_{MM}
\end{matrix}\right],\ \ \ & [\bar{B}]=\left[\begin{matrix}
\bar{B}_{S}\\
\bar{B}_{M}
\end{matrix}\right]\\ 
[\bar{C}^{accel}]=\left[\begin{matrix}
\bar{C}^{accel}_{S} & \bar{C}^{accel}_{M}
\end{matrix}\right],\ \ \ & [\bar{D}^{accel}]=\left[
\begin{matrix}
\bar{D}^{accel}_{S}
\end{matrix}
\right]
\end{split}
\end{alignat}

and over-lined variables represent variables of the coupled state-space model. The coupled state-space matrices can be computed from the diagonal coupled and diagonal connecting elements state-space models (expressions \eqref{eq:simplified_coupled_ss_model_from_uncoupled_ss} and \eqref{eq:disp_state_space_model_Mount_beta_interface_force_2}, respectively) as follows.

\begin{equation}\label{eq:coupled_matrices_value_non_rigid_formulation}
\begin{gathered}
[\bar{A}_{SS}]=A_{S,D}-B_{S,D}B_{C}^{T}(B_{C}D^{accel}_{S,D}B_{C}^{T}+D^{accel}_{M,D})^{-1}B_{C}C^{accel}_{S,D}\\
[\bar{A}_{SM}]=B_{S,D}B_{C}^{T}(B_{C}D^{accel}_{S,D}B_{C}^{T}+D^{accel}_{M,D})^{-1}C^{accel}_{M,D}\\
[\bar{A}_{MS}]=B_{M,D}(B_{C}D^{accel}_{S,D}B_{C}^{T}+D^{accel}_{M,D})^{-1}B_{C}C^{accel}_{S,D}\\
[\bar{A}_{MM}]=A_{M,D}-B_{M,D}(B_{C}D^{accel}_{S,D}B_{C}^{T}+D^{accel}_{M,D})^{-1}C^{accel}_{M,D}\\
[\bar{B}_{S}]=B_{S,D}-B_{S,D}B_{C}^{T}(B_{C}D^{accel}_{S,D}B_{C}^{T}+D^{accel}_{M,D})^{-1}B_{C}D^{accel}_{S,D}\\
[\bar{B}_{M}]=B_{M,D}(B_{C}D^{accel}_{S,D}B_{C}^{T}+D^{accel}_{M,D})^{-1}B_{C}D^{accel}_{S,D}\\
[\bar{C}^{accel}_{S}]=C^{accel}_{S,D}-D^{accel}_{S,D}B_{C}^{T}(B_{C}D^{accel}_{S,D}B_{C}^{T}+D^{accel}_{M,D})^{-1}B_{C}C^{accel}_{S,D}\\
[\bar{C}^{accel}_{M}]=D^{accel}_{S,D}B_{C}^{T}(B_{C}D^{accel}_{S,D}B_{C}^{T}+D^{accel}_{M,D})^{-1}C^{accel}_{M,D}\\
[\bar{D}^{accel}_{S}]=D^{accel}_{S,D}-D^{accel}_{S,D}B_{C}^{T}(B_{C}D^{accel}_{S,D}B_{C}^{T}+D^{accel}_{M,D})^{-1}B_{C}D^{accel}_{S,D}
\end{gathered}
\end{equation}

By using expressions \eqref{eq:coupled_ss_non_rigid_formulation} and \eqref{eq:coupled_matrices_value_non_rigid_formulation} it is possible to compute the coupled acceleration state-space model (a state-space model whose output vector elements are accelerations). To get displacement and velocity coupled state-space models it is sufficient to adjust the output and feed-through state-space matrices in expressions \eqref{eq:coupled_matrices_value_non_rigid_formulation} \textcolor{black}{(see \ref{Extending the LM-SSS-WI formulation to directly compute coupled displacement and velocity state-space models})}.

\textcolor{black}{For CEs, whose diagonal and off diagonal dynamic stiffness terms are not related as given in expression \eqref{eq:Z_IS_relation_1}, equation \eqref{eq:Relation_connecting_forces_lambda}} will not be valid, hence this formulation cannot be used. For those situations, it is suggested to treat the CE as another substructure to be coupled and perform coupling by using the LM-SSS method. Applying this suggestion to the structure depicted in figure \ref{fig:Free_Body_Diagram_Non_Rigid_Coupling}, we would assume two rigid couplings, one between the interface of component $\alpha$ and the interface $M\alpha$ of the CE and the other between the interface of component $\beta$ and the interface $M\beta$ of the CE.

\section{Minimal-order coupled state-space models}\label{Minimal-order coupled state-space models}

The \textcolor{black}{LM-SSS with compatibility relaxation} method developed in section \ref{LM_SSS_Non_Rigid_Formulation} is valid to couple state-space models previously transformed into coupling form. However, the computation of the minimal-order coupled state-space model from the one directly obtained by using \textcolor{black}{this coupling approach} deserves some attention.

Let us assume that the coupling of the components presented in figure \ref{fig:Free_Body_Diagram_Non_Rigid_Coupling} is performed by using \textcolor{black}{LM-SSS with compatibility relaxation}. Let us further assume that the diagonal connecting elements state-space model (see equation \eqref{eq:disp_state_space_model_Mount_beta_interface_force_2}) and the state-space models of the substructures were previously transformed into coupling form. By dropping $\{\bullet\}$ and $(t)$ for ease of readability, the state vectors of the state-space models of the components transformed into coupling form would be given as follows:

    \begin{subequations}\label{eq:state_vector_coupling_form}
        \noindent
        \begin{tabularx}{\linewidth}{XX}
        \begin{equation}
        \{z_{\alpha}\}=\left\{\begin{matrix}
        \dot{y}^{J}_{\alpha} &
        y^{J}_{\alpha} &
        x^{I}_{\alpha}
        \end{matrix}
        \right\}^{T} 
        \label{eq:sv_cf_alpha}
        \end{equation}
        &
        \begin{equation}
        \{z_{\beta}\}=\left\{\begin{matrix}
        \dot{y}^{J}_{\beta} &
        y^{J}_{\beta} &
        x^{I}_{\beta}
        \end{matrix}
        \right\}^{T} 
        \label{eq:sv_cf_beta}
        \end{equation}
    \end{tabularx}
    \end{subequations}

where $z$ represents the state vector of a state-space model transformed into coupling form.

\color{black}

By following the notation used in equation \eqref{eq:disp_state_space_model_Mount_beta_interface_force}, the state vector of the diagonal connecting elements state-space model transformed into coupling form will be given below.

\begin{equation}\label{eq:diagonal_joint_element_ss_coup_form_sv}
\{z_{M}\}=\left\{\begin{matrix}
\dot{y}^{J}_{M\beta} & y^{J}_{M\beta} & x^{I}_{M}
\end{matrix}
\right\}^{T}
\end{equation}

By using \textcolor{black}{LM-SSS with compatibility relaxation} to couple the transformed state-space models and having in mind the identity provided by equation \eqref{eq:relation_yM_and_yM_coupled}, the state vector of the obtained coupled state-space model will be given as follows:

\begin{equation}\label{eq:sv_coupled_ss_model_coup_form}
\{\bar{z}_{\alpha M \beta}\}=\left\{\begin{matrix}
\dot{\bar{y}}^{J}_{\alpha} &
\bar{y}^{J}_{\alpha} &
\bar{x}^{I}_{\alpha} &
\dot{\bar{y}}^{J}_{\beta} &
\bar{y}^{J}_{\beta} &
\bar{x}^{I}_{\beta} &
[B_{C}]\{\dot{\bar{y}}\} &
[B_{C}]\{\bar{y}\} &
\bar{x}^{I}_{M}
\end{matrix}
\right\}^{T}
\end{equation}

where subscript $\alpha M \beta$ denotes a variable related to the assembled structure depicted in figure \ref{fig:Numerical_assembly_2} and $\{\bar{y}\}$ is given below.

\begin{equation}\label{eq:coupled_y_vector}
\{\bar{y}\}=\left\{\begin{matrix}
\dot{\bar{y}}^{I}_{\alpha} &
\bar{y}^{J}_{\alpha} &
\bar{y}^{J}_{\beta} &
\bar{y}^{I}_{\beta}
\end{matrix}
\right\}^{T}
\end{equation}

By analyzing the state vector of the coupled state-space model given by equation \eqref{eq:sv_coupled_ss_model_coup_form}, it is clear that this vector does not contain states that represent the same physical quantity. Instead, thanks to the inclusion of the CE through the use of the diagonal connecting elements state-space model (equation \eqref{eq:disp_state_space_model_Mount_beta_interface_force_2}), this model presents states whose physical meaning is the difference between the interface outputs of the connected structures. Even though these states cannot be considered as redundant ones, since they represent different physical quantities, their contribution can be retrieved from the interface outputs of the connected structures, hence they can be eliminated. 

To obtain the minimal-order form of the computed coupled state-space model, the use of the post-processing procedure that relies on the construction of a Boolean localization matrix, $[L_{T}]$, presented in \citep{RD_2021} still holds. However, the construction of the state mapping matrix $[B_{T}]$ that is used to compute $[L_{T}]$ will depend on the construction of the $[B_{C}]$ matrix. To ease the computation of $[B_{T}]$ let us assume that $[B_{C}]$ is constructed as follows

\begin{equation}\label{eq:B_T_matrix_non_rigid_form_ex}
[B_{C}]\{y\}=\left[\begin{matrix}
0 & I^{J}_{\alpha} & I^{J}_{\beta} & 0
\end{matrix}
\right]\left\{\begin{matrix}
\bar{y}^{I}_{\alpha}\\
\bar{y}^{J}_{\alpha}\\
\bar{y}^{J}_{\beta}\\
\bar{y}^{I}_{\beta}
\end{matrix}
\right\}=\left\{\begin{matrix}
0
\end{matrix}
\right\}
\end{equation}

In equation \eqref{eq:B_T_matrix_non_rigid_form_ex}, $I^{J}_{\alpha}$ and $I^{J}_{\beta}$ are given as follows:

    \begin{subequations}\label{eq:I_J_alpha_beta}
        \noindent
        \begin{tabularx}{\linewidth}{XX}
        \begin{equation}
        [I^{J}_{\alpha}]=\phi [I] \label{eq:I_J_alpha}
        \end{equation}
        &
        \begin{equation}
        [I^{J}_{\beta}]=-\phi [I] \label{eq:I_J_beta}
        \end{equation}
    \end{tabularx}
    \end{subequations}
    
where $[I]$ is an identity matrix with dimension $n_{J}$, being $n_{J}$ the number of connected DOFs. $\phi$ is a numeric coefficient, whose value can be selected to be either $-1$ or $1$. The role of this variable is to adjust the use of the post-processing procedures here outlined in compliance with the computation of $[B_{C}]$.

By using expressions \eqref{eq:B_T_matrix_non_rigid_form_ex}, \eqref{eq:I_J_alpha} and \eqref{eq:I_J_beta}, $[B_{T}]$ can be always constructed in compliance with $[B_{C}]$ by following equation \eqref{eq:B_T_matrix_non_rigid_form}.

\begin{equation}\label{eq:B_T_matrix_non_rigid_form}
[B_{T}]\{\bar{z}_{\alpha M \beta}\}=\left[\begin{matrix}
I^{J}_{\beta} & 0 & 0 & I^{J}_{\alpha} & 0 & 0 & I & 0 & 0\\
0 & I^{J}_{\beta} & 0 & 0 & I^{J}_{\alpha} & 0 & 0 & I & 0 
\end{matrix}
\right]\left\{\begin{matrix}
\dot{\bar{y}}^{J}_{\alpha}\\
\bar{y}^{J}_{\alpha}\\
\bar{x}^{I}_{\alpha}\\
\dot{\bar{y}}^{J}_{\beta}\\
\bar{y}^{J}_{\beta}\\
\bar{x}^{I}_{\beta}\\
\phi \left(\dot{\bar{y}}^{J}_{\alpha}-\dot{\bar{y}}^{J}_{\beta}\right)\\
\phi \left(\bar{y}^{J}_{\alpha}-\bar{y}^{J}_{\beta}\right)\\
\bar{x}^{I}_{M}
\end{matrix}
\right\}=\left\{\begin{matrix}
0
\end{matrix}
\right\}
\end{equation}

\color{black}

After using $[B_{T}]$ to calculate the $[L_{T}]$ matrix and after applying the post-processing procedure presented in \citep{RD_2021}, the coupled minimal-order state-space model would be obtained and its state vector would be given as shown hereafter.

\begin{equation}\label{eq:coupled_state_vector}
\begin{gathered}
\{\bar{z}_{\alpha M \beta}\}=\left\{\begin{matrix}
\dot{\bar{y}}_{\alpha}^{J} &
\bar{y}_{\alpha}^{J} &
\bar{x}_{\alpha}^{I} &
\dot{\bar{y}}_{\beta}^{J} &
\bar{y}_{\beta}^{J} &
\bar{x}_{\beta}^{I} & 
\bar{x}^{I}_{M}
\end{matrix}\right\}^{T}
\end{gathered}
\end{equation}

It is worth noticing that the procedure used above to compute minimal-order coupled state-space models from the ones obtained by using \textcolor{black}{LM-SSS with compatibility relaxation} can be extended to coupled state-space models representing mechanical systems composed of an unlimited number of components and CEs.

On the other hand, the manual post-processing procedure to compute minimal-order coupled state-space models presented in \citep{RD_2021} cannot be applied to the coupled state-space models obtained by using \textcolor{black}{LM-SSS with compatibility relaxation}. Nevertheless, we may outline a manual procedure that is suitable to compute minimal-order models from the coupled state-space models obtained by exploiting this coupling formulation.

To implement a manual post-processing procedure that mimics the one that relies on the use of a Boolean localization matrix, the following procedure should be used for each pair of the coupled interface outputs present on the state vector.

\begin{enumerate}
\item Add the column of $[\bar{A}]$ that is being multiplied by the difference of the coupled interface outputs previously multiplied by $\phi$ to the column that is being multiplied by the coupled interface output of the first substructure state-space model introduced in the diagonal coupled model (see equation \eqref{eq:simplified_coupled_ss_model_from_uncoupled_ss_0});
\item Add the column of $[\bar{A}]$ that is being multiplied by the difference of the coupled interface outputs previously multiplied by $-\phi$ to the column that is being multiplied by the coupled interface output of the second substructure state-space model introduced in the diagonal coupled model (see equation \eqref{eq:simplified_coupled_ss_model_from_uncoupled_ss_0});
\item Repeat the procedures outlined in the first and second points for the matrix $[\bar{C}]$;
\item Eliminate the row correspondent to the difference of the coupled interface outputs and the column of matrix $[\bar{A}]$ that is being multiplied by the same quantity;
\item Eliminate the same row and column of matrices $[\bar{B}]$ and $[\bar{C}]$, respectively;
\item Repeat the procedure for the first derivative of the analyzed pair of coupled interface DOFs.
\end{enumerate}

\color{black}

\section{Determining inverted state-space models representative of diagonal apparent mass terms of CEs}\label{Determining state-space models representative of diagonal apparent mass terms of connecting elements}

The coupling approach developed in section \ref{LM_SSS_Non_Rigid_Formulation} makes it possible to include CEs into the LM-SSS formulation by using state-space models representative of their inverted diagonal apparent mass terms. Depending on the selected strategy to determine these models (i.e. analytically, numerically or experimentally), different approaches must be applied to compute them. 

\subsection{Primal state-space formulation}\label{Primal formulation in time domain}

For an easier understanding on how the inverted state-space models representative of the diagonal apparent mass terms of a CE can be determined experimentally, in this section we will develop a primal state-space formulation. It is important to mention that when the substructures under analysis present just interface DOFs, both dual and primal formulations are equivalent \citep{MH_2020}. Whereas, if the substructures under analysis present both internal and interface DOFs, the use of dual formulations is preferred. Dual formulations just require the inversion of quantities related to the interface DOFs in order to perform coupling or decoupling, while primal ones in general require the inversion of all data, since accelerance is typically the estimated quantity from experimental tests. Nevertheless for specific applications, the use of primal formulations might present important advantages over dual ones. This applies, for instance, when a state-space model of a massless component is intended to be identified from the state-space model of the assembled structure wherein it is included. In fact, if this is the case, the use of a dual formulation would lead to an ill-conditioned matrix inversion, since the feed-through matrix of the identified model would be composed by infinite valued terms. However, an identification without facing any numerical problem would be achieved by using primal disassembly. By following this procedure the state-space models involved in the decoupling process would have to be inverted, but during primal disassembly none matrix inversion is performed (see equation \eqref{eq:impedance_diag_coup_ss_2}). This will prevent the
inversion of ill-conditioned matrices, leading to a successfully characterization of the massless component.

Here, the primal state-space formulation will be established by taking the primal formulation in frequency domain presented in \citep{DK_20081169} as main source of inspiration. To estimate state-space models representative of the dynamic behaviour of real substructures it is common practice to use information identified from the measured FRFs (in general from accelerance). Yet, when performing primal assembly in frequency domain, the substructures to be assembled are represented by their inverse FRFs. Hence, the state-space models to be primally assembled must either be constructed by using information from the inverted FRFs of the substructures to be coupled or by inverting the state-space models constructed from the measured FRFs. In general, the state-space models are constructed by following the second approach (see, for instance \citep{MG_2020}). This procedure makes the estimation of displacement state-space models possible. These models should then be double-differentiated to obtain the associated acceleration state-space models, since to invert a state-space model the inversion of its feed-through matrix is required (see \ref{Inverting_state_space_models}).

Let us set up a diagonal coupled state-space model as follows:

\begin{equation}\label{eq:ss_Mounts_explicit_form_impedance}
\begin{gathered}
\left\{\begin{matrix}
\{\dot{x}_{\alpha}(t)\}\\
\{\dot{x}_{\beta}(t)\}\\
\vdots
\end{matrix}
\right\}=\left[
A_{S,D}
\right]\left\{\begin{matrix}
\{x_{\alpha}(t)\}\\
\{x_{\beta}(t)\}\\
\vdots
\end{matrix}
\right\}+\left[
B_{S,D} 
\right]\left\{\begin{matrix}
\{\ddot{y}_{\alpha}(t)\}\\
\{\ddot{y}_{\beta}(t)\}\\
\vdots
\end{matrix}
\right\}\\
\left\{\begin{matrix}
\{u_{\alpha}(t)\}\\
\{u_{\beta}(t)\}\\
\vdots
\end{matrix}
\right\}+\left\{\begin{matrix}
\{g_{\alpha}(t)\}\\
\{g_{\beta}(t)\}\\
\vdots
\end{matrix}
\right\}=\left[
C_{S,D}
\right]\left\{\begin{matrix}
\{x_{\alpha}(t)\}\\
\{x_{\beta}(t)\}\\
\vdots
\end{matrix}
\right\}+\left[
D_{S,D}
\right]\left\{\begin{matrix}
\{\ddot{y}_{\alpha}(t)\}\\
\{\ddot{y}_{\beta}(t)\}\\
\vdots
\end{matrix}
\right\}
\end{gathered}
\end{equation}

where matrices $[A_{S,D}]$, $[B_{S,D}]$, $[C_{S,D}]$ and $[D_{S,D}]$ are given below.

\begin{alignat}{2}\label{eq:MatrixABCD_impedance}
\begin{split}
[A_{S,D}]=\left[\begin{matrix}
A_{\alpha} & &\\
& A_{\beta} & \\
& & \ddots
\end{matrix}
\right],\ \ \ & [B_{S,D}]=\left[\begin{matrix}
B_{\alpha} & &\\
& B_{\beta} & \\
& & \ddots
\end{matrix}
\right]\\
[C_{S,D}]=\left[\begin{matrix}
C_{\alpha} & &\\
& C_{\beta} & \\
& & \ddots
\end{matrix}
\right],\ \ \ & [D_{S,D}]=\left[\begin{matrix}
D_{\alpha} & &\\
& D_{\beta} & \\
& & \ddots
\end{matrix}
\right]
\end{split}
\end{alignat}

By using a more compact representation, the state-space model given by expression \eqref{eq:ss_Mounts_explicit_form_impedance} can be rewritten as in equation \eqref{eq:impedance_diag_coup_ss}.

\begin{equation}\label{eq:impedance_diag_coup_ss}
\begin{gathered}
\{\dot{x}_{S}(t)\}=[A_{S,D}]\{x_{S}(t)\}+[B_{S,D}]\{\ddot{y}(t)\}\\
\{u(t)\}+\{g(t)\}=[C_{S,D}]\{x_{S}(t)\}+[D_{S,D}]\{\ddot{y}(t)\}
\end{gathered}
\end{equation}

At this point, the Boolean localization matrix $[L]$ presented in \citep{DK_20081169} must be introduced. This matrix can be calculated by computing the nullspace of the mapping matrix $[B_{C}]$, and can be used to establish the following relations \citep{DK_20081169}:   

    \begin{subequations}\label{eq:coupling_conditions}
        \noindent
        \begin{tabularx}{\linewidth}{XX}
        \begin{equation}
        \{y(t)\}=[L]\{\bar{y}(t)\} \label{eq:L_outputs}
        \end{equation}
        &
        \begin{equation}
        [L]^{T}\{g(t)\}=\{0\} \label{eq:L_connecting_forces}
        \end{equation}
    \end{tabularx}
    \end{subequations}

where, $\{\bar{y}(t)\}$ represents the coupled output vector, which is composed by the unique set of displacements. 

In primal formulation, compatibility is enforced by retaining the unique set of DOFs, hence this coupling condition is mathematically given by equation \eqref{eq:L_outputs}. On the other hand, the equilibrium condition is mathematically given by equation \eqref{eq:L_connecting_forces}, because it enforces the mutual cancellation of the connecting forces that guarantees force equilibrium at the interfaces.

Since the inputs of the state-space models under analysis are accelerations (see equation \eqref{eq:impedance_diag_coup_ss}) the compatibility equation \eqref{eq:L_outputs} should be double differentiated.

\begin{equation}\label{eq:L_outputs_acc}
\{\ddot{y}(t)\}=[L]\{\ddot{\bar{y}}(t)\}
\end{equation}

By using equation \eqref{eq:L_outputs_acc}, expression \eqref{eq:impedance_diag_coup_ss} can be rewritten as follows.

\begin{equation}\label{eq:impedance_diag_coup_ss_1}
\begin{gathered}
\{\dot{x}_{S}(t)\}=[A_{S,D}]\{x_{S}(t)\}+[B_{S,D}][L]\{\ddot{\bar{y}}(t)\}\\
(\{u(t)\}+\{g(t)\})=[C_{S,D}]\{x_{S}(t)\}+[D_{S,D}][L]\{\ddot{\bar{y}}(t)\}
\end{gathered}
\end{equation}

By pre-multiplying the output equations of the state-space model given by expression \eqref{eq:impedance_diag_coup_ss_1} by $[L]^{T}$ and by using equation \eqref{eq:L_connecting_forces}, expression \eqref{eq:impedance_diag_coup_ss_1} can be rewritten as follows.

\begin{equation}\label{eq:impedance_diag_coup_ss_2}
\begin{gathered}
\{\dot{x}_{S}(t)\}=[A_{S,D}]\{x_{S}(t)\}+[B_{S,D}][L]\{\ddot{\bar{y}}(t)\}\\
[L^{T}]\{u(t)\}=[L]^{T}[C_{S,D}]\{x_{S}(t)\}+[L]^{T}[D_{S,D}][L]\{\ddot{\bar{y}}(t)\}
\end{gathered}
\end{equation}

The state-space model given by expression \eqref{eq:impedance_diag_coup_ss_2} represents the primally assembled model and can be rewritten as follows:

\begin{equation}\label{eq:impedance_diag_coup_ss_final}
\begin{gathered}
\{\dot{\bar{x}}_{S}(t)\}=[\bar{A}_{S,D}]\{\bar{x}_{S}(t)\}+[\bar{B}_{S,D}]\{\ddot{\bar{y}}(t)\}\\
\{\bar{u}(t)\}=[\bar{C}_{S,D}]\{\bar{x}_{S}(t)\}+[\bar{D}_{S,D}]\{\ddot{\bar{y}}(t)\}
\end{gathered}
\end{equation}

where, $\{\bar{u}(t)\}$ represents the unique set of forces and matrices $[\bar{A}_{S,D}]$, $[\bar{B}_{S,D}]$, $[\bar{C}_{S,D}]$ and $[\bar{D}_{S,D}]$ are given as in equation \eqref{eq:MatrixABCD_impedance_coupled}.

\begin{alignat}{2}\label{eq:MatrixABCD_impedance_coupled}
\begin{split}
[\bar{A}_{S,D}]=\left[A_{S,D}
\right],\ \ \ & 
[\bar{B}_{S,D}]=\left[B_{S,D}
\right][L],\\
[\bar{C}_{S,D}]=[L]^{T}\left[
C_{S,D}
\right],\ \ \ &
[\bar{D}_{S,D}]=[L]^{T}\left[
D_{S,D}
\right][L]
\end{split}
\end{alignat}

Note that, expressions \eqref{eq:MatrixABCD_impedance_coupled} can also be used to perform primal disassembly. However, the state-space models of the substructures to be disassembled must be previously transformed into negative form by following the procedure presented in \ref{Negative form of a state-space model representative of apparent mass}.

It is worth noticing that the assembly/disassembly of state-space models previously transformed into coupling form (see, \citep{SJO_20072697},\citep{AMRDvMTPATMR_2020}) continues to be valid when primal formulation is used. The transformation into coupling form can be either applied to state-space models representative of accelerance, which then must be inverted to be primally assembled, or can be directly applied to state-space models representative of inverted FRFs. Afterwards, the redundant states can be erased by using the post-processing procedures presented in \citep{RD_2021}.

\subsection{Analytical and Numerical determinations}\label{Analytical and Numerical determinations}

In an analytical or numerical context, the mass, stiffness and damping matrices of the mechanical systems under study are usually known. In this way, to compute a state-space model representative of the diagonal apparent mass terms of a given CE, we must start by retaining from its mass, stiffness and damping matrices the coefficients associated with the set of diagonal apparent mass terms to be identified. Then, by following the expressions given in \ref{inverted_ss_matrices}, a state-space model representative of the selected set of diagonal apparent mass terms of the CE is obtained. Finally, by inverting this model (see \ref{Inverting_state_space_models}) we are capable of computing the intended state-space model representative of the inverted diagonal apparent mass terms of the CE.

Note that, if the CE under analysis is massless, the feed-through matrix of the state-space model representative of the diagonal apparent mass terms of the CE will be null. Hence, the inversion of this state-space model will introduce ill-conditioned numerical problems. To overcome this issue, one may add to the $[D]$ matrix of this state-space model an identity matrix multiplied by a small residue. The value of this small residue should be as small as possible: its only intent is preventing ill-conditioning of the $[D]$ matrix. It is clear that adding this matrix is similar to add small virtual masses to the edges of the CE under study.

\subsection{Experimental determination}\label{Experimental determination}

In an experimental context, the CEs are usually tested with attached fixtures to its edges, because performing a direct excitation of an isolated CE is generally unfeasible. Thus, to estimate state-space models representative of the inverted diagonal apparent mass terms we will extend the IS approach into the state-space domain.

\begin{figure}[H]
\centering
    \includegraphics[width=0.4\textwidth]{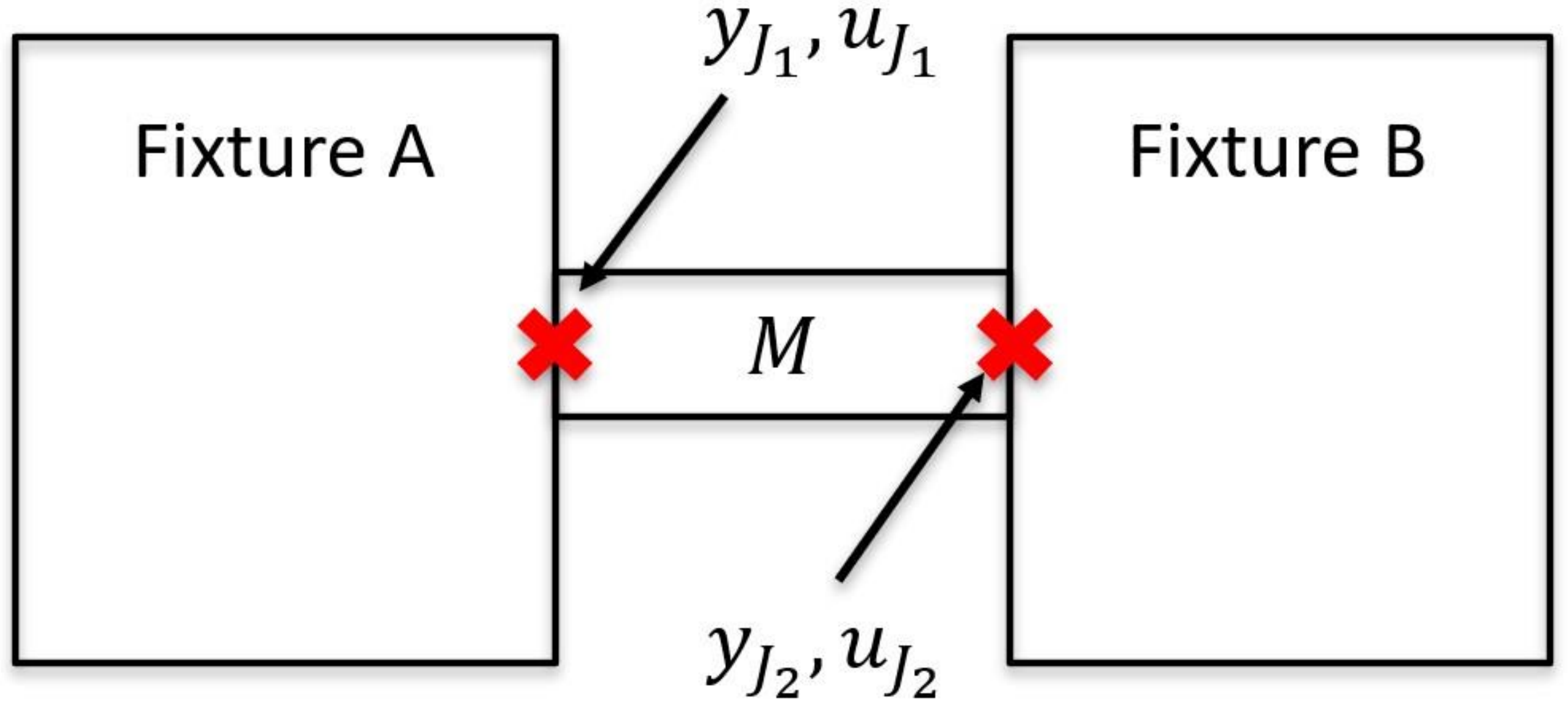}
    \caption{\textcolor{black}{Typical measurement set up to characterize the dynamic behaviour of a connecting element.}}
     \label{fig:Fixture_Mount_Fixture}
\end{figure}

To start, let us consider a generic set-up used to characterize the dynamic behaviour of CEs as represented in figure \ref{fig:Fixture_Mount_Fixture}.

Let us further assume that the inverted acceleration state-space models of the CE $M$ and of the fixtures $A$ and $B$ are available. If the primal formulation presented in section \ref{Primal formulation in time domain} is applied to assemble those state-space models, the assembled model would be given as follows:

\begin{equation}\label{eq:primal_assembled_ssm}
\begin{gathered}
\left\{\begin{matrix}
\dot{x}_{F_{A}}(t)\\
\dot{x}_{M}(t)\\
\dot{x}_{F_{B}}(t)
\end{matrix}
\right\}=\left[\begin{matrix}
A_{F_{A}} & 0 & 0\\
0 & A_{M} & 0\\
0 & 0 & A_{F_{B}}
\end{matrix}\right]\left\{\begin{matrix}
x_{F_{A}}(t)\\
x_{M}(t)\\
x_{F_{B}}(t)
\end{matrix}
\right\}+\left[\begin{matrix}
B_{F_{A},J_{1}} & 0\\
B_{M,J_{1}} & B_{M,J_{2}}\\
0 & B_{F_{B},J_{2}}
\end{matrix}\right]\left\{\begin{matrix}
\ddot{y}_{J_{1}}(t)\\
\ddot{y}_{J_{2}}(t)
\end{matrix}
\right\}\\
\left\{\begin{matrix}
u_{J_{1}}(t)\\
u_{J_{2}}(t)
\end{matrix}
\right\}=\left[\begin{matrix}
C_{F_{A},J_{1}} & C_{M,J_{1}} & 0\\
0 & C_{M,J_{2}} & C_{F_{B},J_{2}}
\end{matrix}\right]\left\{\begin{matrix}
x_{F_{A}}(t)\\
x_{M}(t)\\
x_{F_{B}}(t)
\end{matrix}
\right\}+\left[\begin{matrix}
D_{F_{A},J_{1}J_{1}}+D_{M,J_{1}J_{1}} & D_{M,J_{1}J_{2}} \\
D_{M,J_{2}J_{1}} & D_{F_{B},J_{2}J_{2}}+D_{M,J_{2}J_{2}} \\
\end{matrix}\right]\left\{\begin{matrix}
\ddot{y}_{J_{1}}(t)\\
\ddot{y}_{J_{2}}(t)
\end{matrix}
\right\}\\
\end{gathered}
\end{equation}

where, subscripts $F_{A}$ and $F_{B}$ denote variables related to the fixtures $A$ and $B$ (see figure \ref{fig:Fixture_Mount_Fixture}), respectively, while subscripts $J_{1}$ and $J_{2}$ denote inputs or/and outputs associated to the DOFs responsible for modeling the connection of the CE with fixtures $A$ and $B$, respectively. 

By applying a Fourier transformation, the state-space model given by expression \eqref{eq:primal_assembled_ssm} becomes:

\begin{equation}\label{eq:equationfirstformorderLaplace1}
\begin{gathered}
j\omega\{X_{F_{A}MF_{B}}(j\omega)\}=[A_{F_{A}MF_{B}}]\{X_{F_{A}MF_{B}}(j\omega)\} + [B_{F_{A}MF_{B}}](-\omega^{2}\{Y_{F_{A}MF_{B}}(j\omega)\})\\
\{U_{F_{A}MF_{B}}(j\omega)\}=[C_{F_{A}MF_{B}}]\{X_{F_{A}MF_{B}}(j\omega)\}+[D_{F_{A}MF_{B}}](-\omega^{2}\{Y_{F_{A}MF_{B}}(j\omega)\})
\end{gathered}
\end{equation}

where, subscript $F_{A}MF_{B}$ denote variables related to the assembly representative of the fixtures A and B linked by the connecting element M.

After some mathematical manipulations, the following expression to calculate the apparent mass of the mechanical system presented in figure \ref{fig:Fixture_Mount_Fixture} is obtained:

\begin{equation}\label{eq:equationfirstformorderLaplace}
\begin{gathered}
[Z^{A}(j\omega)]=\frac{\{U_{F_{A}MF_{B}}(j\omega)\}}{-\omega^{2}\{Y_{F_{A}MF_{B}}(j\omega)\}}=[C_{F_{A}MF_{B}}](j\omega[I]-[A_{F_{A}MF_{B}}])^{-1}[B_{F_{A}MF_{B}}]+[D_{F_{A}MF_{B}}]
\end{gathered}
\end{equation}

where, variable $Z^{A}$ represents an apparent mass matrix.

By using equation \eqref{eq:primal_assembled_ssm} and by performing some mathematical manipulations, equation \eqref{eq:equationfirstformorderLaplace} can be rewritten as follows:

\begin{equation}\label{eq:dynamic_mass}
\begin{gathered}
Z^{A}_{F_{A}MF_{B}}(j\omega)=\left[\begin{matrix}
Z^{A}_{F_{A}MF_{B},J_{1}J_{1}}(j\omega) & Z^{A}_{F_{A}MF_{B},J_{1}J_{2}}(j\omega)\\
Z^{A}_{F_{A}MF_{B},J_{2}J_{1}}(j\omega) & Z^{A}_{F_{A}MF_{B},J_{2}J_{2}}(j\omega)
\end{matrix}\right]
\end{gathered}
\end{equation}

where, matrices $[Z^{A}_{F_{A}MF_{B},J_{1}J_{1}}(j\omega)]$, $[Z^{A}_{F_{A}MF_{B},J_{1}J_{2}}(j\omega)]$, $[Z^{A}_{F_{A}MF_{B},J_{2}J_{1}}(j\omega)]$ and $[Z^{A}_{F_{A}MF_{B},J_{2}J_{2}}(j\omega)]$ are given by expression \eqref{eq:dynamic_mass_with_ssmt}.

\begin{equation}\label{eq:dynamic_mass_with_ssmt}
\begin{gathered}
[Z^{A}_{F_{A}MF_{B},J_{1}J_{1}}(j\omega)]=\left[\begin{matrix}
C_{F_{A},J_{1}} & C_{M,J_{1}} & 0
\end{matrix}\right](j\omega[I]-[A_{F_{A}MF_{B}}])^{-1}\left[\begin{matrix}
B_{F_{A},J_{1}}\\
B_{M,J_{1}}\\
0
\end{matrix}\right]+D_{F_{A},J_{1}J_{1}}+D_{M,J_{1}J_{1}}\\ 
[Z^{A}_{F_{A}MF_{B},J_{1}J_{2}}(j\omega)]=\left[\begin{matrix}
C_{F_{A},J_{1}} & C_{M,J_{1}} & 0
\end{matrix}\right](j\omega[I]-[A_{F_{A}MF_{B}}])^{-1}\left[\begin{matrix}
0\\
B_{M,J_{2}}\\
B_{F_{B},J_{2}}
\end{matrix}\right]+D_{M,J_{1}J_{2}}\\ 
[Z^{A}_{F_{A}MF_{B},J_{2}J_{1}}(j\omega)]=\left[\begin{matrix}
0 & C_{M,J_{2}} & C_{F_{B},J_{2}}
\end{matrix}\right](j\omega[I]-[A_{F_{A}MF_{B}}])^{-1}\left[\begin{matrix}
B_{F_{A},J_{1}}\\
B_{M,J_{1}}\\
0
\end{matrix}\right]+D_{M,J_{2}J_{1}}\\
[Z^{A}_{F_{A}MF_{B},J_{2}J_{2}}(j\omega)]=\left[\begin{matrix}
0 & C_{M,J_{2}} & C_{F_{B},J_{2}}
\end{matrix}\right](j\omega[I]-[A_{F_{A}MF_{B}}])^{-1}\left[\begin{matrix}
0\\
B_{M,J_{2}}\\
B_{F_{B},J_{2}}
\end{matrix}\right]+D_{F_{B},J_{2}J_{2}}+D_{M,J_{2}J_{2}}
\end{gathered}
\end{equation}

By analyzing expressions in \eqref{eq:dynamic_mass_with_ssmt} and remembering that $[A_{F_{A}MF_{B}}]$ is a block diagonal matrix, we conclude that the \textcolor{black}{off} diagonal apparent mass terms of the assembled structure can be retrieved from the state-space model of the CE. Therefore, those terms are also the \textcolor{black}{off} diagonal apparent mass terms of the CE as suggested in \citep{JM_2015}. To obtain a state-space model representative of one of those terms, we must retain from the assembled state-space model the inputs and outputs of the \textcolor{black}{off} diagonal term to be identified, while the others must be eliminated. This means that the columns of matrices $[B]$ and $[D]$ that are not being multiplied by the inputs of the \textcolor{black}{off} diagonal apparent mass term to be identified must be eliminated. Similarly, the rows of matrices $[C]$ and $[D]$ that do not correspond to the outputs of the \textcolor{black}{off} diagonal apparent mass term to be identified must also be eliminated.

By performing the identification of state-space models representative of the \textcolor{black}{off} diagonal apparent mass terms as described above, we will be including states from the fixtures $A$ and $B$, which are in theory not necessary to describe the terms intended to be identified. However, when the measured FRFs of the assembly used to test the CE are processed to estimate the respective state-space model, the states come mixed and not arranged as depicted in equation \eqref{eq:primal_assembled_ssm}. Consequently, the distinction between the states of the fixtures and of the CEs becomes hard to be performed. Hence, in general, during the construction of the state-space models representative of the \textcolor{black}{off} diagonal terms, all the states of the state-space model estimated from the measured FRFs of the assembled system are kept.

If inverse substructuring method is valid to characterize the CE under analysis, the diagonal apparent mass terms can be obtained by multiplying the \textcolor{black}{off} diagonal ones by $-1$ (see equation \eqref{eq:Z_IS_relation_1}). Hence, to obtain a state-space model representative of the diagonal apparent mass terms, the identified state-space model representative of the \textcolor{black}{off} diagonal terms must be transformed into negative form, which is performed by multiplying its output and feed-through matrices by $-1$ as demonstrated in \ref{Negative form of a state-space model representative of apparent mass}. Finally, by inverting the obtained state-space model, we may compute the intended model representative of the inverted diagonal apparent mass terms of the CE.

\color{black}

\section{\textcolor{black}{Numerical validation}}\label{Numerical Example}

\textcolor{black}{In this section, a numerical example will be used to validate the LM-SSS with compatibility relaxation developed in section \ref{LM_SSS_Non_Rigid_Formulation} and the state-space realization of IS presented in section \ref{Experimental determination}}. The section will start by describing the exploited numerical model (in section \ref{Description of the Numerical Model}). \textcolor{black}{Then, in section} \ref{Identification of State-Space Models}, the identification of the state-space models from the FRFs of the components and mounts contaminated with artificial noise is performed. Afterwards, in section \ref{Coupling by Relaxing the Compatibility Conditions} the coupling of the numerical substructures (see figure \ref{fig:Assembled_Numerical_structure_componnets_separated}) is performed by using \textcolor{black}{LM-SSS with compatibility relaxation}. Two coupled models are computed: one is obtained by coupling the untransformed models and the other by coupling the same models previously transformed into the coupling form presented in \citep{SJO_20072697} (from now on labelled as OCF, standing for Original Coupling Form). In section \ref{Coupling by Treating the Mounts as a Regular Component}, the components are coupled by assuming the mounts as another substructure to be coupled. To obtain the state-space model of the mounts primal state-space disassembly (see section \ref{Primal formulation in time domain}) is applied to decouple the attached fixtures, then the identified state-space models are coupled by using LM-SSS formulation (see \citep{RD_2021}). The FRFs of the coupled model are then compared with the ones of the coupled state-space model obtained from \textcolor{black}{LM-SSS with compatibility relaxation}. \textcolor{black}{Finally, in section \ref{Results and discussion} a discussion on the obtained results is performed.}

\subsection{Numerical test case}\label{Description of the Numerical Model} 
\begin{figure}[ht]
\centering
    \includegraphics[width=0.8\textwidth]{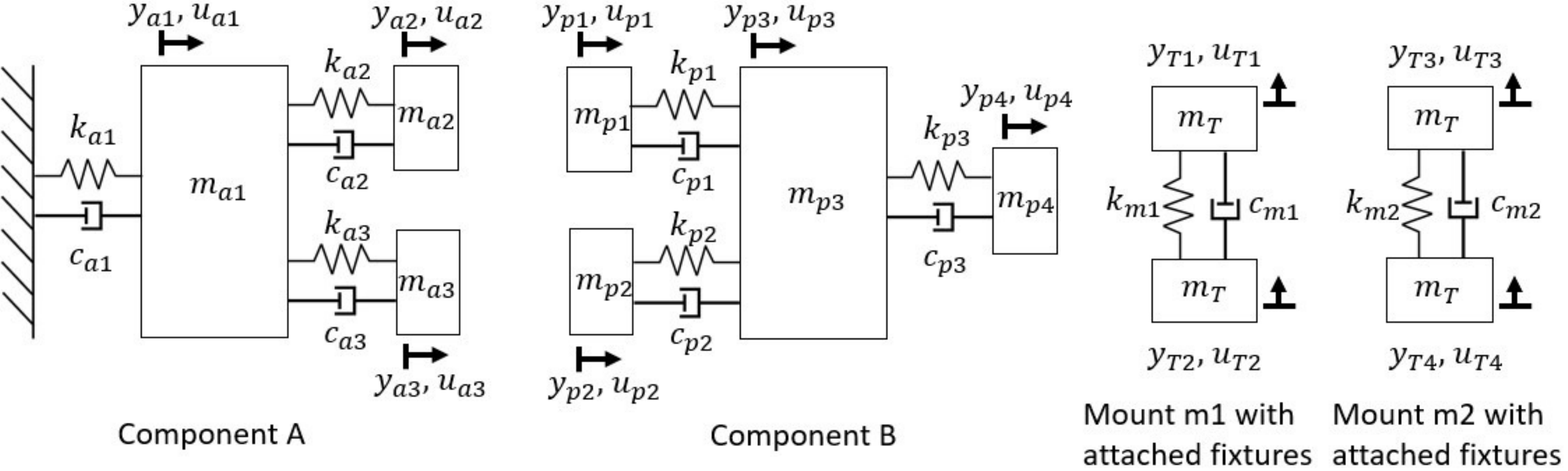}
    \caption{Isolated components.}
     \label{fig:Assembled_Numerical_structure_componnets_separated}
\end{figure}

\begin{figure}[ht]
\centering
    \includegraphics[width=0.6\textwidth]{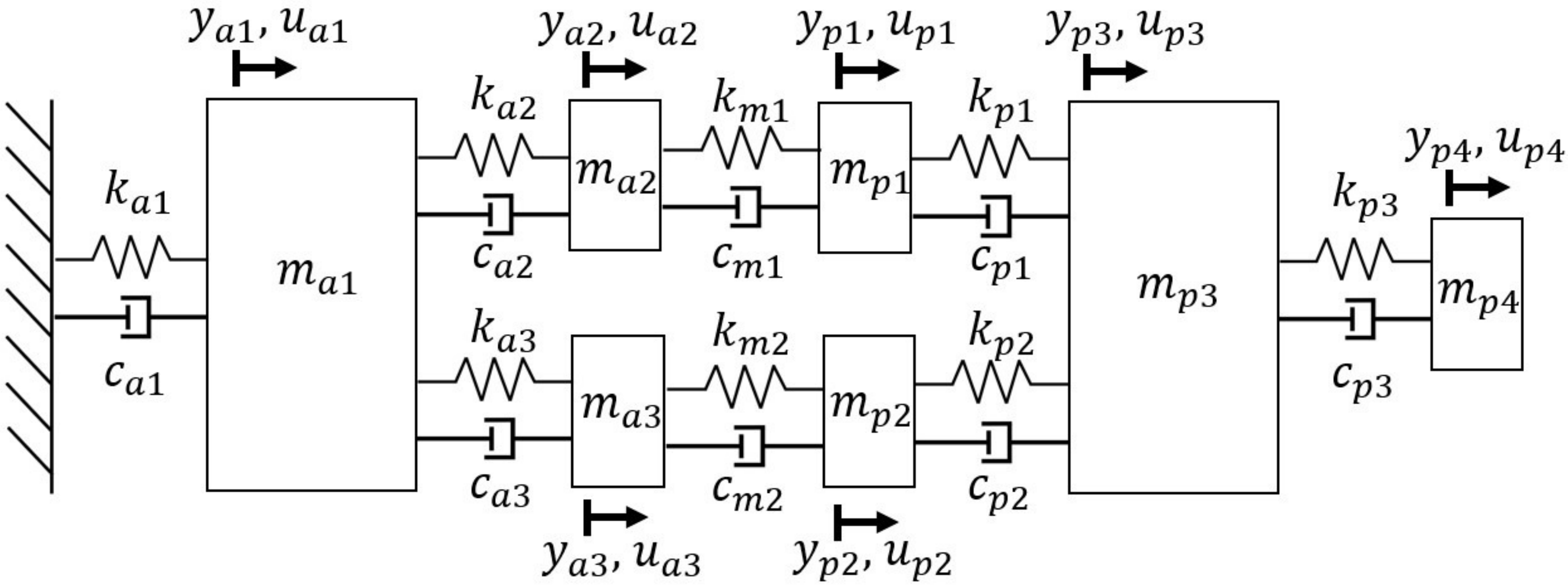}
    \caption{Assembled Structure.}
     \label{fig:Assembled_Numerical_structure}
\end{figure}

The numerical example here analyzed is composed by two components (see figure \ref{fig:Assembled_Numerical_structure_componnets_separated}) connected by two mounts (see figure \ref{fig:Assembled_Numerical_structure}). Table \ref{table:T1.1} reports the values assigned to the parameters addressed in figures \ref{fig:Assembled_Numerical_structure_componnets_separated} and \ref{fig:Assembled_Numerical_structure}. Those parameters are the mass, the damping and the stiffness of the substructures $A$ and $B$, of the mounts $m1$ and $m2$ and of the fixtures $T$.

\begin{table}[H]
\centering
\caption{Physical parameter values.}
\begin{tabular}{@{}llll@{}}
\toprule
$i$ \ & $m_{i}$ (kg) \ & $c_{i}$ (N s m$^{-1}$) \ & $k_{i}$ (N m$^{-1}$)  \\ \midrule
$a1$ & 10 & 30             & $1.5\times 10^{5}$\\
$a2$ & 3 & 50             & $5\times 10^{5}$\\
$a3$ & 3 & 50             & $4.5\times 10^{5}$\\
$p1$ & 5 & 50             & $1\times 10^{5}$\\
$p2$ & 7 & 50             & $1.5\times 10^{5}$\\
$p3$ & 10 & 10             & $5\times 10^{3}$\\
$p4$ & 1 & -             & -\\ 
$m1$ & - & 20 & $1 \times 10^{5}$\\ 
$m2$ & - & 20 & $2 \times 10^{5}$\\ 
$T$  & 2 & - & -\\ \midrule  
\label{table:T1.1}
\end{tabular}
\end{table}

\subsection{Identification of State-Space Models}\label{Identification of State-Space Models}

To compute the exact state-space models of the components, mounts and assembled structure, Lagrange equations \citep{HRH_1997} were firstly used to compute the mass and stiffness matrices of each mechanical system. The stiffness matrix of each component was used to built the respective damping one by replacing the stiffness terms with the damping ones. Finally, the exact accelerance state-space models were constructed from the computed matrices in compliance with \citep{FL_1988}.

Noise was also added to the FRFs of the components, in order to simulate realistic experimental data that would be affected by noise. This noise introduction followed the procedure used in \citep{SJO_20072697}, which consists in perturbing the real and imaginary parts of each FRF per frequency line by adding Gaussian distributed independent stochastic variables in accordance with equation \eqref{eq:perturbing_FRFs}:

\begin{equation}\label{eq:perturbing_FRFs}
H_{p,ij}(j\omega_{k})=H_{p,ij}(j\omega_{k})+\gamma_{ijk}+j\theta_{ijk}
\end{equation}

where, subscript $p$ denotes the structure to which a variable is associated. Subscripts $i$, $j$ and $k$ denote the output, input and the discrete frequency of the FRF term that is being perturbed. Finally, variable $j$ represents the imaginary unit, while variables $\gamma$ and $\theta$ are the Gaussian distributed variables responsible for perturbing the FRFs of the components, presenting zero mean and a standard deviation of $5 \times 10^{-3}$ $ms^{-2}N^{-1}$. \textcolor{black}{It is worth mentioning that by perturbing the FRFs as suggested in equation \eqref{eq:perturbing_FRFs} only the output part of the FRFs is contaminated with noise.} 

\begin{figure}[H]
\centering
    \includegraphics[width=1\textwidth]{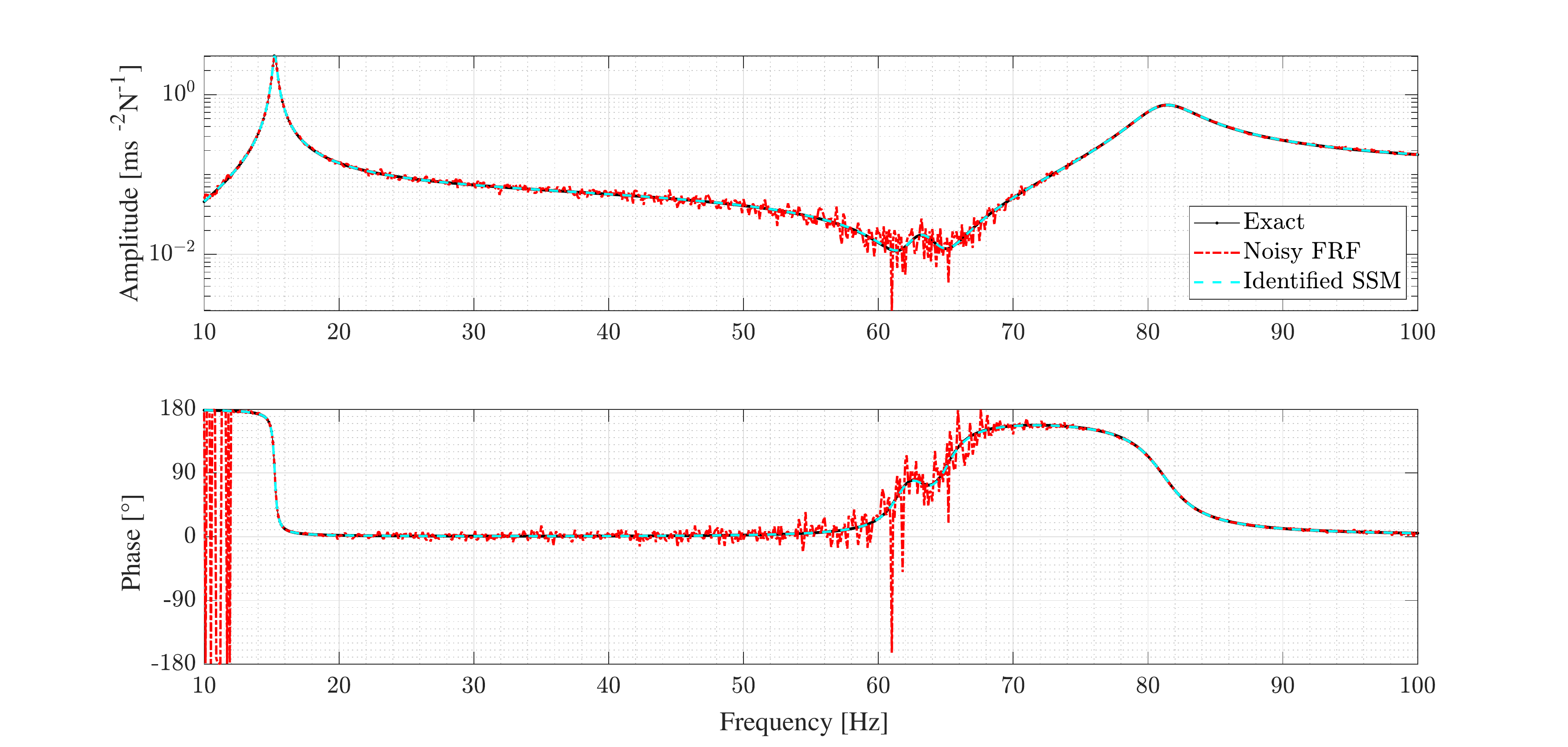}
    \caption{Accelerance FRF of the component $A$, whose output is the DOF $a1$ and the input is the DOF $a1$.}
     \label{fig:component_A_identified}
\end{figure}

\begin{figure}[H]
\centering
    \includegraphics[width=1\textwidth]{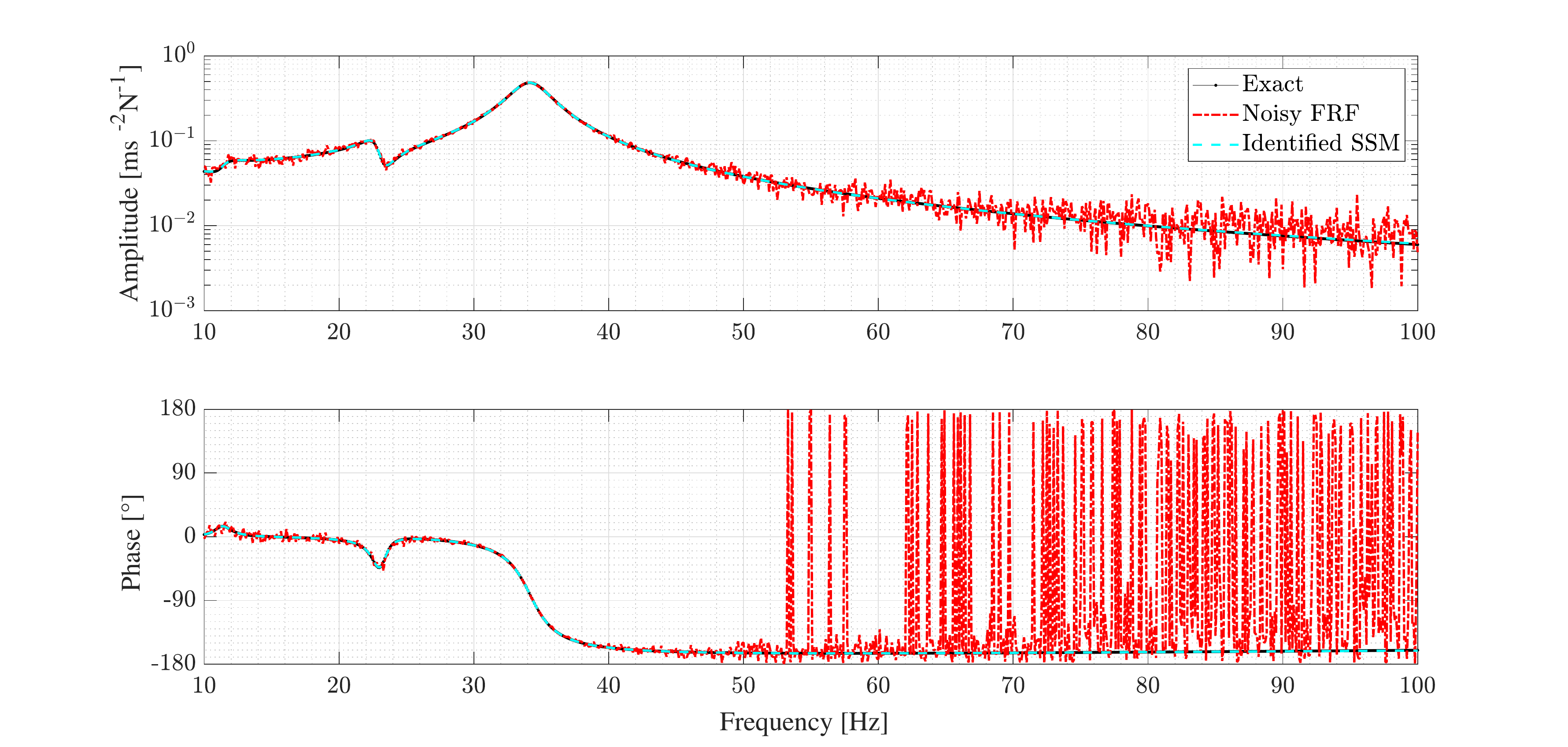}
    \caption{Accelerance FRF of the component $B$, whose output is the DOF $p1$ and the input is the DOF $p3$.}
     \label{fig:component_B_identified}
\end{figure}

\textcolor{black}{After the estimation of the FRFs, the Simcenter Testlab\textsuperscript{\textregistered} implementations of the PolyMAX (see, \cite{BPet_2004395}) and of the ML-MM (see, \cite{MEL_2015567}) methods were exploited to estimate modal parameters. It was found that the identified state-space models representative of components $A$ and $B$  were composed by 18 and 22 states, respectively. These models} were then double-differentiated (see \citep{FL_1988}) in order to compute the respective acceleration models. In figure \ref{fig:component_A_identified}, we observe the comparison of a noisy FRF of the component $A$ with the respective FRF of the exact and identified models, the same can be observed in figure \ref{fig:component_B_identified} for one FRF of component $B$.

After identifying the state-space models representative of the components, the same identification procedure was applied to estimate state-space models from the computed noisy FRFs representative of the mounts with attached fixtures (see figure \ref{fig:Assembled_Numerical_structure_componnets_separated}). \textcolor{black}{Each of the computed state-space models was found to be composed by 10 states.} Then, the obtained state-space model was double differentiated and inverted. To obtain a model representative of an \textcolor{black}{off} diagonal apparent mass term of the mount, the output and input associated to that term must be retained, while the others must be cancelled out (see section \textcolor{black}{\ref{Experimental determination}}). In this case, the output $y_{T1}$ and the input $u_{T2}$ were retained for the mount $m1$, while the output $y_{T3}$ and the input $u_{T4}$ were retained for the mount $m2$. After obtaining the state-space models representative of the \textcolor{black}{off} diagonal term of each mount, these models were transformed into negative form (see \ref{Negative form of a state-space model representative of apparent mass}) in order to obtain the state-space models representative of the diagonal apparent mass terms. Afterwards, those models were again inverted to obtain the state-space models required to perform coupling by using \textcolor{black}{LM-SSS with compatibility relaxation} (see section \ref{LM_SSS_Non_Rigid_Formulation}). \textcolor{black}{Moreover, by inverting the noisy FRFs of mounts $m1$, the noisy off diagonal apparent mass term, whose output is $y_{T1}$ and the input is $u_{T2}$ was retained. Then, by multiplying the retained off diagonal apparent mass term by -1 the noisy diagonal apparent mass term of mount $m1$ was obtained. The same procedures were applied with the noisy FRFs of mount $m2$ to identify the noisy off diagonal term, whose output is $y_{T3}$ and the input is $u_{T4}$. Afterwards, by multiplying the identified off diagonal term by -1 the noisy diagonal apparent mass term of mount $m2$ was computed.}

\begin{figure}[H]
\centering
    \includegraphics[width=1\textwidth]{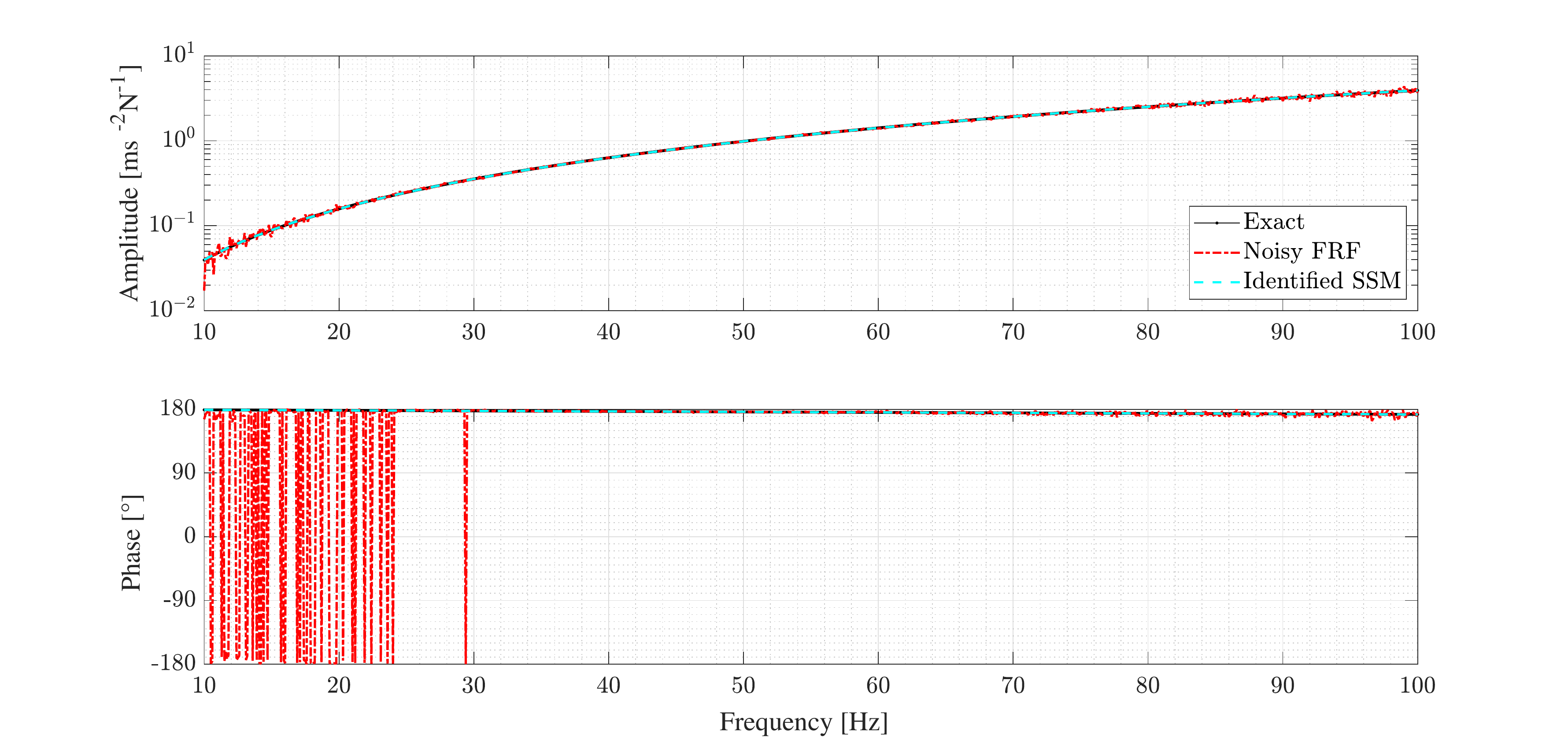}
    \caption{Inverted diagonal apparent mass term of mount $m1$.}
     \label{fig:component_M1_identified}
\end{figure}

\begin{figure}[H]
\centering
    \includegraphics[width=1\textwidth]{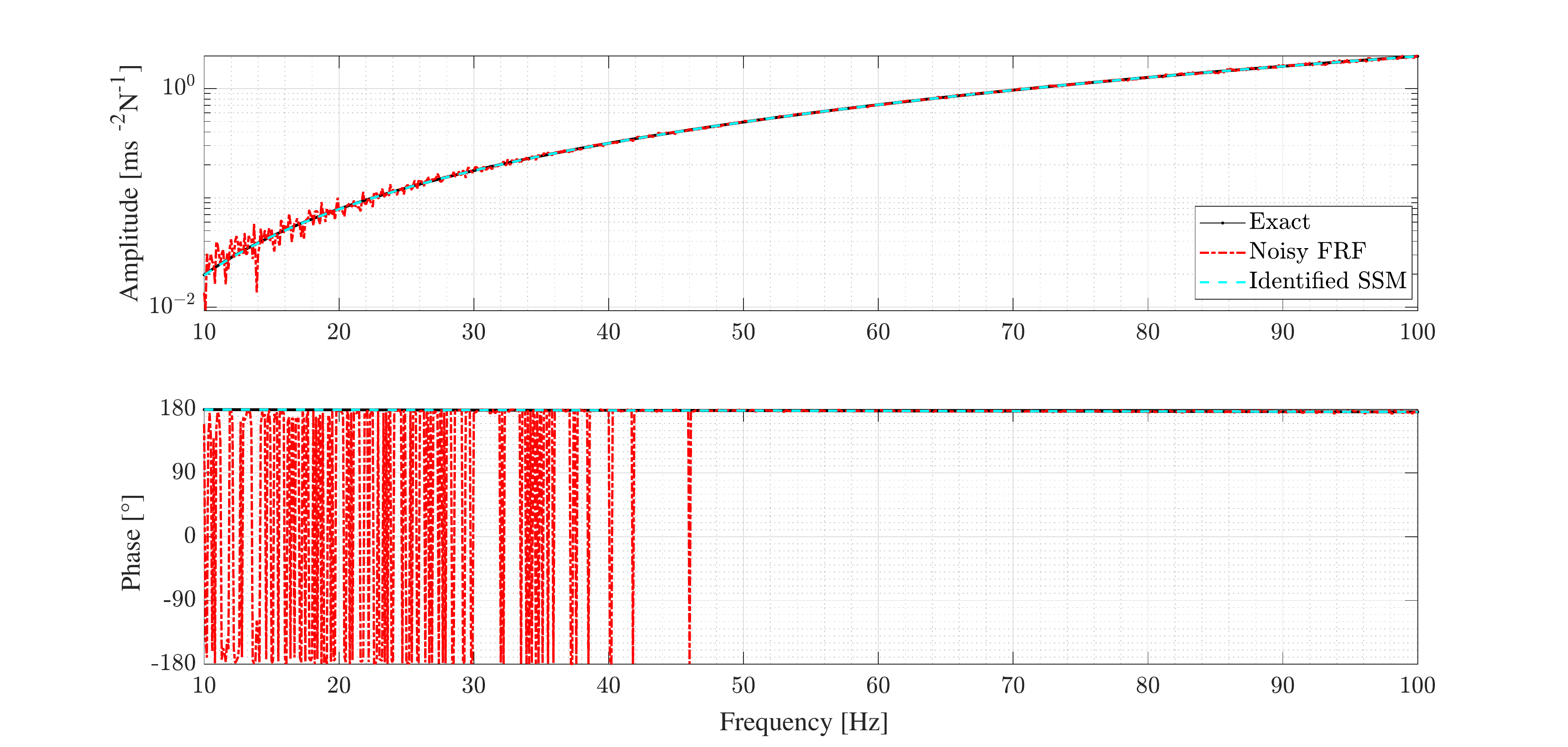}
    \caption{Inverted diagonal apparent mass term of mount $m2$.}
     \label{fig:component_M2_identified}
\end{figure}

Figure \ref{fig:component_M1_identified} reports the accelerance FRFs of the noisy inverted diagonal apparent mass term of mount m1 as well as those of the exact and identified SSM models. Figure \ref{fig:component_M2_identified} reports the same comparison for mount m2. 

\textcolor{black}{Analyzing figures \ref{fig:component_A_identified} and \ref{fig:component_B_identified}, it is clear that the system identification method has performed an accurate identification of both components. Moreover, observing figures \ref{fig:component_M1_identified} and \ref{fig:component_M2_identified} we may conclude that by exploiting the state-space realization of IS we are capable of computing reliable inverted state-space models representative of the diagonal apparent mass terms of the CEs from the models representative of the structures where they are included.} The estimated state space models were found to be not passive.

\subsection{Coupling by Relaxing the Compatibility Conditions}\label{Coupling by Relaxing the Compatibility Conditions}

To validate \textcolor{black}{LM-SSS with compatibility relaxation}, the state-space models estimated in section \ref{Identification of State-Space Models} were coupled. Then, to also validate the post-processing procedures to eliminate the extra states originated from coupling described in section \ref{Minimal-order coupled state-space models}, the estimated models were transformed into OCF and were coupled by using \textcolor{black}{LM-SSS with compatibility relaxation} to compute \textcolor{black}{a} minimal-order coupled model. \textcolor{black}{The state-space model obtained by using \textcolor{black}{LM-SSS with compatibility relaxation} was composed by 60 states, while the computed minimal-order coupled model was composed by 56 states.}

It is worth mentioning that the state-space models of the inverted diagonal apparent mass of the mounts written in OCF were obtained from the state-space model of the mount with attached fixtures, previously, transformed into coupling form by following the procedure described in section \ref{Identification of State-Space Models}. This choice is justified by the fact that to transform a state-space model into OCF its displacement output matrix is required. This matrix is available for the state-space model representative of the mount with attached fixtures (calculated by the system identification method). However, for the state-space model representative of the inverted diagonal apparent mass of the mounts this matrix is not available, being known just its acceleration output matrix. Since the computation of the displacement output matrix from the acceleration one demands the inversion of the state matrix (see \citep{FL_1988}), it is not recommended the direct transformation of the inverted diagonal apparent mass of the mounts into OCF.

\begin{figure}[H]
\centering
    \includegraphics[width=1\textwidth]{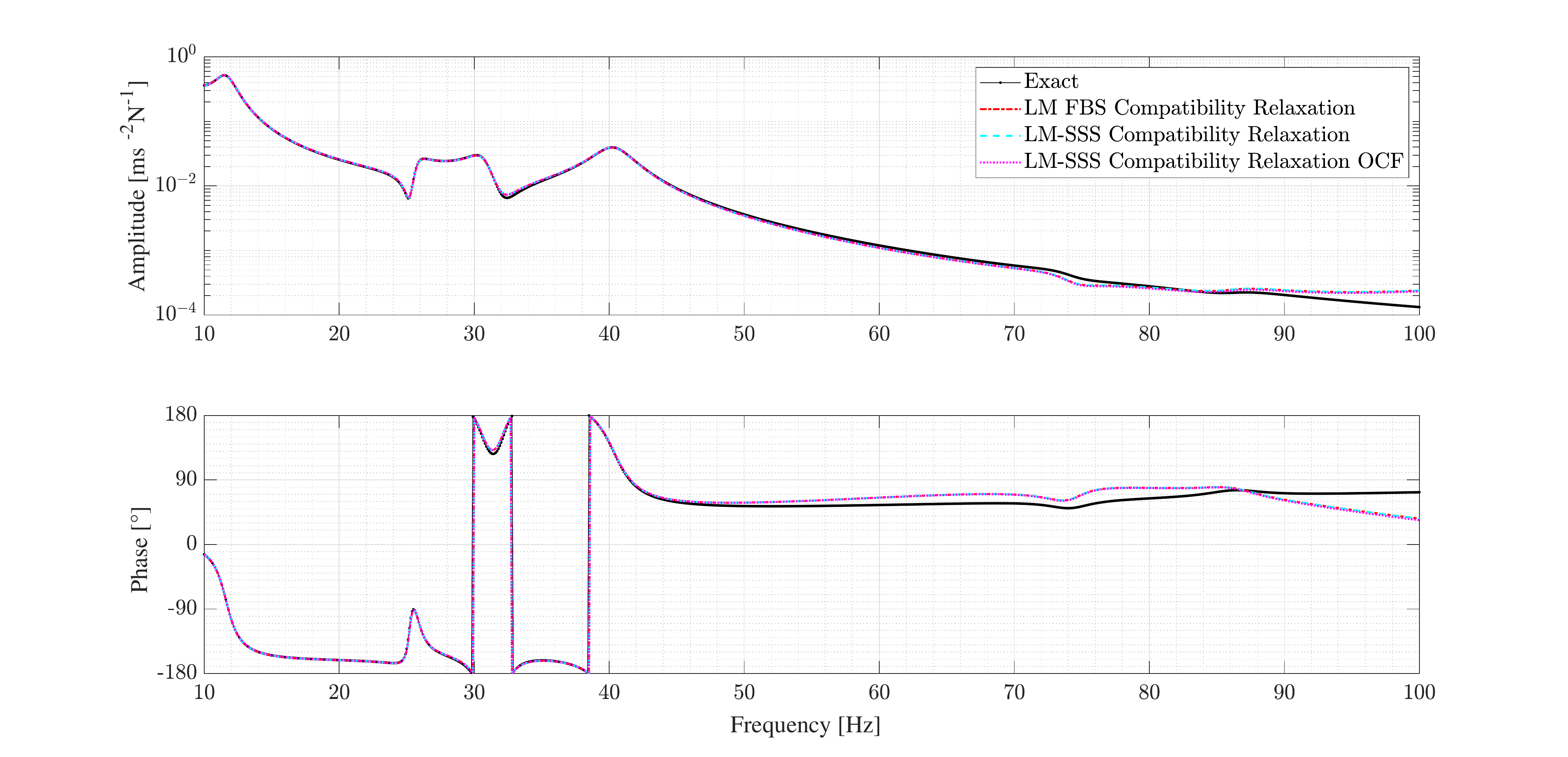}
    \caption{Accelerance FRF of the assembled structure, whose output is the DOF $p4$ and the input is the DOF $p1$.}
     \label{fig:coupled_results}
\end{figure}

Figure \ref{fig:coupled_results} shows the FRF of the exact coupled model as well as the FRF obtained by LM-FBS with compatibility relaxation \citep{EB_2014} and the FRFs of the coupled model obtained by coupling the estimated state space models untransformed and transformed into OCF.

By observing figure \ref{fig:coupled_results} we can conclude that the FRFs of the coupled state-space models obtained by \textcolor{black}{LM-SSS with compatibility relaxation} and the ones obtained by LM FBS with compatibility relaxation perfectly match. Furthermore, the FRF of the coupled model obtained by \textcolor{black}{LM-SSS with compatibility relaxation} is very well matching with the exact one. Hence, we may conclude that \textcolor{black}{LM-SSS with compatibility relaxation} is a reliable technique to couple substructures linked by CEs that are suitable to be characterized by IS. Moreover, from figure \ref{fig:coupled_results} we may conclude that the FRFs of the coupled state-space models obtained by using \textcolor{black}{LM-SSS with compatibility relaxation} and \textcolor{black}{LM-SSS with compatibility relaxation} OCF are very well matching, showing that the post-processing procedures outlined in section \ref{Minimal-order coupled state-space models} are valid to eliminate the extra states originated during coupling with \textcolor{black}{this method}. 



\subsection{Coupling by treating the Mounts as Regular Components}\label{Coupling by Treating the Mounts as a Regular Component}

\begin{figure}[H]
\centering
    \includegraphics[width=1\textwidth]{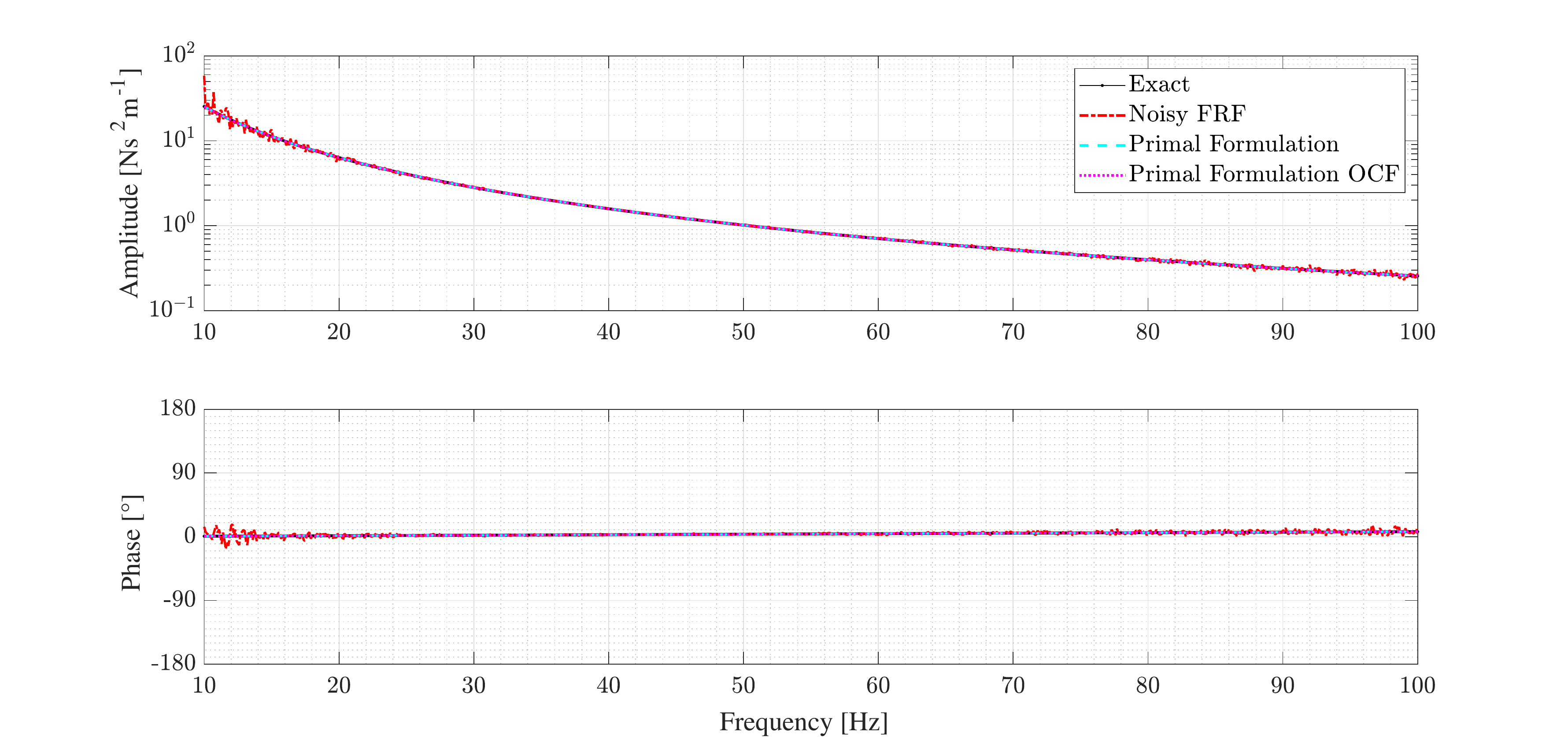}
    \caption{Apparent mass of the mount $m1$, whose output is the DOF $T1$ and the input is the DOF $T2$.}
     \label{fig:identified_mount_1}
\end{figure}

\begin{figure}[H]
\centering
    \includegraphics[width=1\textwidth]{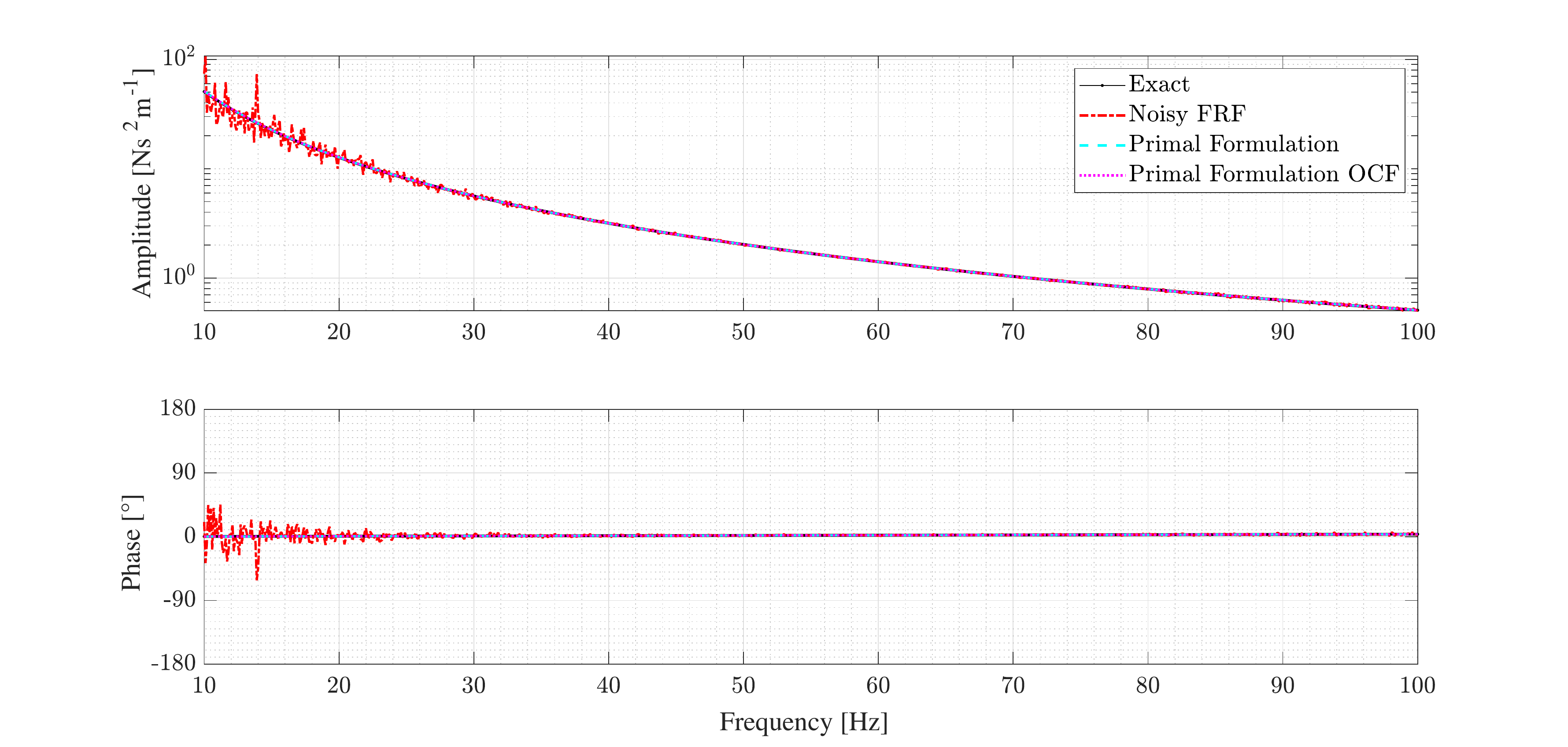}
    \caption{Apparent mass of the mount $m2$, whose output is the DOF $T3$ and the input is the DOF $T4$.}
     \label{fig:identified_mount_2}
\end{figure}

To validate the primal state-space formulation presented in section \ref{Primal formulation in time domain}, this technique was used to compute state-space models representative of mounts $m1$ and $m2$ by decoupling the attached fixtures at its edges (see figure \ref{fig:Assembled_Numerical_structure_componnets_separated}). Due to the lack of stiffness, the masses $T$ attached to each mount (see figure \ref{fig:Assembled_Numerical_structure_componnets_separated}) did not present any flexible mode. Consequently, it did not make sense to introduce their FRFs in the routines of the system identification algorithm to identify modes that actually do not exist. Hence, the state-space model of the attached masses was calculated directly from its mass matrix. Figures \ref{fig:identified_mount_1} and \ref{fig:identified_mount_2} show, respectively for mounts $m1$ and $m2$, the transfer (DOF $T2$ over DOF $T1$ for figure \ref{fig:identified_mount_1}; DOF $T3$ over DOF $T4$ for figure \ref{fig:identified_mount_2}) apparent mass function obtained for the following cases: i) exact model; ii) identification from the inverted noisy FRFs by using IS, primal disassembly iii) by using untransformed state-space models and iv) by using the state-space-models transformed into OCF \textcolor{black}{together with the post-processing presented in \cite{RD_2021} that relies on the use of a Boolean localization matrix to compute minimal-order coupled models.}

By observing figures \ref{fig:identified_mount_1} and \ref{fig:identified_mount_2} it is evident that the models of both mounts were successfully identified by using primal disassembly. Furthermore, the FRF of the identified model, obtained by using the models previously transformed into OCF, is also very well-matching with the exact solution. Hence, these results validate the primal formulation developed in section \ref{Primal formulation in time domain} and the hints given in the same section to obtain minimal-order models by using primal disassembly. \textcolor{black}{It was found that the state-space model obtained through primal disassembly was composed by 14 states, while the minimal-order model obtained by exploiting the same approach presented 10 states.}

\color{black}

By using LM-SSS, the inverted state-space models representative of the mounts M1 and M2 identified from the decoupling operations are coupled with substructures A and B. Furthermore, the same coupling is performed, but with the models previously transformed into OCF, together with the post-processing procedures presented in \cite{RD_2021} to compute a minimal-order coupled model. These coupled state-space models are composed by 60 and 52 states, respectively.

\begin{figure}[ht]
\centering
    \includegraphics[width=1\textwidth]{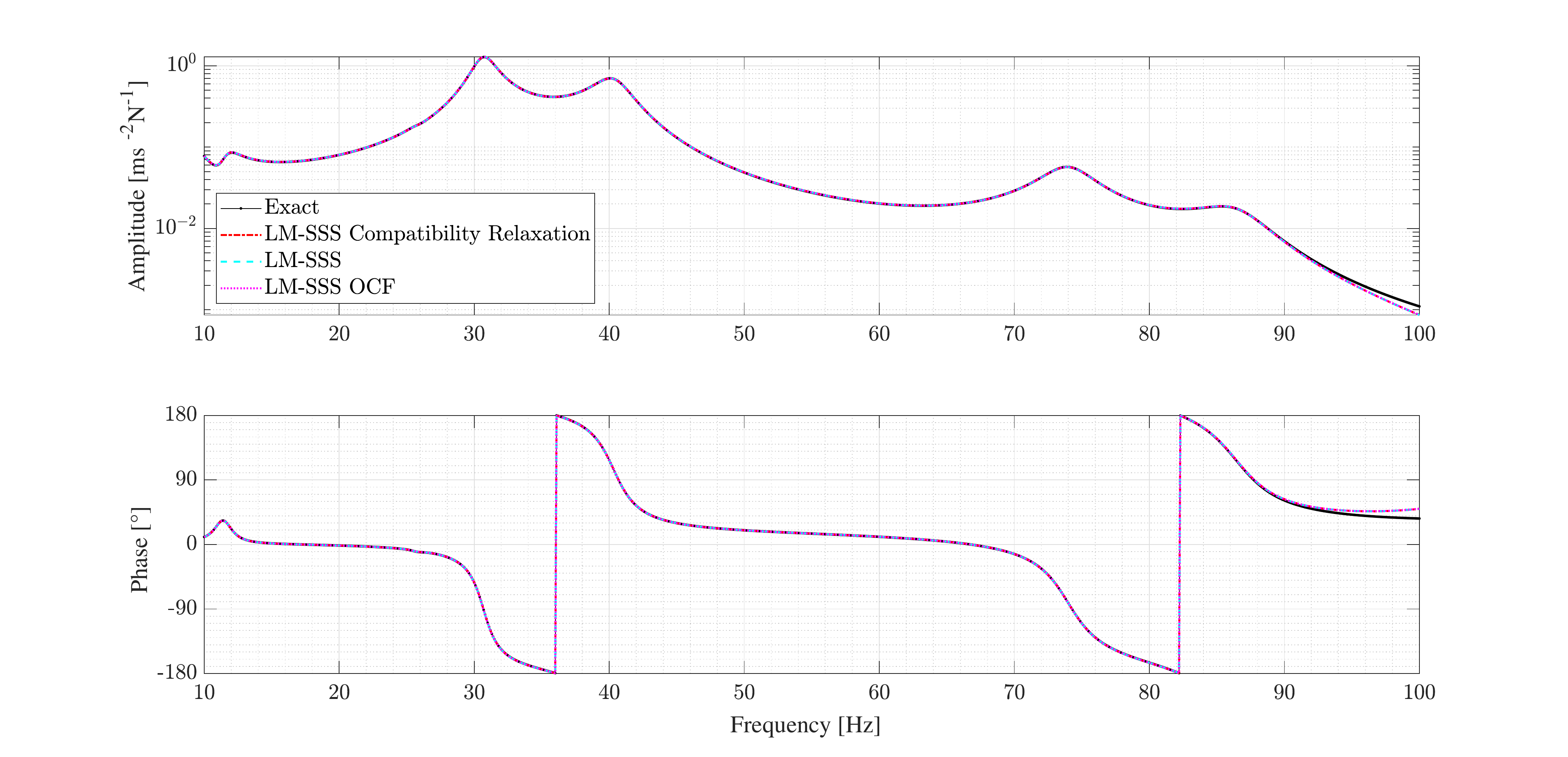}
    \caption{Accelerance FRF of the assembled structure, whose output is the DOF $p1$ and the input is the DOF $p2$.}
     \label{fig:coupled_mounts_with_LM_SSS}
\end{figure}

Figure \ref{fig:coupled_mounts_with_LM_SSS} shows the FRF of the theoretical model and those retrieved by the \textcolor{black}{LM-SSS with compatibility relaxation}, LM-SSS and LM-SSS OCF approaches. This comparison aimed at proving that all the techniques could lead to coupled state-space models presenting the same FRFs, provided that the connecting elements under study obey the IS assumptions.

\subsection{Results and discussion}\label{Results and discussion}

In sections \ref{Coupling by Relaxing the Compatibility Conditions} and \ref{Coupling by Treating the Mounts as a Regular Component}, it was found that the identified models representative of mounts M1 and M2 obtained by exploiting the state-space realization of the IS approach or by performing decoupling operations are reliable to characterize the dynamics of both mounts. Nonetheless, it was clear that the number of states presented by the models computed with primal disassembly was higher than the one computed through IS. This increment of states is due to the use of primal disassembly to decouple the masses attached to the edges of the mounts. 

However, analyzing the number of states of the minimal-order model obtained by exploiting primal disassembly together with tailored post-processing procedures, we conclude that this model is composed by the same number of states as the model obtained with IS. This is a direct consequence of the models representative of the masses be computed analytically and solely composed by interface DOFs (see section \ref{Coupling by Treating the Mounts as a Regular Component}). Analytically computed state-space models are composed by a number of states that is twice the number of DOFs, but as the models representative of the masses are just composed by interface DOFs, they are composed by a number of states which is twice the number of their interface DOFs. Therefore, as the post-processing procedures used to compute minimal-order models from primal assembly/disassembly operations eliminate a number of states that is twice the number of interface DOFs, all the states introduced through the performed decoupling operations are eliminated.  

Evaluating the number of states of the computed coupled models, we conclude that the models computed by \textcolor{black}{LM-SSS with compatibility relaxation} and LM-SSS are composed by 60 and 68 states, respectively. Once again, the state-space model obtained by using LM-SSS is composed by a higher number of states due to the inclusion of the mounts M1 and M2 through models identified by using decoupling operations. The minimal-order coupled state-space models originated with \textcolor{black}{LM-SSS with compatibility relaxation} OCF and LM-SSS OCF present 56 and 52 states, respectively. The larger number of states of the state-space model obtained by \textcolor{black}{LM-SSS with compatibility relaxation} OCF is justified by the fact that due to the inclusion of the mounts dynamics by an inverted model representative of one of their diagonal apparent mass terms, we are only enforcing compatibility by using half of their interface DOFs, hence from the coupled state-space models we are only able to eliminate a number of extra states equal to the number of the interface DOFs of the mounts. Whereas, by exploiting LM-SSS we are enforcing compatibility by using the full-set of interface DOFs of the CEs, thus we are able to eliminate a number of redundant states equal to twice the number of the interface DOFs of the mounts.

After analyzing the computed coupled models, we can conclude that the model computed by LM-SSS OCF is the most interesting, because it presents the same FRFs of the other coupled models and is composed by a lower number of states. Nevertheless, this happens only because we are dealing with an analytical example. In fact, if we were dealing with an experimental substructuring case, the state-space models of the masses attached to the mounts could not be obtained analytically. Therefore, they would have to be estimated by using system identification algorithms, which in general leads to the identification of state-space models composed by a number of states much larger than the number of outputs of the structures under analyze. As consequence, by performing decoupling operations we would contaminate the identified state-space models representative of the CEs with spurious states, making the coupled models computed with LM-SSS composed by a dramatically higher number of states. This is better illustrated in section \ref{Experimental Validation}, where state-space models estimated from experimentally acquired FRFs are coupled by exploiting \textcolor{black}{LM-SSS with compatibility relaxation}, primal state-space formulation and LM-SSS.

\section{Experimental validation}\label{Experimental Validation}

In this section, the approaches discussed in this paper are experimentally validated. We start by briefly presenting in section \ref{Testing Campaign} the experimentally characterized mechanical systems and how the testing campaign was conducted. In section \ref{Estimated components state-space models}, the FRFs collected from the performed experimental tests are exploited to compute state-space models representative of the components experimentally characterized. Then, in section \ref{Identified rubber mount state-space models} different state-space models representative of a CE are computed by performing decoupling with LM-SSS (see \cite{RD_MSSP_2022}), by using primal state-space disassembly (see section \ref{Primal formulation in time domain}) and by exploiting the state-space realization of IS (see section \ref{Experimental determination}) on the model of one of the assemblies composed by the CE to be characterized. Furthermore, two additional state-space models of the CE are computed by exploiting primal state-space disassembly together with post-processing procedures enabling the computation of a minimal-order model and by exploiting IS on the model of one of the assemblies composed by the CE to be characterized, previously transformed into OCF. Finally, in section \ref{Coupling Results} by using the identified CE models, different coupled state-space models are computed by exploiting LM-SSS, primal state-space assembly and \textcolor{black}{LM-SSS with compatibility relaxation} (see section \ref{LM_SSS_Non_Rigid_Formulation}). Afterwards, two minimal-order coupled models are also computed by exploiting both primal state-space assembly and \textcolor{black}{LM-SSS with compatibility relaxation} approaches together with tailored post-processing procedures. 


\subsection{Testing Campaign}\label{Testing Campaign}

The experimental data that will be exploited in this section was collected from the experimental modal characterization of the following substructures/structures:

\begin{itemize}
    \item Two aluminum crosses (from now on labelled as cross aluminum A and B);
    \item Two steel crosses (from now on labelled as cross steel A and B);
    \item Assembly A composed by the two aluminum crosses connected by a rubber mount;
    \item Assembly B composed by the two steel crosses connected by a rubber mount.
\end{itemize}

To perform dynamic substructuring operations, we are required to estimate state-space models, whose outputs and inputs are located at the interface of the analyzed components. However, placing sensors or/and exciting at the interface of the components under analyze is infeasible. Therefore, to overcome this issue, the aluminium and steel crosses were designed to behave as rigid bodies in the frequency band of interest (from 20 Hz to 500 Hz). Hence, we may excite and place sensors at accessible locations, then by using the well-known virtual point transformation (VPT) (see \cite{MV_2013}) approach we are capable of estimating the desired state-space models.  

\begin{figure}
\centering
\begin{subfigure}{0.5\textwidth}
  \centering
  \includegraphics[width=0.6\linewidth]{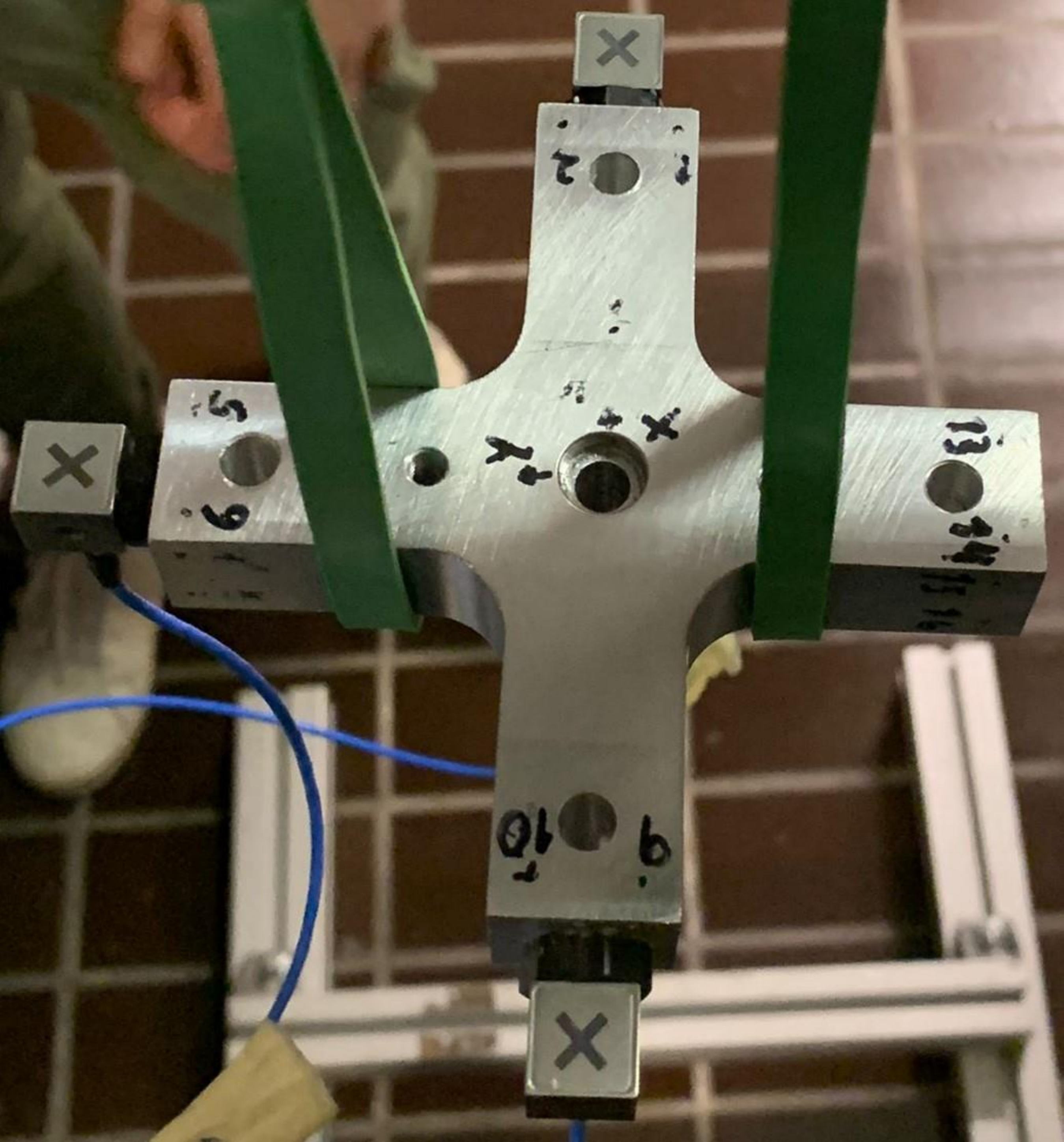}
  \caption{\textcolor{black}{Crosses.}}
  \label{fig:cross}
\end{subfigure}%
\begin{subfigure}{0.5\textwidth}
  \centering
  \includegraphics[width=0.6\linewidth]{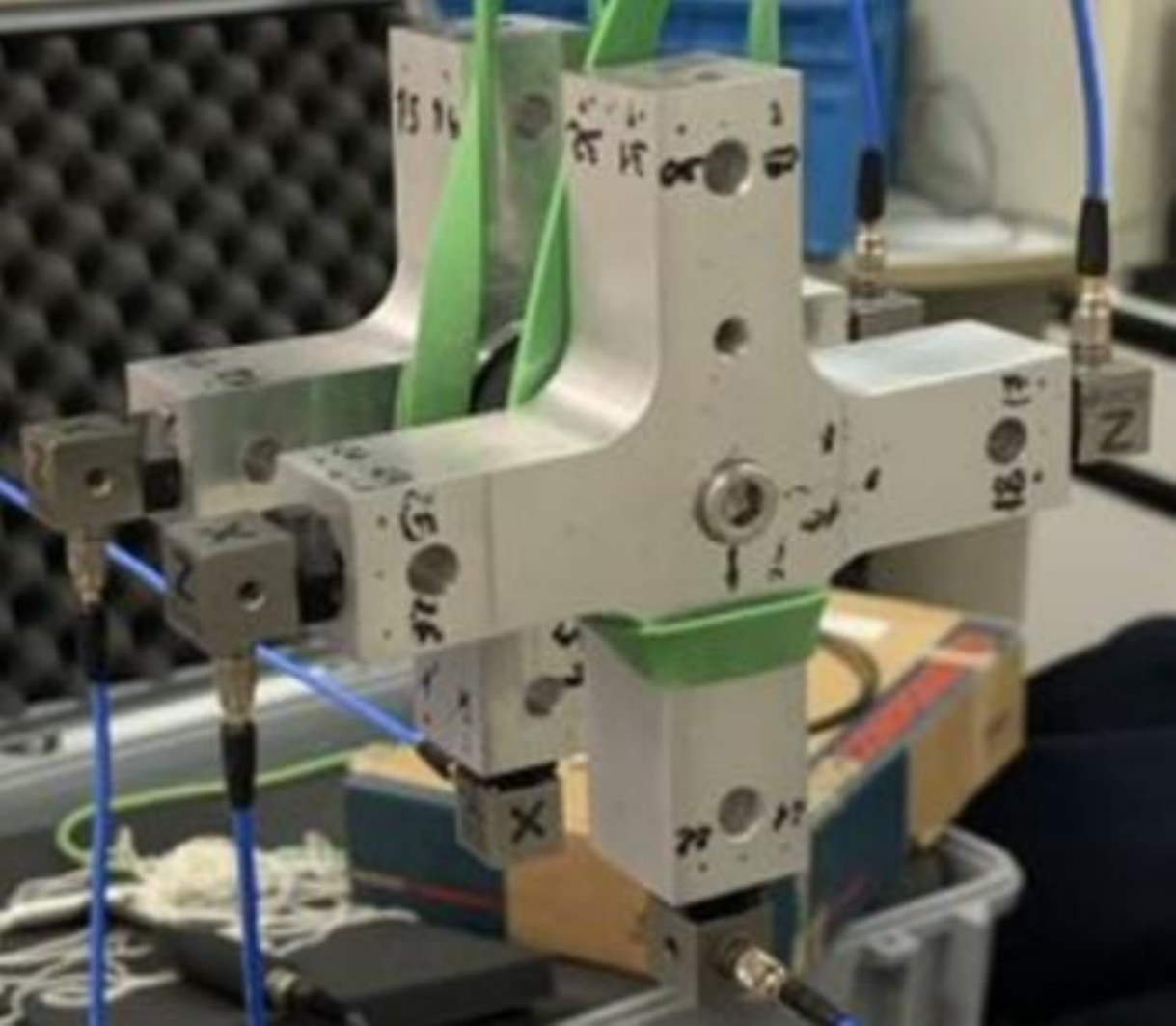}
  \caption{\textcolor{black}{Assemblies.}}
  \label{fig:Assembly}
\end{subfigure}
\caption{\textcolor{black}{Test set up used to experimentally characterize the isolated crosses and assemblies.}}
\label{fig:test}
\end{figure}

\begin{figure}[ht]
\centering
    \includegraphics[width=0.7\textwidth]{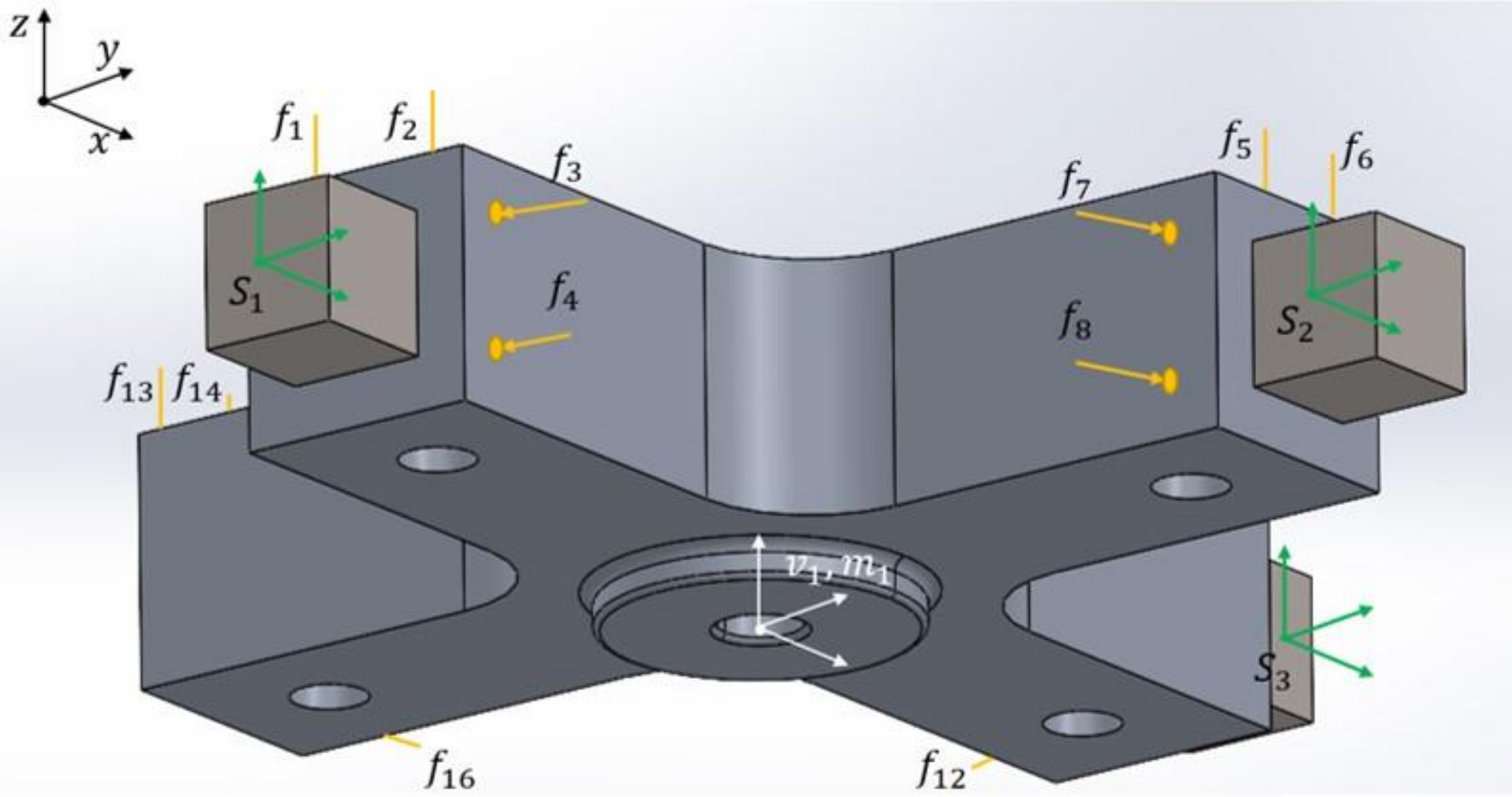}
    \caption{\textcolor{black}{Locations of measurement accelerometers (green), hammer impact directions (orange arrows) and virtual point (white).}}
     \label{fig:VPT_Cross}
\end{figure}

The test set-up exploited to characterize the isolated crosses and the assemblies can be observed in figures \ref{fig:cross} and \ref{fig:Assembly}. To perform the experimental tests the roving hammer testing approach was used. On each cross (either isolated or included in the assemblies), sixteen hits at different locations were applied and three accelerometers were placed as depicted in figure \ref{fig:VPT_Cross}. For further details on the performed testing campaign, the readers should refer to \cite{RD_MSSP_2022}.

\subsection{Estimated components state-space models}\label{Estimated components state-space models}

To estimate state-space models representative of the experimentally characterized components, we started by exploiting the Simcenter Testlab\textsuperscript{\textregistered} implementation of the PolyMAX and ML-MM methods to estimate modal parameters from the measured sets of FRFs. Then, the VPT-SS approach proposed in \cite{RD_ISMA_2022} (which, represents the extension of VPT into the state-space domain) was exploited to transform the outputs and the inputs of the estimated state-space models into the interface of the structures under study. It is worth mentioning that we could have followed a different approach to estimate the intended state-space models. We could have, firstly, applied the VPT approach on the measured FRFs to obtain the FRFs at the interface of the structures under study and then estimate the intended state-space models from these sets of FRFs. However, in \cite{RD_ISMA_2022} this approach demonstrated to lead to the estimation of less accurate state-space models composed by a larger amount of states. Furthermore, possible errors on the construction of the transformation matrices involved in the VPT approach can be easily compensated if the approach proposed in \cite{RD_ISMA_2022} is exploited.

\begin{figure}
  \begin{subfigure}[t]{1\textwidth}
    \centering
    \includegraphics[width=.8\textwidth]{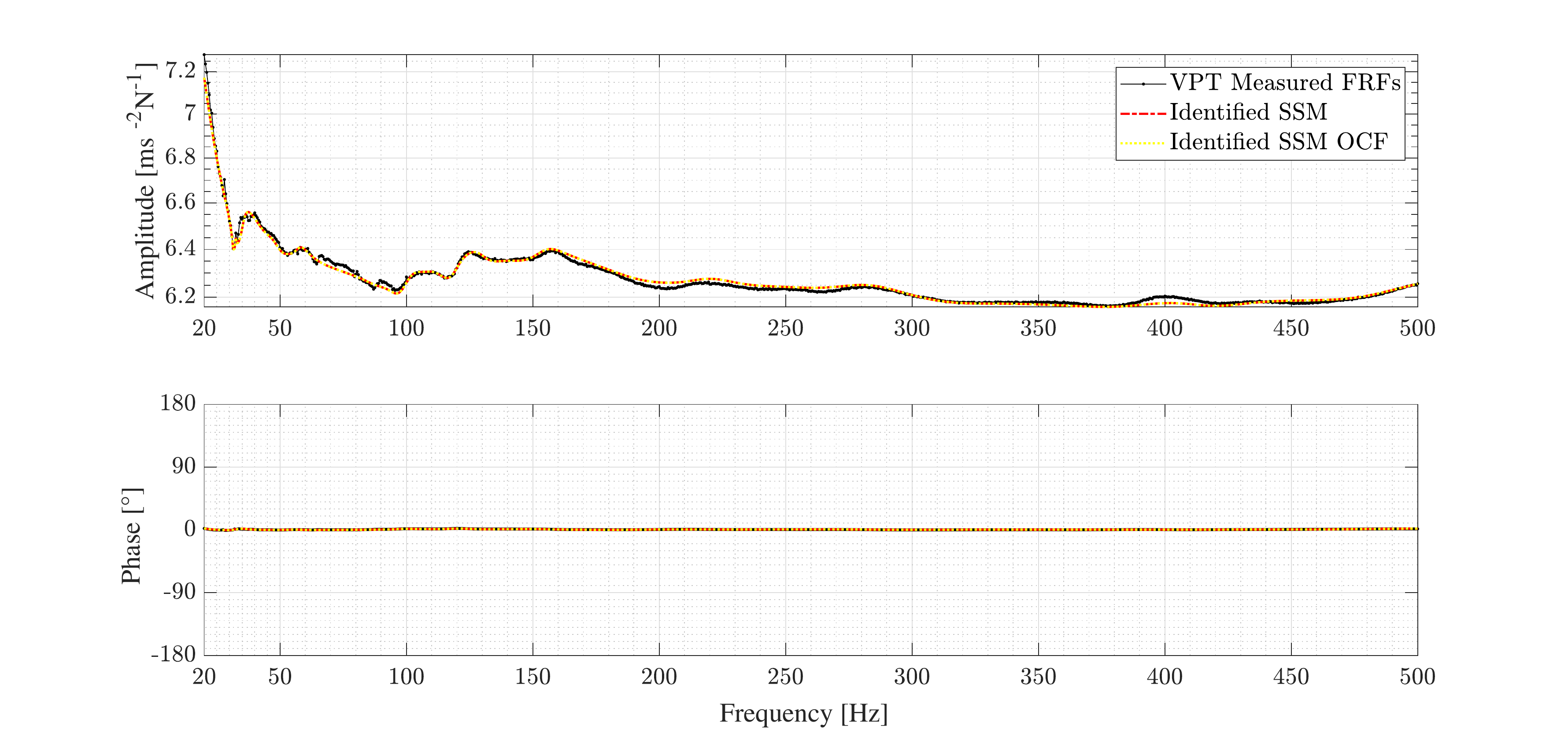}
    \caption{}
     \label{fig:Identified_VPT_Cross_Aluminum_A}
  \end{subfigure}
  \\
  \medskip
  \begin{subfigure}[t]{1\textwidth}
    \centering
    \includegraphics[width=.8\textwidth]{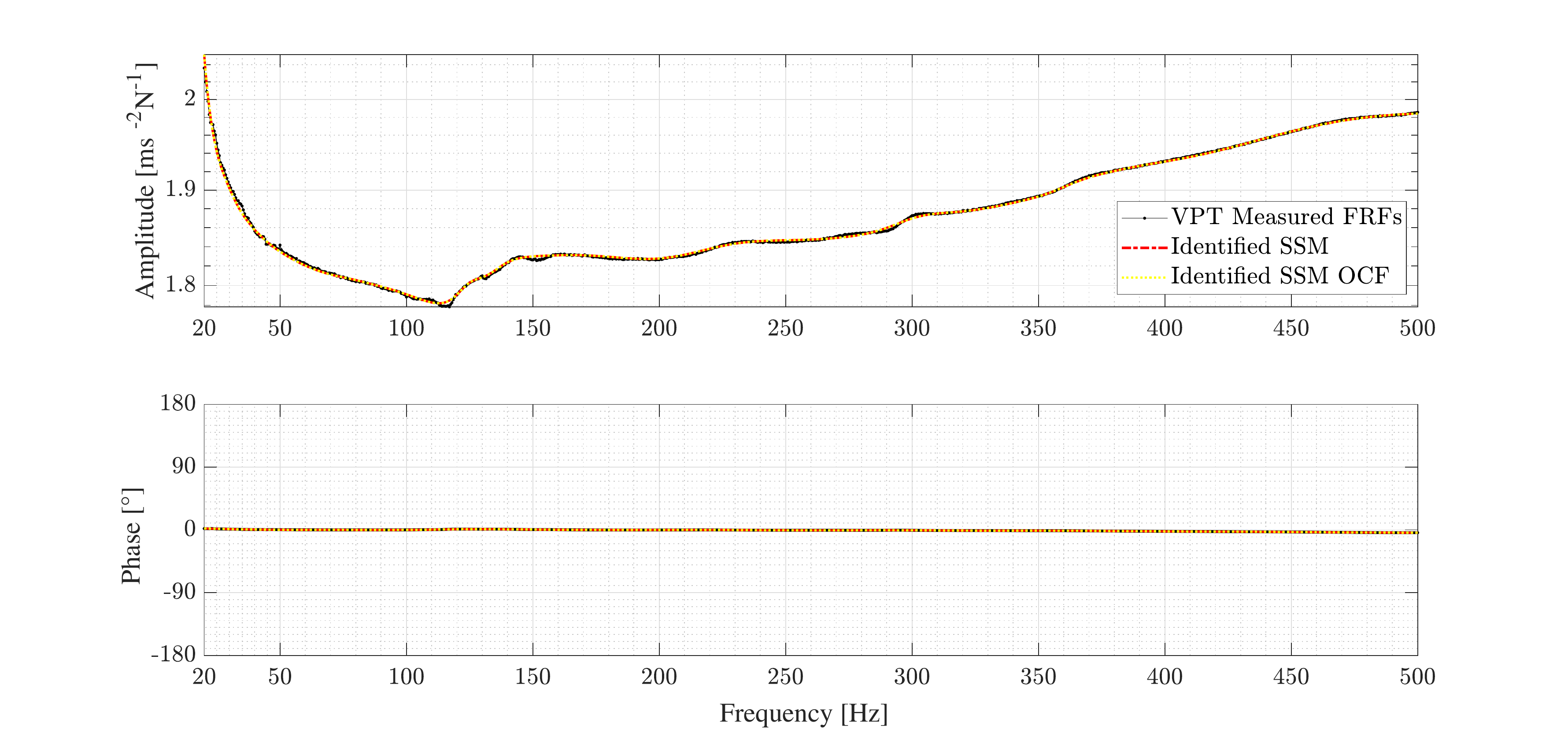}
    \caption{}
     \label{fig:Identified_VPT_Cross_Steel_B}
  \end{subfigure}
\caption{\textcolor{black}{Comparison, for each estimated model, of interface FRFs obtained by applying VPT on the measured FRFs with the FRFs of the estimated state-space model and with the FRFs of the same model previously transformed into OCF: a) FRF of the aluminum cross A, whose output is the DOF $v_{1}^{y}$ and the input is the DOF $m_{1}^{y}$; b) FRF of the steel cross B, whose output is the DOF $v_{2}^{z}$ and the input is the DOF $m_{2}^{z}$.}}
  \label{fig:Identified_VPT_State_Space_Models}
\end{figure}

\begin{figure}
  \begin{subfigure}[t]{1\textwidth}
    \centering
    \includegraphics[width=0.8\textwidth]{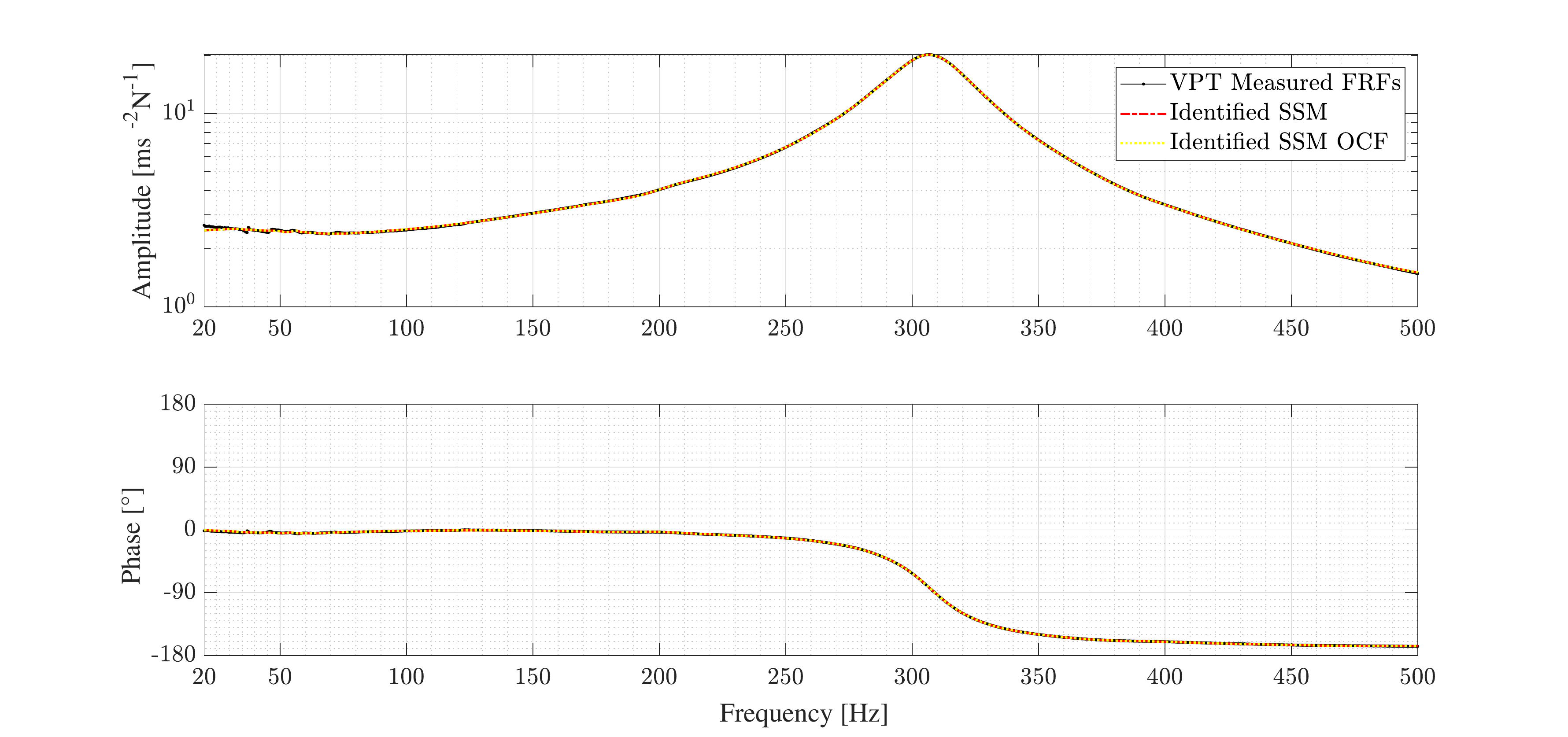}
    \caption{}
     \label{fig:Identified_VPT_Assembly_Aluminum_translational}
  \end{subfigure}
  \\
  \medskip
  \begin{subfigure}[t]{1\textwidth}
    \centering
    \includegraphics[width=0.8\textwidth]{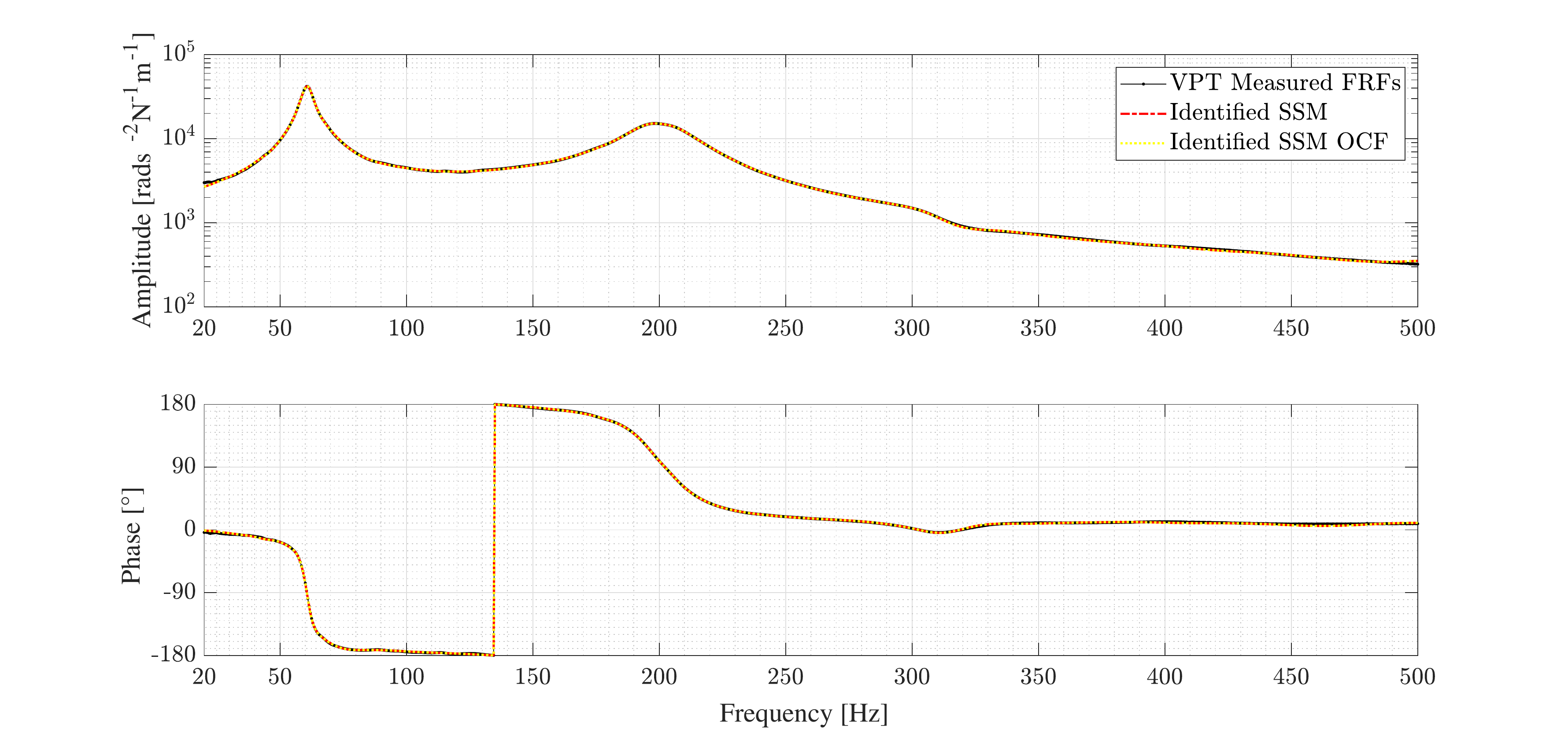}
    \caption{}
     \label{fig:Identified_VPT_Assembly_Aluminum_rotational}
  \end{subfigure}
  
\caption{\textcolor{black}{Comparison, for each estimated model, of interface FRFs obtained by applying VPT on the measured FRFs with the FRFs of the estimated state-space model and with the FRFs of the same model previously transformed into OCF: a) FRF of the assembly A, whose output is the DOF $v_{1}^{z}$ and the input is the DOF $m_{2}^{z}$; b) FRF of the assembly A, whose output is the DOF $v_{2}^{R_{x}}$ and the input is the DOF $m_{1}^{R_{x}}$.}}
  \label{fig:Identified_VPT_State_Space_Models_2}
\end{figure}

Figures \ref{fig:Identified_VPT_State_Space_Models} and \ref{fig:Identified_VPT_State_Space_Models_2} show, for each estimated model, the comparison of the interface FRFs obtained by applying VPT on the measured FRFs with the FRFs of the estimated state-space model and with the FRFs of the same model previously transformed into OCF. In the captions of figures \ref{fig:Identified_VPT_State_Space_Models} and \ref{fig:Identified_VPT_State_Space_Models_2} and through the rest of the paper, the interface outputs of the crosses aluminum and steel A will be tagged as $v_{1}$, while its interface inputs will be labelled as $m_{1}$. The interface outputs of the crosses aluminum and steel B will be tagged as $v_{2}$, whereas its interface inputs will be addressed as $m_{2}$. It is worth noticing that no state-space model representative of assembly B was estimated, because in section \ref{Coupling Results} we aim at computing coupled state-space models representative of this structure by using the state-space models here estimated. In the same section, the accuracy of these state-space models will be evaluated to assess the performance of the approaches used to compute them.

By observing figures \ref{fig:Identified_VPT_State_Space_Models} and \ref{fig:Identified_VPT_State_Space_Models_2}, it is possible to conclude that for all the estimated state-space models, the interface FRFs obtained by applying VPT on the measured FRFs are very well matching the FRFs of the estimated state-space model and the FRFs of the same model transformed into OCF. Therefore, we may conclude that the state-space models were successfully estimated, being suitable to be transformed into coupling form. 

It is worth mentioning that the computed state-space models representative of the aluminium cross A and B are composed by 195 and 201 states, respectively, the estimated models representative of the steel crosses A and B are composed by 197 and 183 states, respectively, and the state-space model of assembly A is composed by 242 states. It is worth mentioning that the estimated state-space models representative of the crosses aluminum A and B were expected to present the same number of states (the same was expected for the cross steel A and B). However, the FRFs of these components were experimentally acquired, as consequence they did not perfectly match, which lead to the estimation of state-space models composed by a different number of states. The mismatch between the measured FRFs of aluminum crosses A and B can be justified, for example due to slightly different boundary conditions during the experimental tests and due to small differences between the mechanical properties of both crosses (e.g. different mass) (the same justifications apply to explain the differences between the measured FRFs of the crosses steel A and B). Hence, the number of identified in-band modes was different, leading to the computation of state-space models composed by a different number of states.

\subsection{Identified rubber mount state-space models}\label{Identified rubber mount state-space models}

In this section, starting from the state-space model representative of assembly A, state-space models representative of the rubber mount included in assemblies A and B will be estimated by following three different approaches based on LM-SSS decoupling, primal state-space disassembly and exploiting IS on the estimated model representative of assembly A. To evaluate the accuracy of the mentioned approaches, the rubber mount FRFs obtained by exploiting LM FBS decoupling and IS on the interface FRFs of the assembly A will be taken as reference. Both approaches have already demonstrated to be valid to identify dynamic models of rubber isolators (see \cite{MH_2020}).

\begin{figure}
  \begin{subfigure}[t]{1\textwidth}
    \centering
    \includegraphics[width=.8\textwidth]{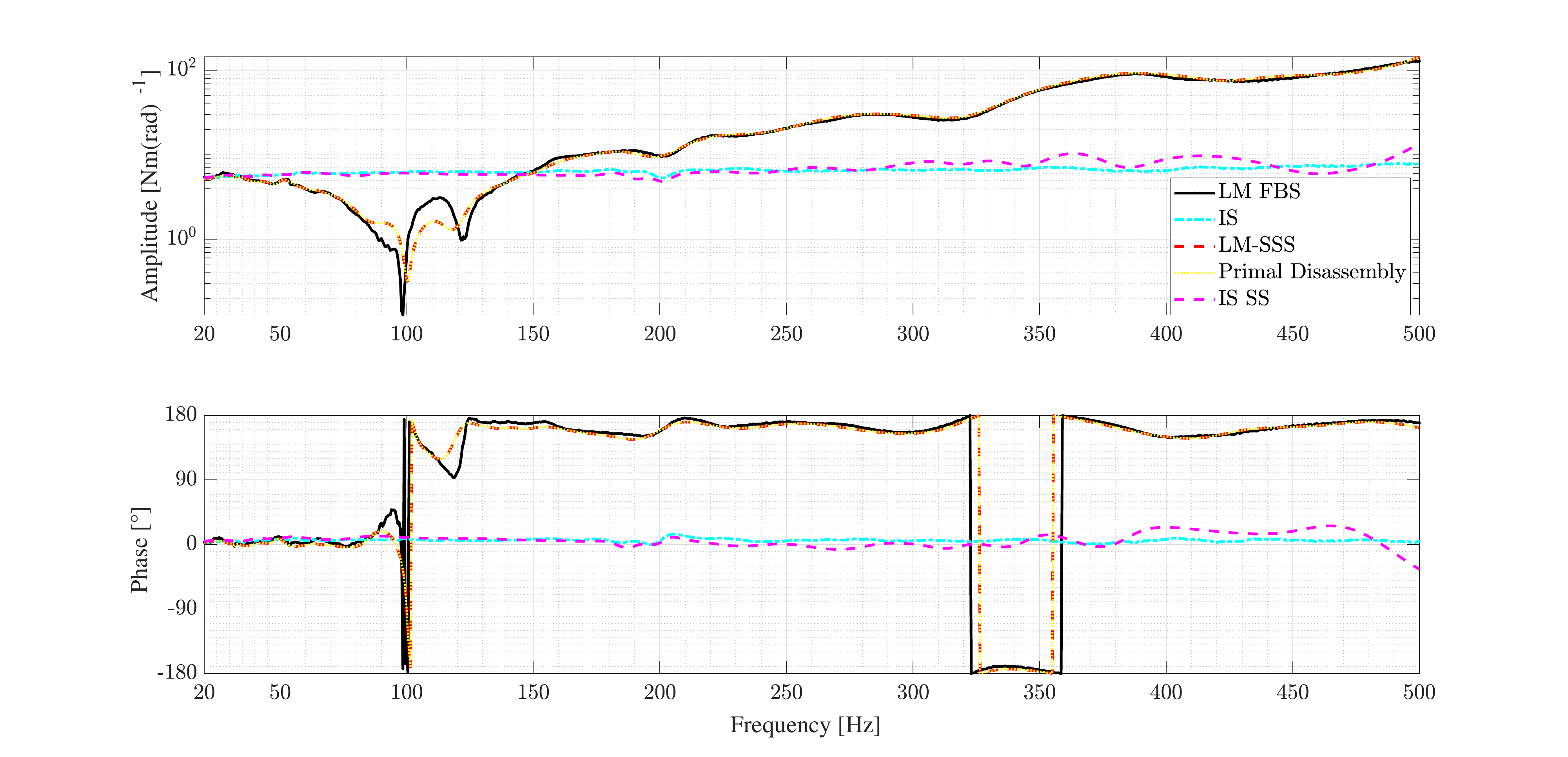}
   \caption{\textcolor{black}{Dynamic stiffness of the rubber mount, whose output is the DOF $v_{1}^{R_{z}}$ and the input is the DOF $m_{1}^{R_{z}}$.}}
     \label{fig:Identified_RM_SSMs_1}
  \end{subfigure}\\
  \medskip  
  \begin{subfigure}[t]{1\textwidth}
    \centering
    \includegraphics[width=.8\textwidth]{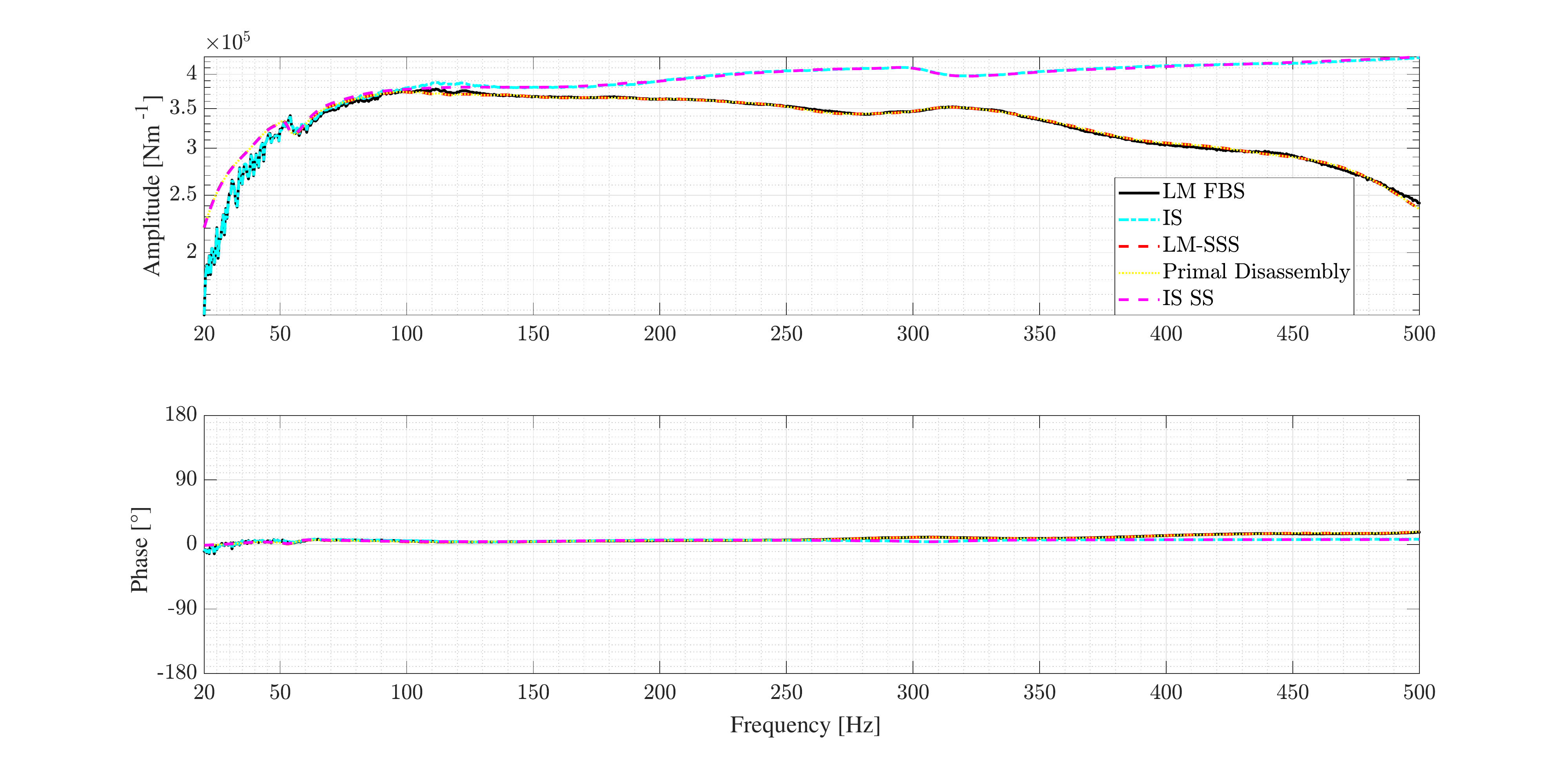}
    \caption{\textcolor{black}{Dynamic stiffness of the rubber mount, whose output is the DOF $v_{2}^{z}$ and the input is the DOF $m_{2}^{z}$.}}
     \label{fig:Identified_RM_SSMs_2}
  \end{subfigure}
  
\caption{\textcolor{black}{Comparison of some rubber mount FRFs obtained by applying the following approaches: i) LM FBS decoupling (represented by the black solid line), ii) IS on the apparent mass of assembly A (represented by the cyan dashdot line), iii) FRFs of the state-space model obtained with LM-SSS decoupling (represented by the red dashed line), iv) FRFs of the state-space model obtained with primal state-space disassembly (represented by the yellow dotted line) and v) IS on the inverted state-space model representative of assembly A (represented by the magenta dashed line).}}
  \label{fig:Identified_RM_SSMs}
\end{figure}

Figure \ref{fig:Identified_RM_SSMs} shows the comparison of a rubber mount driving point translation and rotation dynamic stiffness obtained by exploiting the following approaches: i) decoupling with LM FBS the accelerance interface FRFs of crosses aluminum A and B from the interface FRFs of assembly A, ii) by using IS approach on the apparent mass of assembly A to estimate the diagonal apparent mass terms of the rubber mount by multiplying by -1 the off diagonal apparent mass terms associated with the outputs placed on cross aluminum A and the inputs placed on cross aluminum B, iii) dynamic stiffness of the state-space model obtained by decoupling with LM-SSS the state-space models of the crosses aluminum A and B from the state-space model of assembly A, iv) dynamic stiffness of the model obtained by primarily disassemble the state-space models representative of the apparent mass of crosses aluminum A and B from the state-space model representative of the apparent mass of assembly A and v) dynamic stiffness of the model obtained by applying the IS approach on the inverted model of assembly A to identify a model representative of the diagonal apparent mass terms of the rubber mount by transforming into negative form (see \ref{Negative form of a state-space model representative of apparent mass}) the state-space model representative of its off diagonal terms, whose outputs are placed on cross aluminum A and the inputs are placed on cross aluminum B.

Note that, to ease the interpretation of figure \ref{fig:Identified_RM_SSMs}, we have decided to show the correspondent dynamic stiffness matrices of all the computed solutions. Hence, the accelerance FRFs directly obtained with approaches i and iii were inverted and multiplied by $-\omega^{2}$ to calculate the correspondent dynamic stiffness matrices. Whereas, the apparent mass directly computed by approaches ii, iv and v were multiplied by $-\omega^{2}$ to compute the correspondent dynamic stiffness.

By observing figure \ref{fig:Identified_RM_SSMs}, it is straightforward that as expected the dynamic stiffness of the rubber mount state-space models identified by LM-SSS decoupling and by primal state-space disassembly are perfectly matching. Hence, we may conclude that the primal state-space formulation developed in section \ref{Primal formulation in time domain} is accurate when dealing with state-space models estimated from measured data. Furthermore, both solutions are well-matching the dynamic stiffness identified with LM FBS decoupling, which demonstrates that the state-space models computed by exploiting LM-SSS decoupling and the primal state-space disassembly are reliable. 

By further analyzing figure \ref{fig:Identified_RM_SSMs}, it is evident that the dynamic stiffness matrices computed by exploiting IS on the measured FRFs of assembly A and on the estimated state-space model representative of assembly A are well-matching. Hence, we conclude that the procedures outlined in section \ref{LM_SSS_Non_Rigid_Formulation} to identify state-space models representative of off diagonal terms is valid even when dealing with state-space models estimated from measured FRFs. 

Comparing the dynamic stiffness matrices obtained from the three approaches relying on a decoupling operation (i.e. approaches i,iii and iv) with the two approaches relying on IS (i.e. approaches ii and v), we observe that at low frequencies the results provided by all the approaches are similar. However, as we analyze higher frequencies, a more pronounced bias between the approaches starts to occur. Indeed, the observed bias at higher frequencies is expected, because by exploiting IS we are assuming that the rubber mount is massless. Thus, the observed bias indicates that for higher frequencies this assumption does not hold and that the mass of the rubber mount starts to have a non negligible effect over its dynamics.

\begin{figure}
  \begin{subfigure}[t]{1\textwidth}
    \centering
    \includegraphics[width=.8\textwidth]{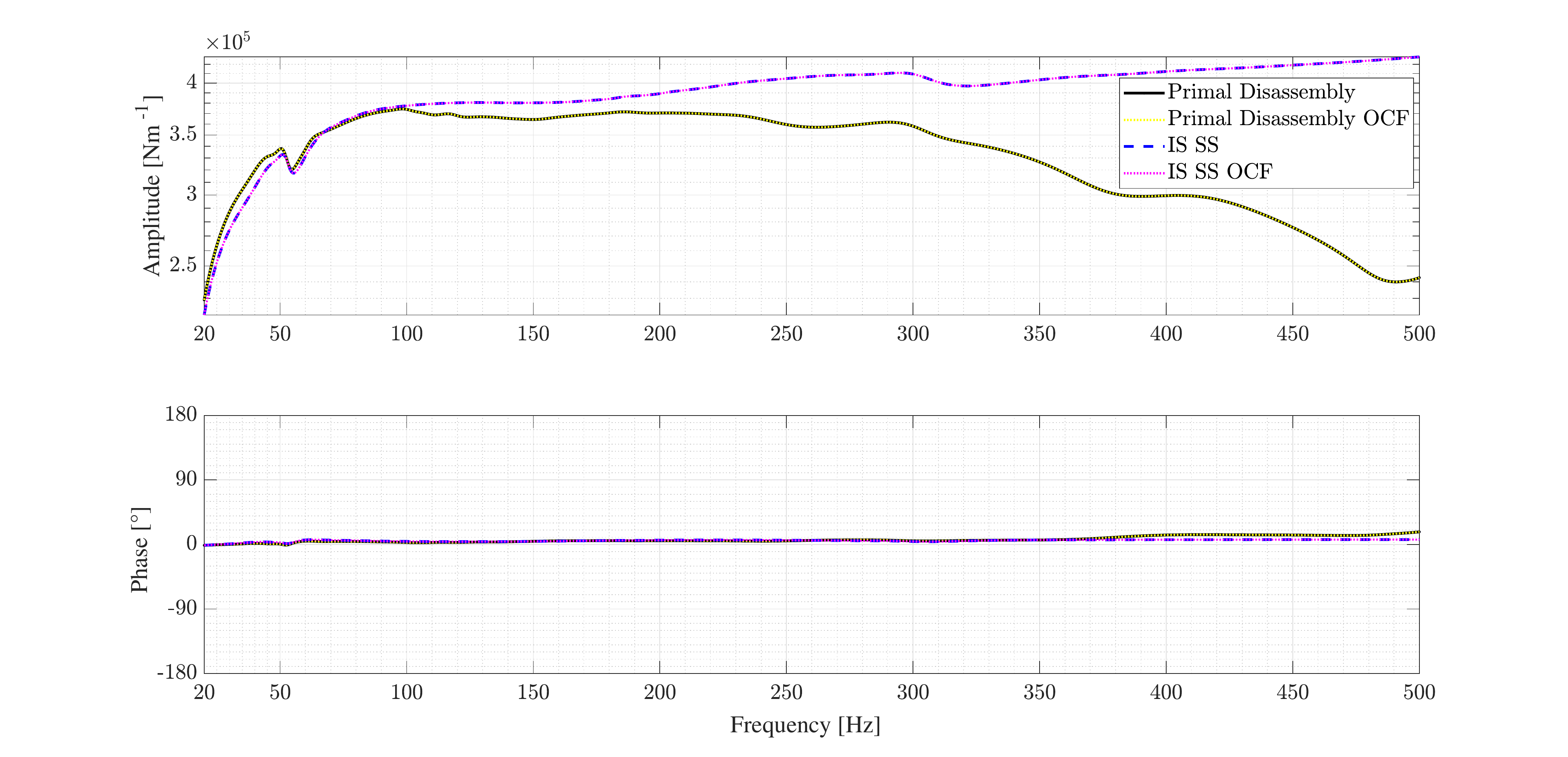}
   \caption{\textcolor{black}{Dynamic stiffness of the rubber mount, whose output is the DOF $v_{1}^{z}$ and the input is the DOF $m_{1}^{z}$.}}
     \label{fig:Identified_RM_SSMs_CF_1}
  \end{subfigure}\\
  \medskip
  \begin{subfigure}[t]{1\textwidth}
    \centering
    \includegraphics[width=.8\textwidth]{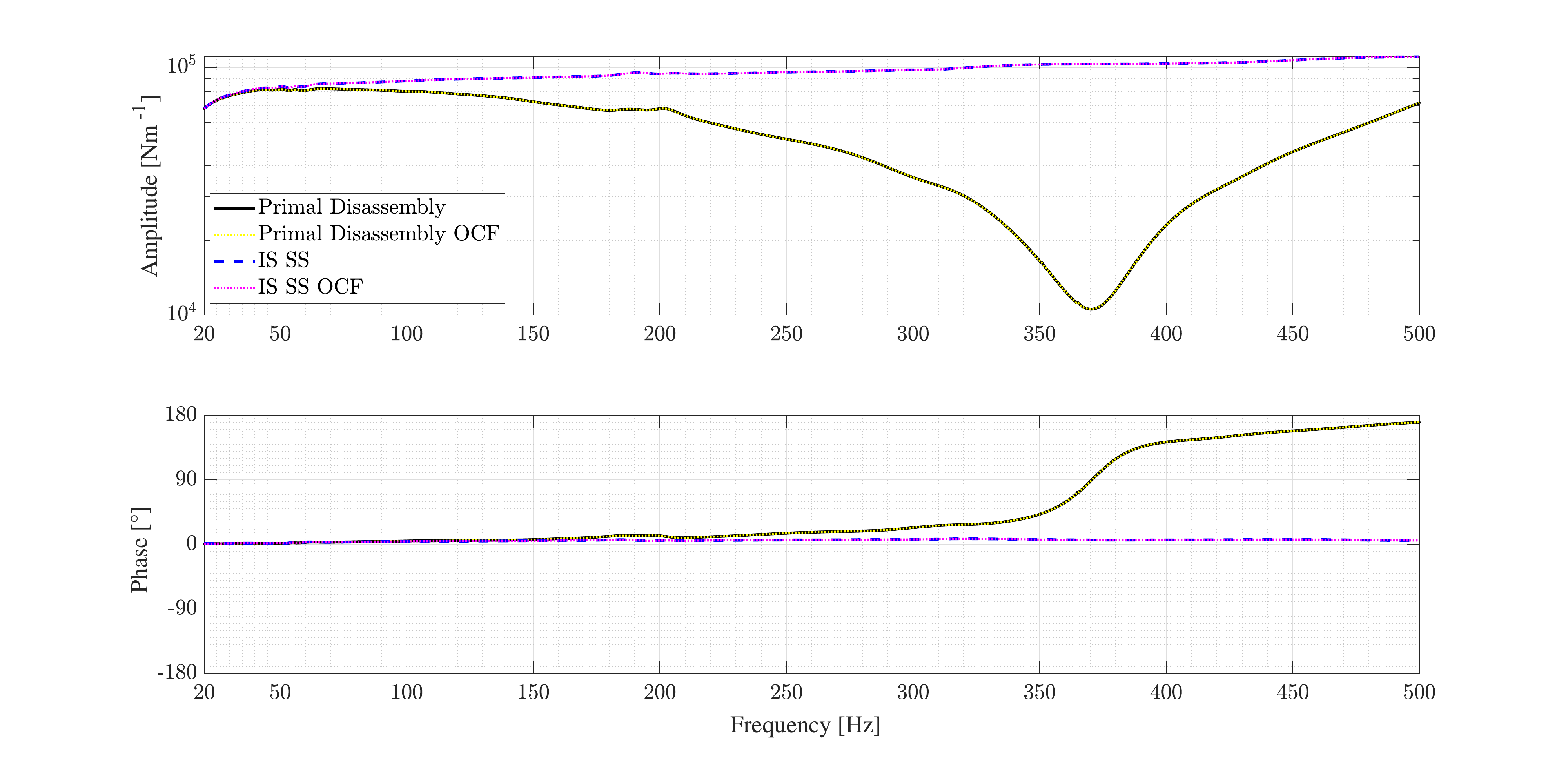}
    \caption{\textcolor{black}{Dynamic stiffness of the rubber mount, whose output is the DOF $v_{2}^{x}$ and the input is the DOF $m_{2}^{x}$.}}
     \label{fig:Identified_RM_SSMs_CF_2}
  \end{subfigure}
  
\caption{\textcolor{black}{Comparison of the FRFs of the identified rubber mount dynamic stiffness by exploiting the following methodologies: dynamic stiffness of the state-space model computed by primal state-space disassembly (approach iv), dynamic stiffness of the minimal-order coupled state-space model obtained by primal state-space disassembly (approach vi), dynamic stiffness obtained by exploiting IS on the inverted state-space model representative of assembly A (approach v) and the dynamic stiffness obtained by exploiting IS on the inverted model of assembly A previously transformed into OCF (approach vii).}}
  \label{fig:Identified_RM_SSMs_CF}
\end{figure}

Additional state-space models were computed by using the primal state-space disassemble with tailored post-processing procedures to compute a minimal-order state-space model representative of the rubber mount and by exploiting IS on the inverted assembly A model transformed into coupling form to compute state-space models transformed into OCF representative of the off diagonal terms of the rubber mount. Figure \ref{fig:Identified_RM_SSMs_CF} reports the comparison between the dynamic stiffness computed by the following procedures: dynamic stiffness of the state-space model computed by exploiting the primal state-space disassembly technique (approach iv), dynamic stiffness of the minimal-order state-space model computed by using primal state-space disassembly to decouple the aluminum crosses models transformed into OCF from the assembly A model transformed into OCF with the post-processing procedure suggested in \cite{RD_2021} that relies on the use of a Boolean localization matrix (approach vi), dynamic stiffness of the model representative of the diagonal apparent mass terms obtained by applying IS approach on the inverted model of assembly A (approach v) and the dynamic stiffness of the model obtained by applying the IS approach on the inverted model of assembly A previously transformed into OCF to identify a model transformed into OCF representative of the diagonal apparent mass terms of the rubber mount by transforming into negative form (see \ref{Negative form of a state-space model representative of apparent mass}) the state-space model representative of its off diagonal terms, whose outputs are placed on cross aluminum A and the inputs are placed on cross aluminum B (approach vii). Once again, to ease the interpretation of figure \ref{fig:Identified_RM_SSMs_CF_2} we have decided to show the correspondent dynamic stiffness matrices of all the computed solutions. 

By observing figure \ref{fig:Identified_RM_SSMs_CF}, it is evident that both dynamic stiffness functions computed by exploiting the primal state-space disassembly are perfectly matching. Hence, we can state that we have validated the use of the primal state-space formulation together with the post-processing procedures outlined in \cite{RD_2021} to compute minimal-order coupled models when dealing with models estimated from measured FRFs. Moreover, the dynamic stiffness obtained by applying IS approch on the untransformed inverted assembly A model and on the same model transformed into OCF are also perfectly matching. This further demonstrates that IS can be also applied on state-space models previously transformed into coupling form.

As final note, we must reflect on the number of states of the state-space models representative of the rubber mount here identified. By following the approaches iii and iv, we ended up with a state-space model representative of the rubber mount composed by 638 states, while by following approach vi, we have identified a state-space model composed by 614 states. Whereas, by using approaches v and vii, the identified state-space models representative of the off diagonal terms are composed by 242 states. The models identified by using the approaches iii, iv and vi are composed by a number of states substantially larger than the models identified by exploiting approaches v or vii due to the performance of decoupling operations to remove the dynamics of crosses aluminum A and B from the model of assembly A. Indeed, as discussed in \cite{MS_2015} and in the introduction of this paper, the performance of decoupling operations promotes the computation of a state-space model, which double includes the dynamics of the decoupled structures. Consequently, this state-space model will be spoiled with pairs of spurious modes, which are difficult to identify and eliminate, when there is no previous knowledge about the structure to be identified (which is the common situation in practice).

Thus, we may state that the LM-SSS decoupling and the primal state-space disassembly hold the advantage of not assuming that the CE to be identified is massless and that there are no cross couplings between DOFs. Consequently, these approaches are valid to identify any structure/CE, whereas the application of IS approach is restricted to identify models representative of structures that are suitable to be characterized by IS (for example CEs, such as rubber mounts). However, when the CE is intended to be experimentally characterized, we are only required to experimentally characterize the assembly were the CE is included, avoiding dismounting and mounting operations and the performance of more experimental tests to characterize the structures to which the CE is coupled. Furthermore, by characterizing the CE by exploiting IS, we are not required to perform decoupling operations, leading to the identification of models composed by a substantially lower number of states. This advantage will be more significant as the number of state-space models to be decoupled in order to identify the intended CE model increases and as these state-space models are composed by more states, thereby as the structures to be decoupled get more complex. As the number of states of the state-space models identified by IS is in general dramatically lower, a significantly lower computational effort is required. Hence, they are interesting, for example, when real-time-domain analyses are intended to be performed with the identified CEs state-space models.

\subsection{Coupling Results}\label{Coupling Results}

\begin{figure}
  \begin{subfigure}[t]{1\textwidth}
    \centering
    \includegraphics[width=.8\textwidth]{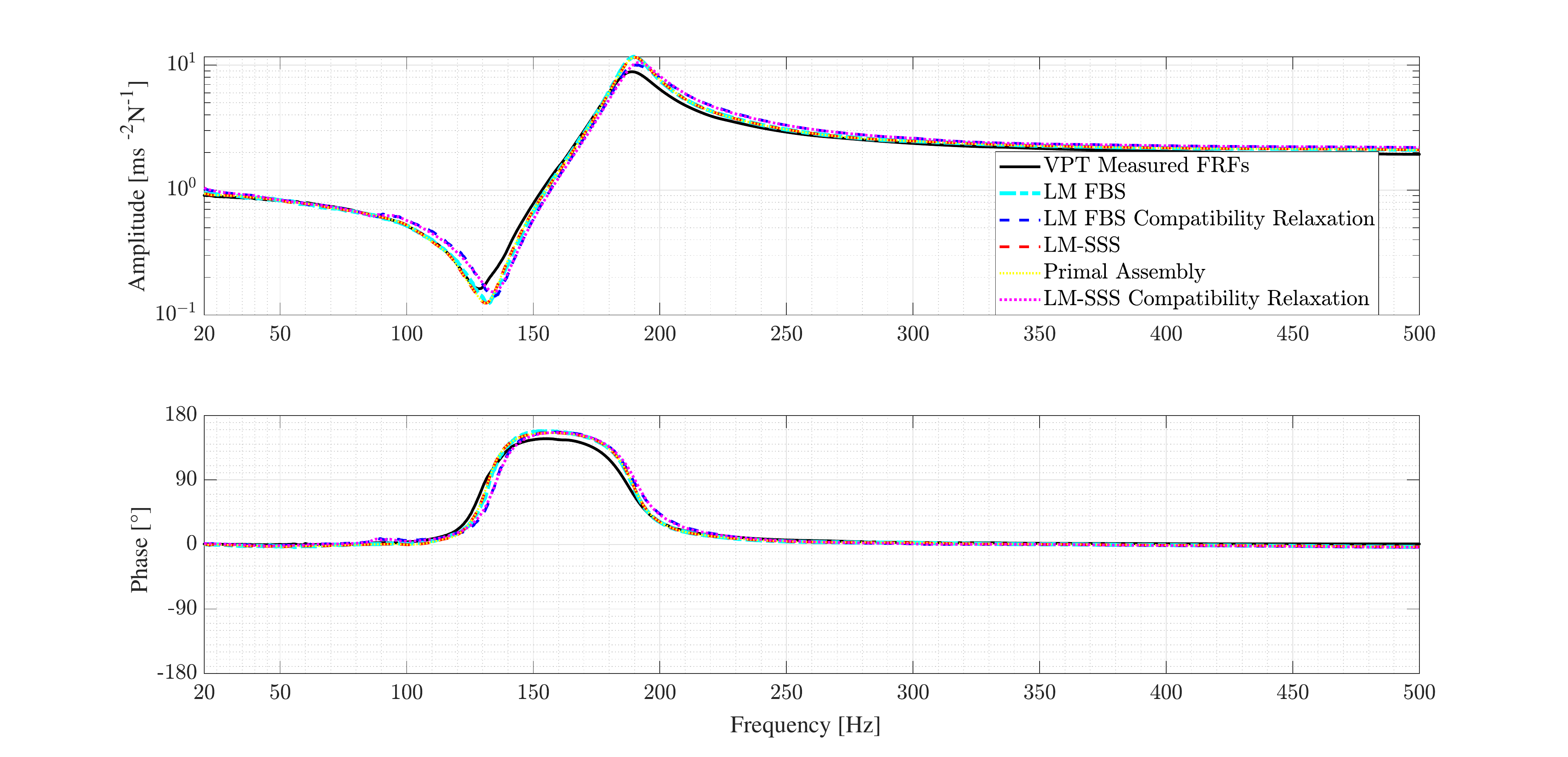}
    \caption{}
     \label{fig:Coupled_Assembly_B_SSMs_1}
  \end{subfigure}\\
  \medskip
  \begin{subfigure}[t]{1\textwidth}
    \centering
    \includegraphics[width=.8\textwidth]{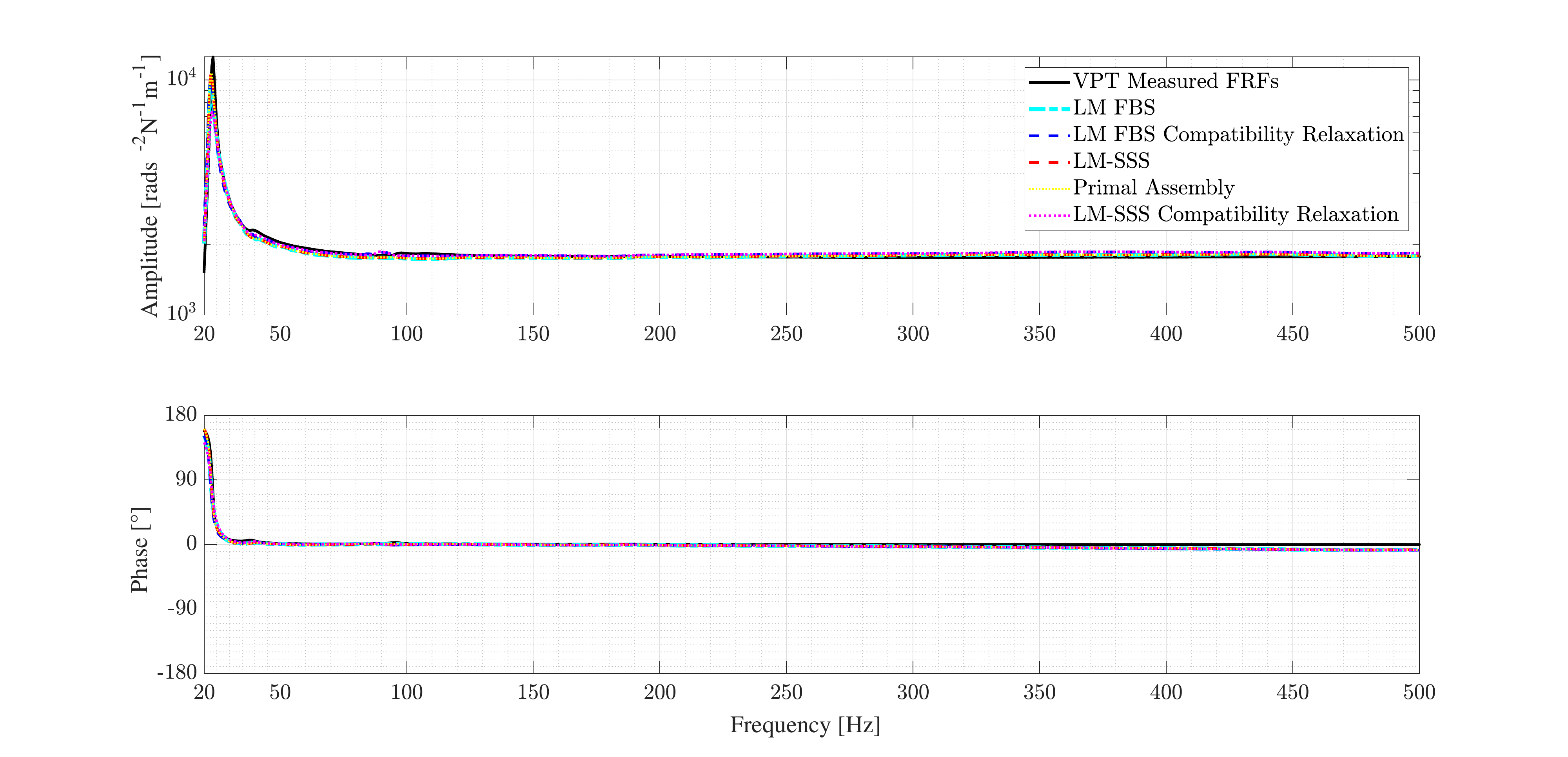}
    \caption{}
     \label{fig:Coupled_Assembly_B_SSMs_2}
  \end{subfigure}
  \caption{\textcolor{black}{Comparison of the interface FRFs of assembly B with the coupled FRFs obtained by using the following methodologies: LM FBS (approach a), LM FBS with compatibility relaxation (approach b), LM SSS (approach c), primal state-space assembly (approach d) and LM-SSS with compatibility relaxation (approach e): a) FRF of the assembly B, whose output is the DOF $v_{2}^{z}$ and the input is the DOF $m_{2}^{z}$; b) FRF of the assembly B, whose output is the DOF $v_{1}^{R_{z}}$ and the input is the DOF $m_{1}^{R_{z}}$.}}
  \label{fig:Coupled_Assembly_B_State_Space_Models}
  \end{figure}

\begin{figure}
    \begin{subfigure}[t]{1\textwidth}
    \centering
    \includegraphics[width=.8\textwidth]{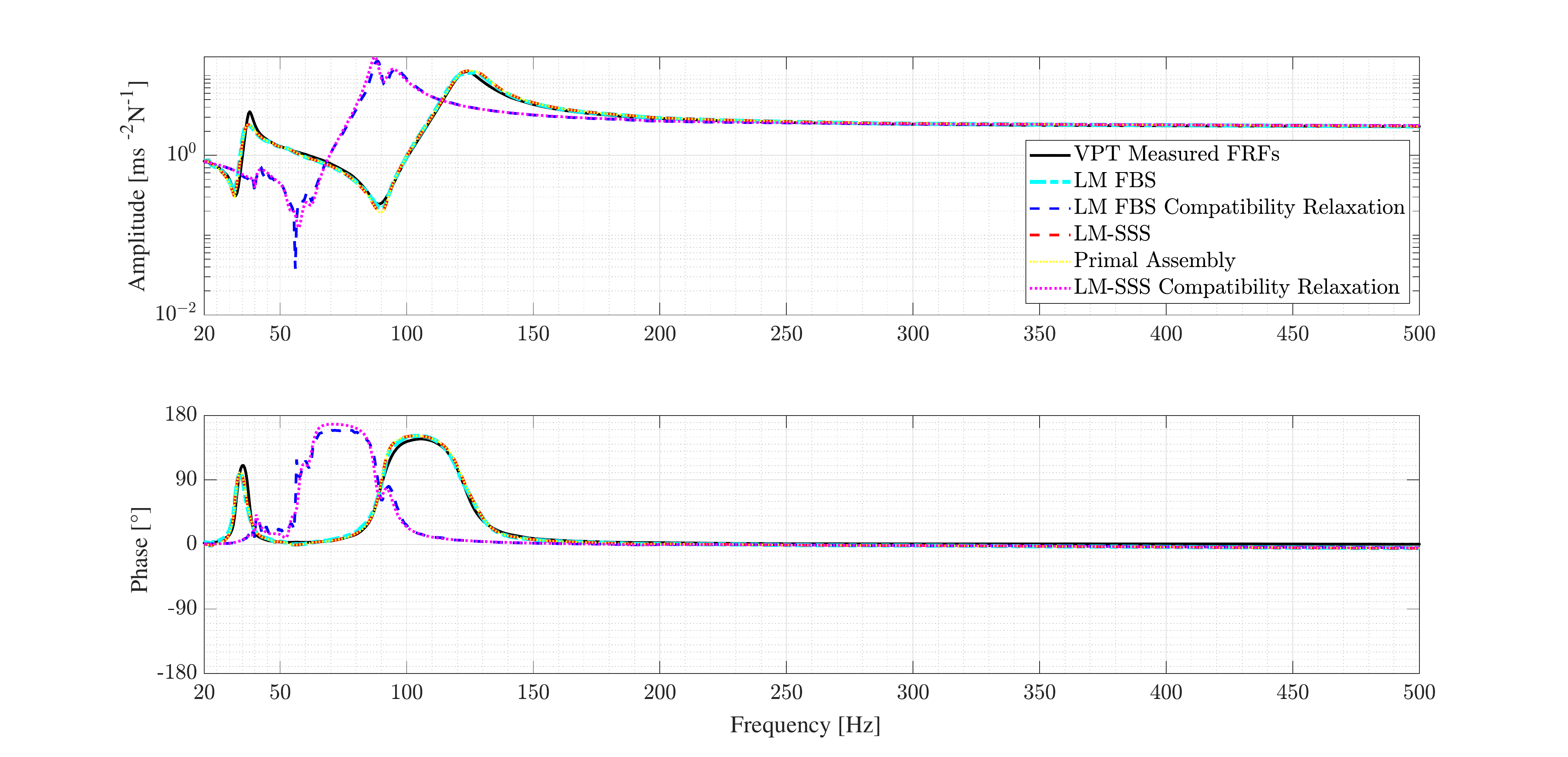}
    \caption{}
    \label{fig:Coupled_Assembly_B_SSMs_3}
 \end{subfigure}
  \hfill
  \begin{subfigure}[t]{1\textwidth}
    \centering
    \includegraphics[width=.8\textwidth]{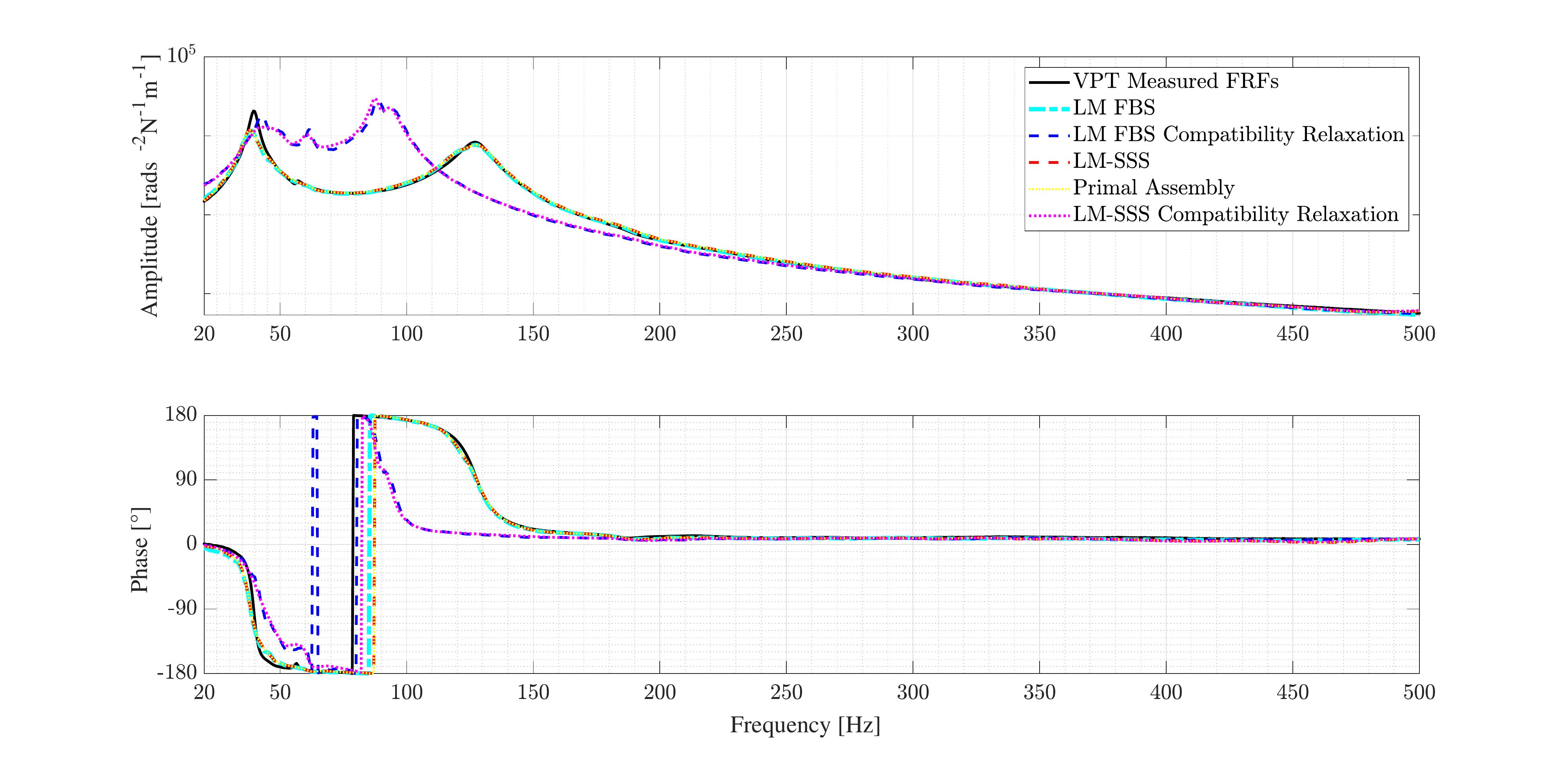}
    \caption{}
     \label{fig:Coupled_Assembly_B_SSMs_4}
  \end{subfigure}
  
\caption{\textcolor{black}{Comparison of the interface FRFs of assembly B with the coupled FRFs obtained by using the following methodologies: LM FBS (approach a), LM FBS with compatibility relaxation (approach b), LM SSS (approach c), primal state-space assembly (approach d) and LM-SSS with compatibility relaxation (approach e): a) FRF of the assembly B, whose output is the DOF $v_{1}^{x}$ and the input is the DOF $m_{1}^{x}$; b) FRF of the assembly B, whose output is the DOF $v_{2}^{R_{x}}$ and the input is the DOF $m_{1}^{R_{x}}$.}}
  \label{fig:Coupled_Assembly_B_State_Space_Models_2}
\end{figure}

In this section, we will start by coupling the state-space models representative of steel crosses A and B with the rubber mount models identified by exploiting LM-SSS decoupling, primal state-space disassembly and IS on the inverted state-space model of assembly A (see section \ref{Identified rubber mount state-space models}). In figures \ref{fig:Coupled_Assembly_B_State_Space_Models} and \ref{fig:Coupled_Assembly_B_State_Space_Models_2} we may observe the comparison of some interface accelerance FRFs of assembly B (obtained by applying VPT on the measured FRFs) with the following solutions: a) the FRFs obtained by coupling with LM-FBS the rubber mount FRFs obtained by using LM-FBS decoupling with the interface FRFs of crosses steel A and B, b) the FRFs obtained by using LM-FBS with compatibility relaxation to couple together the interface FRFs of steel crosses A and B with the inverted diagonal rubber mount terms obtained by exploiting IS on the interface apparent mass of assembly A, c) the FRFs of the coupled state-space model obtained by using LM-SSS to couple the models of crosses steel A and B with the rubber mount model identified with LM-SSS decoupling, d) the FRFs of the inverted coupled state-space model obtained by using primal state-space assembly to couple the inverted models of the steel crosses A and B with the rubber mount model identified by using primal state-space disassembly and e) the FRFs of the coupled state-space model obtained by exploiting \textcolor{black}{LM-SSS with compatibility relaxation} to couple the steel crosses models with the inverted state-space model of the diagonal rubber mount terms obtained by exploiting the state-space realization of IS on the inverted state-space model of assembly A.

By observing figures \ref{fig:Coupled_Assembly_B_State_Space_Models} and \ref{fig:Coupled_Assembly_B_State_Space_Models_2}, we may conclude that the coupled FRFs computed by LM FBS (approach a), the FRFs of the coupled state-space models obtained by exploiting LM-SSS (approach b) and the FRFs of the coupled model computed with primal state-space assembly (approach c) are well-matching the interface FRFs of assembly B. Thus, we may conclude that the coupled state-space models computed by exploiting LM-SSS and primal state-space assembly are accurate. Moreover, as expected, the FRFs of the coupled state-space models obtained by using LM-SSS and primal state-space assembly are perfectly matching, which further validates the primal formulation developed in section \ref{Primal formulation in time domain}.

Figure \ref{fig:Coupled_Assembly_B_SSMs_1} also shows that the axial (z axis direction in figure \ref{fig:VPT_Cross}) interface FRF of assembly B is well-matching the FRF obtained by exploiting \textcolor{black}{LM-SSS with compatibility relaxation} (see figure \ref{fig:Coupled_Assembly_B_SSMs_1}); the same holds for the FRF associated with the rotations around the axial direction (see figure \ref{fig:Coupled_Assembly_B_SSMs_2}). This is very likely due to the light weight of  the rubber mount (27 g) compared to the weight of a single steel cross (560 g). Furthermore, as discussed in \cite{MH_2020}, where the experimental characterization of a cylindrical shape rubber isolator (same shape of the rubber mount under analyze) is performed, there are no cross couplings between the z direction and the rotations around the x and y axis (axis orientation can be observed in figure \ref{fig:VPT_Cross}). The same can be stated about the cross couplings between the x and y directions and the rotation around the z axis. For these reasons, the underlying assumptions of identifying the rubber mount by applying IS have a small effect over the quality of the results associated with both z direction and rotation around z. However, by observing figures \ref{fig:Coupled_Assembly_B_SSMs_3} and \ref{fig:Coupled_Assembly_B_SSMs_4}, we can conclude that the FRFs associated with the x and y directions and with the rotations around x and y axis computed by \textcolor{black}{LM-SSS with compatibility relaxation} are far from being a good estimate of these radial (associated with x or y axis direction in figure \ref{fig:VPT_Cross}) interface FRFs of assembly B. The low quality results for these directions can be explained by the existence of important cross couplings between the displacements on one of the radial directions and the rotations around the other axial direction (see \cite{MH_2020}), which are ignored when the rubber mount is characterized by exploiting IS approach. Comparing the FRFs obtained by exploiting LM FBS with compatibility relaxation with the FRFs of the state-space model computed with \textcolor{black}{LM-SSS with compatibility relaxation}, it is evident that the FRFs obtained by both methodologies are well-matching. The two solutions are not perfectly fitting because the identification of the state-space models from the measured sets of FRFs is not completely perfect. Thus, we may claim that the \textcolor{black}{LM-SSS with compatibility relaxation} method is validated to perform coupling with state-space models estimated from experimentally acquired data. 

\begin{figure}
  \begin{subfigure}[t]{1\textwidth}
    \centering
    \includegraphics[width=.8\textwidth]{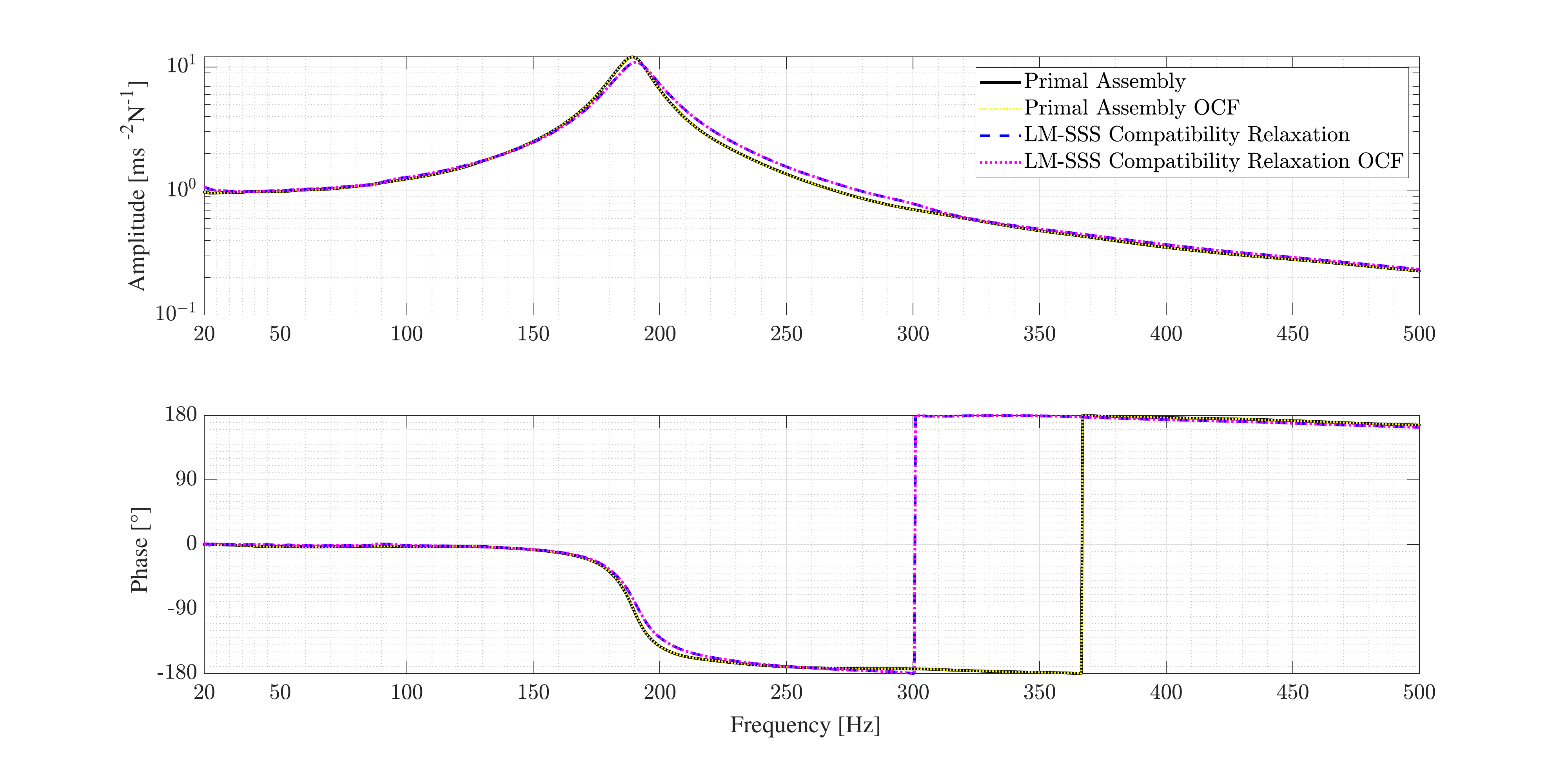}
   \caption{\textcolor{black}{FRF of the assembly B, whose output is the DOF $v_{1}^{z}$ and the input is the DOF $m_{2}^{z}$.}}
     \label{fig:Coupled_Assembly_B_SSMs_CF_1}
  \end{subfigure}\\
  \medskip
  \begin{subfigure}[t]{1\textwidth}
    \centering
    \includegraphics[width=.8\textwidth]{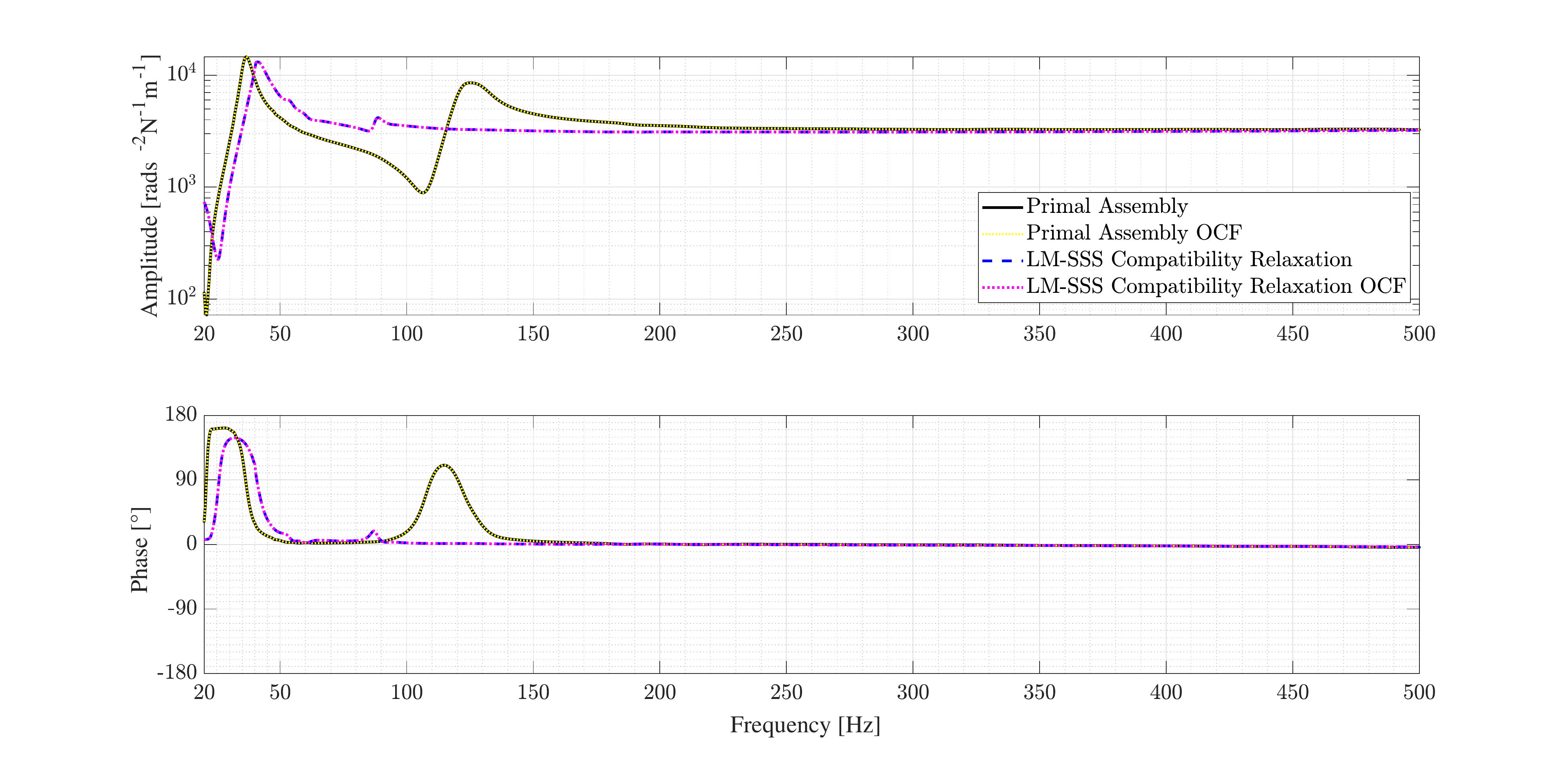}
    \caption{\textcolor{black}{FRF of the assembly B, whose output is the DOF $v_{2}^{R_{y}}$ and the input is the DOF $m_{2}^{R_{y}}$.}}
     \label{fig:Coupled_Assembly_B_SSMs_CF_2}
  \end{subfigure}

\caption{\textcolor{black}{Comparison of the FRFs of the following state-space models: state-space model obtained by primal state-space assembly (approach d), minimal-order coupled state-space model obtained by primal state-space assembly (approach f), state-space model obtained with \textcolor{black}{LM-SSS with compatibility relaxation} (approach e) and the minimal-order coupled state-space model computed with \textcolor{black}{LM-SSS with compatibility relaxation} (approach g)}}
  \label{fig:Coupled_Assembly_B_State_Space_Models_CF}
\end{figure}

At this point, it is worth evaluating the suitability of the primal state-space assembly formulation and of the \textcolor{black}{LM-SSS with compatibility relaxation} technique to compute minimal-order coupled models, when dealing with state-space models computed from experimentally acquired data. Figures \ref{fig:Coupled_Assembly_B_State_Space_Models_CF} compare the FRFs of different coupled models computed by exploiting the following procedures: inverted coupled state-space model computed by exploiting the primal state-space assembly technique (approach d), inverted minimal-order coupled state-space model obtained by using primal state-space assembly to couple the inverted models of the steel crosses A and B previously transformed into OCF with the rubber mount model identified by using primal state-space disassembly with state-space models transformed into OCF together with the post-processing procedure proposed in \cite{RD_2021} that relies on the use of a Boolean localization matrix (approach f), state-space models obtained by using \textcolor{black}{LM-SSS with compatibility relaxation} to couple the models of the steel crosses and of the inverted diagonal apparent mass terms of the rubber mount (approach e) and the minimal-order coupled state-space model computed by using \textcolor{black}{LM-SSS with compatibility relaxation} to couple the steel crosses models previously transformed into OCF with the inverted state-space model of the diagonal apparent mass terms of the rubber mount obtained from the inverted model of assembly A previously transformed into OCF together with the post-processing procedure outlined in section \ref{Minimal-order coupled state-space models} that relies on the use of a Boolean localization matrix (approach g).

By analyzing figures \ref{fig:Coupled_Assembly_B_State_Space_Models_CF}, it is evident that the FRFs of the inverted state-space model obtained by primal state-space assembly are perfectly matching the FRFs of the inverted minimal-order coupled state-space model originated by applying primal state-space assembly with a post-processing procedure to eliminate the redundant states originated from the coupling operation. Thus, we may conclude that the use of primal state-space assembly  together with the post-processing procedures presented in \cite{RD_2021} to compute minimal-order coupled models is valid to deal with state-space models estimated from measured data.

Moreover, further observing figures \ref{fig:Coupled_Assembly_B_State_Space_Models_CF}, we can claim that the FRFs of the state-space model computed by \textcolor{black}{LM-SSS with compatibility relaxation} are perfectly matching the FRFs of the minimal-order coupled state-space model obtained by using \textcolor{black}{LM-SSS with compatibility relaxation} with the post-processing procedure suggested in section \ref{Minimal-order coupled state-space models} that relies on the use of a Boolean localization matrix. Hence, we may state that the use of \textcolor{black}{LM-SSS with compatibility relaxation} together with the post-processing procedures suggested in section \ref{Minimal-order coupled state-space models} is a reliable approach to compute minimal-order coupled models representative of structures connected by CEs suitable to be characterized by exploiting IS. 

As final analysis, we must compare the number of states of the coupled state-space models computed in this section. The coupled models obtained by exploiting approaches c and d, lead to the estimation of models composed by 1018 states, while by following approach f a coupled model with 970 states was obtained. Approaches e and g enabled the computation of coupled state-space models composed by 622 and 610 states, respectively. The coupled models computed by approaches c, d and f are composed by a number of states significantly larger that the models obtained by approaches e and g, because methodologies c, d and f include the dynamics of the rubber mount through state-space models obtained from decoupling operations. In this way, we may conclude that by computing state-space models by including CEs as regular structures there is no need to make assumptions on their physical behavior. Consequently, the method does not suffer on any restriction, neither on frequency ranges nor on the type of CEs characterizing the target structure. However, this approach always lead to the computation of coupled models composed by spurious states and, consequently, composed by a greater number of states.

On the other hand, as \textcolor{black}{LM-SSS with compatibility relaxation} relies on the characterization of CEs by exploiting IS, its applicability is restricted to CEs suitable to be characterized by IS. Furthermore, its application is limited up to frequencies for which the mass of the CEs presents a negligible contribution to their dynamics. Nonetheless, if the mass of the CEs is negligible when compared with the structures to which it is coupled, the massless assumption will not cause a significant degradation of the coupled results (see \cite{MH_2020}). Moreover, by exploiting \textcolor{black}{LM-SSS with compatibility relaxation} we have the possibility of computing coupled state-space models without presenting spurious states, leading in general to the computation of models composed by a significantly lower number of states. Hence, the state-space models computed by \textcolor{black}{LM-SSS with compatibility relaxation} are particularly useful when time-domain analyses are intended to be performed with the coupled models, for example in a real-time context. Therefore, when dealing with CEs that can be characterized by IS or that their characterization by exploiting IS leads to reliable coupling results, there is no advantage of characterizing the CEs by performing decoupling operations. 

\color{black}

\section{Conclusion}\label{Conclusion}

\color{black}
The state-space realization of IS presented in this article showed to be accurate to identify state-space models representative of diagonal apparent mass terms of CEs that are suitable to be characterized by IS. \textcolor{black}{LM-SSS with compatibility relaxation} showed to be an accurate approach to embed the dynamics of this kind of CEs into the LM-SSS formulation. Moreover, the post-processing procedures outlined in section \ref{Minimal-order coupled state-space models} showed to be accurate to eliminate the extra states originated from the performance of coupling operations with LM-SSS via compatibility relaxation. It was also demonstrated that the FRFs of the coupled state-space models obtained by exploiting LM-SSS with compatibility relaxation perfectly match the coupled FRFs computed with LM FBS via compatibility relaxation (see \cite{EB_2014}), provided that the FRFs of the state-space models involved in the coupling operation with LM-SSS via compatibility relaxation match the ones involved on the coupling operation with LM FBS via compatibility relaxation (see section \ref{Coupling by Relaxing the Compatibility Conditions}). 

Discussing the results obtained from the experimental substructuring application exploited in section \ref{Experimental Validation}, it turned out that by exploiting the state-space realization of IS to identify models representative of experimentally characterized CEs, decoupling operations can be avoided. In this way, state-space models free of spurious states and, thus, composed by a significantly lower number of states can be computed. Furthermore, by embedding these identified models in LM-SSS with compatibility relaxation, coupled state-space models composed by a dramatically lower number of states can be computed. We may then claim that the approaches here discussed are interesting to compute state-space models (representative of CEs or of assemblies composed by CEs) tailored to perform time-domain analyses. Nevertheless, it was found that the applicability of both state-space realization of IS and LM-SSS with compatibility relaxation is limited in frequency and restricted to CEs suitable to be characterized by IS, i.e. CEs that do not present important cross couplings between their DOFs and in a frequency range for which their mass presents a negligible contribution to their dynamic behaviour (see sections \ref{Identified rubber mount state-space models} and \ref{Coupling Results}).    

\color{black}

\section*{Credit authorship contribution statement}

\textbf{R.S.O. Dias}: Conceptualization, Investigation, Methodology, Software, Formal analysis, Validation, Data curation, Writing - original draft, Writing - review \& editing. \textbf{M. Martarelli}: Conceptualization, Methodology, Resources, Funding Acquisition, Writing - review \& editing, Supervision, Project administration. \textbf{P. Chiariotti}: Conceptualization, Methodology, Resources, Funding Acquisition, Writing - review \& editing, Supervision, Project administration.

\section*{Declaration of Competing Interest}

The authors declare that they have no known competing financial interests or personal relationships that could have appeared to influence the work reported in this paper.

\section*{Acknowledgements}

The authors gratefully acknowledge Dr. Mahmoud El-Khafafy from Siemens Industry Software NV for supporting the research by providing the system identification algorithm to estimate state-space models from FRFs.

\section*{Funding}

This project has received funding from the European Union's Framework Programme for Research and Innovation Horizon 2020 (2014-2020) under the Marie Sklodowska-Curie Grant Agreement nº 858018.

\setcounter{figure}{0}

\appendix

\section{Inverting state-space models}\label{Inverting_state_space_models}

When a state-space model is inverted, its output vector becomes the input vector of the inverted model and vice-versa. To demonstrate how the inverted state-space model can be computed, let us consider the following generic acceleration state-space model.

\begin{equation}\label{eq:equationfirstformorder}
\begin{gathered}
\{\dot{x}(t)\}=[A]\{x(t)\}+[B]\{u(t)\}\\
\{\ddot{y}(t)\}=[C]\{x(t)\}+[D]\{u(t)\}
\end{gathered}
\end{equation}

By solving the output equations of \eqref{eq:equationfirstformorder} to find the value of $\{u(t)\}$, the equation below is obtained.

\begin{equation}\label{eq:explicitCD12}
\{u(t)\}=[D^{-1}](\{\ddot{y}(t)\}-[C]\{x(t)\})
\end{equation}

By using equation \eqref{eq:explicitCD12}, the state equations of the state-space model given by expression \eqref{eq:equationfirstformorder} can be rewritten as in equation \eqref{eq:equationAB}. 

\begin{equation}\label{eq:equationAB}
\{\dot{x}(t)\}=[A]\{x(t)\}+[B][D^{-1}](\{\ddot{y}(t)\}-[C]\{x(t)\})\\
\end{equation}

By using equations \eqref{eq:equationAB} and \eqref{eq:explicitCD12}, the inverted state-space model can be computed as follows:

\begin{equation}\label{eq:equationfirstformorder3}
\begin{gathered}
\{\dot{x}(t)\}=[A^{inv}]\{x(t)\}+[B^{inv}]\{\ddot{y}(t)\}\\
\{u(t)\}=[C^{inv}]\{x(t)\}+[D^{inv}]\{\ddot{y}(t)\}
\end{gathered}
\end{equation}

where, matrices $[A^{inv}]$, $[B^{inv}]$, $[C^{inv}]$ and $[D^{inv}]$ are given by expression \eqref{eq:inv_matrix_value}.

\begin{alignat}{2}\label{eq:inv_matrix_value}
\begin{split}
[A^{inv}]=[A]-[B][D^{-1}][C],\ \ \ & [B^{inv}]=[B][D^{-1}]\\ 
[C^{inv}]=-[D^{-1}][C],\ \ \ & [D^{inv}]=[D^{-1}]
\end{split}
\end{alignat}

\color{black}

\section{Extending the \textcolor{black}{LM-SSS with compatibility relaxation} formulation to directly compute coupled displacement and velocity state-space models}\label{Extending the LM-SSS-WI formulation to directly compute coupled displacement and velocity state-space models}

In this section, the expressions \eqref{eq:coupled_matrices_value_non_rigid_formulation} that enable the computation of acceleration state-space models by exploiting \textcolor{black}{LM-SSS with compatibility relaxation} (see section \ref{LM_SSS_Non_Rigid_Formulation}) will be adjusted to directly compute displacement and velocity coupled models. As the state and input matrices are equal for displacement, velocity or acceleration state-space models, we are only required to adjust the output and feed-through state-space matrices in expressions \eqref{eq:coupled_matrices_value_non_rigid_formulation}.

To start, let us rewrite the expressions to compute $[\bar{C}^{accel}]$ and $[\bar{D}^{accel}]$ (see equation \eqref{eq:coupled_matrices_value_non_rigid_formulation}) by using displacement state-space matrices (see \citep{FL_1988}).

\begin{equation}\label{eq:coupled_matrices_value_from_disp}
\begin{gathered}
[\bar{C}_{S}^{accel}]=C^{disp}_{S,D}A_{S,D}A_{S,D}-(C^{disp}_{S,D}A_{S,D}B_{S,D})B_{C}^{T}(B_{C}(C^{disp}_{S,D}A_{S,D}B_{S,D})B_{C}^{T}+\\
(C^{disp}_{M,D}A_{M,D}B_{M,D}))^{-1}B_{C}C^{disp}_{S,D}A_{S,D}A_{S,D}\\
[\bar{C}_{M}^{accel}]=(C^{disp}_{S,D}A_{S,D}B_{S,D})B_{C}^{T}(B_{C}(C^{disp}_{S,D}A_{S,D}B_{S,D})B_{C}^{T}+(C^{disp}_{M,D}A_{M,D}B_{M,D}))^{-1}C^{disp}_{M,D}A_{M,D}A_{M,D}\\
[\bar{D}_{S}^{accel}]=C^{disp}_{S,D}A_{S,D}B_{S,D}-(C^{disp}_{S,D}A_{S,D}B_{S,D})B_{C}^{T}(B_{C}(C^{disp}_{S,D}A_{S,D}B_{S,D})B_ {C}^{T}+\\
(C^{disp}_{M,D}A_{M,D}B_{M,D}))^{-1}B_{C}(C^{disp}_{S,D}A_{S,D}B_{S,D})
\end{gathered}
\end{equation}

From the relation between the acceleration and displacement output and feed-through matrices, expressions to compute $[\bar{C}^{disp}]$ and $[\bar{D}^{disp}]$ can be established from expressions \eqref{eq:coupled_matrices_value_from_disp} as in equation \eqref{eq:coupled_matrices_value_from_disp_2}.

\begin{equation}\label{eq:coupled_matrices_value_from_disp_2}
\begin{gathered}
[\bar{C}_{S}^{disp}]=C^{disp}_{S,D}-(C^{disp}_{S,D}A_{S,D}B_{S,D})B_{C}^{T}(B_{C}(C^{disp}_{S,D}A_{S,D}B_{S,D})B_{C}^{T}+(C^{disp}_{M,D}A_{M,D}B_{M,D}))^{-1}B_{C}C^{disp}_{S,D}\\
[\bar{C}_{M}^{disp}]=(C^{disp}_{S,D}A_{S,D}B_{S,D})B_{C}^{T}(B_{C}(C^{disp}_{S,D}A_{S,D}B_{S,D})B_{C}^{T}+(C^{disp}_{M,D}A_{M,D}B_{M,D}))^{-1}C^{disp}_{M,D}\\
[\bar{D}_{S}^{disp}]=C^{disp}_{S,D}A_{S,D}B_{S,D}-(C^{disp}_{S,D}A_{S,D}B_{S,D})B_{C}^{T}(B_{C}(C^{disp}_{S,D}A_{S,D}B_{S,D})B_{C}^{T}+\\
(C^{disp}_{M,D}A_{M,D}B_{M,D}))^{-1}B_{C}(C^{disp}_{S,D}A_{S,D}B_{S,D})-\bar{C}^{disp}\bar{A}\bar{B}=0
\end{gathered}
\end{equation}

From the relation between the acceleration and velocity output and feed-through matrices, expressions to compute $[\bar{C}^{vel}]$ and $[\bar{D}^{vel}]$ can be established from equations \eqref{eq:coupled_matrices_value_from_disp} as follows:

\begin{equation}\label{eq:coupled_matrices_value_from_vel}
\begin{gathered}
[\bar{C}_{S}^{vel}]=C^{disp}_{S,D}A_{S,D}-(C^{disp}_{S,D}A_{S,D}B_{S,D})B_{C}^{T}(B_{C}(C^{disp}_{S,D}A_{S,D}B_{S,D})B_{C}^{T}+(C^{disp}_{M,D}A_{M,D}B_{M,D}))^{-1}B_{C}C^{disp}_{S,D}A_{S,D}\\
[\bar{C}_{M}^{vel}]=(C^{disp}_{S,D}A_{S,D}B_{S,D})B_{C}^{T}(B_{C}(C^{disp}_{S,D}A_{S,D}B_{S,D})B_{C}^{T}+(C^{disp}_{M,D}A_{M,D}B_{M,D}))^{-1}C^{disp}_{M,D}A_{M,D}\\
[\bar{D}_{S}^{vel}]=C^{disp}_{S,D}B_{S,D}-(C^{disp}_{S,D}A_{S,D}B_{S,D})B_{C}^{T}(B_{C}(C^{disp}_{S,D}A_{S,D}B_{S,D})B_{C}^{T}+\\
(C^{disp}_{M,D}A_{M,D}B_{M,D}))^{-1}B_{C}(C^{disp}_{S,D}B_{S,D})
\end{gathered}
\end{equation}

where, superscript $vel$ denotes a state-space matrix of a velocity state-space model (a state-space model whose output vector elements are velocities). If the coupled state-space model obeys the Newton's second law $[C^{disp}_{S,D}B^{disp}_{S,D}]=[0]$ \citep{FL_1988},\citep{SJO_20072697}, hence $[\bar{D}_{S}^{vel}]=[0]$.

As previously mentioned, the state and input matrices are equal for displacement, velocity or acceleration state-space models and are given by equations \eqref{eq:coupled_matrices_value_non_rigid_formulation}. However, those equations can also be easily rewritten by using displacement state-space matrices.

\color{black}

\section{Negative form of a state-space model representative of apparent mass}\label{Negative form of a state-space model representative of apparent mass}

This section demonstrates how to establish the negative form of a state-space model representative of apparent mass. This demonstration can be performed by starting from a numerical acceleration state-space model, which can be established by using the mass, stiffness and damping matrices of the mechanical system under analysis as presented in \citep{FL_1988}.

Let us consider the following acceleration state-space model:

\begin{equation}\label{eq:accel_ss_model}
\begin{gathered}
\{\dot{x}(t)\}=[A]\{x(t)\}+[B]\{u(t)\}\\
\{\ddot{y}(t)\}=[C]\{x(t)\}+[D]\{u(t)\}
\end{gathered}
\end{equation}

where, the state-space matrices are given as follows:

\begin{alignat}{2}\label{eq:MatrixABCD_numerical}
\begin{split}
[A]=\left[
\begin{matrix}
-M^{-1}V & -M^{-1}K\\
I & 0\\
\end{matrix}
\right],\ \ \ & [B]=\left[
\begin{matrix}
M^{-1}\\
0
\end{matrix}
\right]\\ 
[C]=\left[\begin{matrix}
-M^{-1}V & -M^{-1}K
\end{matrix}
\right],\ \ \ & [D]=\left[
\begin{matrix}
M^{-1}
\end{matrix}
\right]
\end{split}
\end{alignat}

and matrices $[M]$, $[K]$ and $[V]$ are the mass, stiffness and damping matrices respectively. 

By following expression \eqref{eq:inv_matrix_value}, the inversion of the state-space model given by expression \eqref{eq:accel_ss_model} can be computed. The state-space matrices of the resultant model are given as follows:

\begin{equation}\label{inverted_ss_matrices}
\begin{gathered}
[A^{inv}]=[A]-[B][D^{-1}][C]=\left[
\begin{matrix}
-M^{-1}V & -M^{-1}K\\
I & 0\\
\end{matrix}
\right]-\left[
\begin{matrix}
M^{-1}\\
0\\
\end{matrix}
\right]\left[
\begin{matrix}
M^{-1}
\end{matrix}
\right]^{-1}\left[
\begin{matrix}
-M^{-1}V & -M^{-1}K
\end{matrix}
\right]\\
=\left[
\begin{matrix}
-M^{-1}V & -M^{-1}K\\
I & 0\\
\end{matrix}
\right]-\left[
\begin{matrix}
-M^{-1}V & -M^{-1}K\\
0 & 0
\end{matrix}
\right]=\left[
\begin{matrix}
0 & 0\\
I & 0\\
\end{matrix}
\right]\\
\\
[B^{inv}]=[B][D^{-1}]=\left[
\begin{matrix}
M^{-1}\\
0\\
\end{matrix}
\right]\left[
\begin{matrix}
M^{-1}
\end{matrix}
\right]^{-1}=\left[
\begin{matrix}
I\\
0\\
\end{matrix}
\right]\\
\\
[C^{inv}]=-[D^{-1}][C]=-\left[
\begin{matrix}
M^{-1}
\end{matrix}
\right]^{-1}\left[
\begin{matrix}
-M^{-1}V & -M^{-1}K
\end{matrix}
\right]=\left[
\begin{matrix}
V & K
\end{matrix}
\right]\\
\\
[D^{inv}]=[D^{-1}]=\left[
\begin{matrix}
M^{-1}
\end{matrix}
\right]^{-1}=[M]
\end{gathered}
\end{equation}

As presented in \citep{MS_2015},\citep{MEL_2019325}, the computation of the negative form of a state-space model can be performed by multiplying its $[M]$, $[K]$ and $[V]$ matrices by $-1$. Hence, the negative form of the state-space matrices presented in the set of equations \eqref{inverted_ss_matrices} can be obtained as shown below:

\begin{align}\label{eq:negativestatespacemodel}
\begin{split}
[A^{inv}_{neg}]=\left[
\begin{matrix}
0 & 0\\
I & 0
\end{matrix}
\right],\ \ \  & [B_{neg}^{inv}]=\left[
\begin{matrix}
I\\
0
\end{matrix}
\right]\\
[C_{neg}^{inv}]=[
\begin{matrix}
-V & -K
\end{matrix}
],\ \ \ &  [D_{neg}^{inv}]=\left[
\begin{matrix}
-M
\end{matrix}
\right]
\end{split}
\end{align}

where, subscript $neg$ denotes a state-space matrix transformed into negative form.

By analyzing expression \eqref{eq:negativestatespacemodel}, we may conclude that to obtain the negative form of a state-space model representative of apparent mass its output and feed-through matrices must be multiplied by $-1$. This procedure is also valid to compute the negative form of state-space models established in modal domain.

 \bibliographystyle{elsarticle-num} 
 \bibliography{cas-refs}





\end{document}